\documentclass[11pt, a4paper]{article}
\usepackage[english]{babel}
\usepackage[margin=.7in]{geometry}

\usepackage{amsthm, ,amssymb}
\theoremstyle{plain} 
\newtheorem{proposition}{Proposition}[section] 
\newtheorem{lemma}{Lemma}[section]
\theoremstyle{definition} 

\usepackage{natbib}
\usepackage{ragged2e, color}
\usepackage{multicol}
\usepackage{colortbl}
\usepackage{multirow} 
\usepackage{pifont}	
\usepackage{adjustbox}
\usepackage{caption}
\usepackage{subcaption}
\usepackage{booktabs}
\usepackage{amssymb}
\usepackage{amsmath}
\usepackage{mathtools}
\usepackage{dsfont}
\usepackage{bm,bbm}
\usepackage{rotating}
\usepackage{lscape}
\usepackage{adjustbox}
\usepackage[skins,theorems]{tcolorbox}
\tcbset{highlight math style={enhanced,
		colframe=blue,colback=white,arc=0pt,boxrule=1pt}}

\usepackage[customcolors]{hf-tikz}
\usepackage{setspace}
\usepackage{tikz} 
\usetikzlibrary{arrows,shapes,fit,backgrounds}
\usetikzlibrary{arrows.meta}
\usetikzlibrary{shapes.geometric}
\usetikzlibrary{fit}

\newcolumntype{a}{>{\columncolor{red!30}}p{3.3cm}}
\makeatletter
\def\blfootnote{\xdef\@thefnmark{}\@footnotetext}
\makeatother
\definecolor{green}{rgb}{0.0, 0.26, 0.15}
\definecolor{azure}{rgb}{0.0, 0.5, 1.0}
\definecolor{transparent}{rgb}{1,.8,.3}
\definecolor{cambridgeblue}{rgb}{0.64, 0.76, 0.68}
\definecolor{ashgrey}{rgb}{0.7, 0.75, 0.71}
\definecolor{beaublue}{rgb}{0.74, 0.83, 0.9}

\usepackage[skins,theorems]{tcolorbox}
\tcbset{highlight math style={enhanced,
		colframe=azure,colback=white,arc=0pt,boxrule=1pt}}
\usepackage{xcolor}

\newtheorem*{definition*}{Definition}
\usepackage[ruled]{algorithm}
\usepackage[]{algorithmic}
\usepackage[makeroom]{cancel}
\renewcommand{\algorithmicrequire}{\textbf{Input:}}
\renewcommand{\algorithmicensure}{\textbf{Output:}}
\floatname{algorithm}{Algorithm}
\algsetup{indent=2em}
\newcommand{\STATEnonum}{\item[]}
\newcommand{\bunderline}[1]{\underline{#1\mkern-2mu}\mkern2mu }
\newcounter{storesubequations}

\newcommand{\KPDP}{KPDP}
\newcommand{\KPDPs}{KPDPs}

\newcommand{\cycleCap}{K}
\newcommand{\cycleSet}{\mathcal{C}_{\cycleCap}}
\newcommand{\chainCap}{L}
\newcommand{\chainSet}{\mathcal{C}_{\chainCap}}
\newcommand{\cyclechainSet}{\mathcal{C}}
\newcommand{\pairSet}{P}
\newcommand{\nddSet}{N}
\newcommand{\vertexSet}{V}
\newcommand{\arcSet}{A}
\newcommand{\graph}{D}
\newcommand{\graphKEP}{\graph = (\vertexSet, \arcSet)}
\newcommand{\containedVxt}{V(\cdot)}
\newcommand{\containedArc}{A(\cdot)}
\newcommand{\match}{c}
\newcommand{\matchSet}{M}

\newcommand{\matchSpace}{\matchSet_D}
\newcommand{\varFirstc}{\tilde{\varFirst}}
\newcommand{\cycleSetSndS}{\cycleSet^{\policy}(\varFirst, \oneCe)}
\newcommand{\chainSetSndS}{\chainSet^{\policy}(\varFirst, \oneCe)}
\newcommand{\cycleSetSndTr}{\cycleSet^{\policy}(\varFirst)}
\newcommand{\chainSetSndTr}{\chainSet^{\policy}(\varFirst)}
\newcommand{\cycleSetSndTrc}{\cycleSet^{\policy}(\varFirstc)}
\newcommand{\chainSetSndTrc}{\chainSet^{\policy}(\varFirstc)}
\newcommand{\digraphSndS}{D^{\policy}(\varFirst, \oneCe)}

\newcommand{\digraphSndSc}{D^{\policy}(\varFirstc, \cdteOne)}
\newcommand{\digraphSndSpr}{D^{\policy}(\varFirstc, \oneCeNoPrime)}

\newcommand{\digraphSndTr}{D^{\policy}(\varFirst)}
\newcommand{\digraphSndTrc}{D^{\policy}(\varFirstc)}
\newcommand{\vertexSetSndBe}{V^{\policy, \varFirst}}\newcommand{\vertexSetSndBeF}{V^{\oneCe}}
\newcommand{\arcSetSndBe}{A^{\policy, \varFirst}}
\newcommand{\arcSetSndBeF}{A^{\oneCe}}

\newcommand{\vertexa}{u}
\newcommand{\vertexb}{v}

\newcommand{\parc}{(\vertexa,\vertexb)}
\newcommand{\parcrev}{(\vertexb,\vertexa)}

\newcommand{\narc}{\vertexa\vertexb}
\newcommand{\narcrev}{\vertexb\vertexa}

\newcommand{\oneCeS}{\gamma}
\newcommand{\oneCe}{\bm{\oneCeS}}
\newcommand{\HatoneCe}{\bar{\bm{\gamma}}}
\newcommand{\oneCeStar}{\bm{\gamma}^{\star}}
\newcommand{\bigsCe}{\Gamma}

\newcommand{\vertex}{u}
\newcommand{\smallsCeVx}{\gamma^{\text{v}}}
\newcommand{\smallsCeArc}{\gamma^{\text{a}}}
\newcommand{\tildesmallsCeVx}{\tilde{\gamma}^{\text{v}}}
\newcommand{\tildesmallsCeArc}{\tilde{\gamma}^{\text{a}}}
\newcommand{\UbigsCe}{\bigsCe^{\text{v}}}
\newcommand{\AbigsCe}{\bigsCe^{\text{a}}}
\newcommand{\vBudget}{r^{\text{v}}}
\newcommand{\aBudget}{r^{\text{a}}}
\newcommand{\charVecFirst}{\bm{\varFirstSel_{\matchSet}}}
\newcommand{\charVecSecond}{\bm{\varRecop_{\matchSet}}}

\newcommand{\SelVerXc}{V(\varFirstc)}
\newcommand{\SelVerX}{V(\varFirst)}
\newcommand{\SelPairsX}{\pairSet(\varFirstc)}

\newcommand{\cycleSetSndSc}{\cycleSet^{\policy}(\varFirstc, \cdteOne)}
\newcommand{\chainSetSndSc}{\chainSet^{\policy}(\varFirstc, \cdteOne)}
\newcommand{\cychSetSndSc}{\mathcal{C}_{K,L}^{\policy}(\varFirstc, \cdteOne)}
\newcommand{\cychSetSndScdte}{\mathcal{C}_{K,L}^{\policy}(\varFirstc, \cdteOne)}
\newcommand{\cychSetSndS}{\mathcal{C}_{K,L}^{\policy}(\varFirstc, \oneCe)}

\newcommand{\policy}{\pi}
\newcommand{\policySet}{\Pi}
\newcommand{\varFirstSel}{x}
\newcommand{\varFirstPlain}{\bm{\varFirstSel}}
\newcommand{\varFirst}{\varFirstPlain}
\newcommand{\setFirst}{\mathcal{X}}

\newcommand{\varReco}{\mathbf{y}}

\newcommand{\varRecop}{y}

\newcommand{\setRecoSym}{\mathcal{Y}}
\newcommand{\setReco}{\setRecoSym^{\policy}(\varFirst, \oneCe)}
\newcommand{\setRecoTrcdte}{\setRecoSym_{\text{tr}}^{\policy}(\varFirstc)}

\newcommand{\setRecoExpc}{\setRecoSym_{\text{tr}}^{\policy}(\varFirstc, \cdteOne)}

\newcommand{\setRecocdte}{\setRecoSym^{\policy}(\varFirstc, \cdteOne)}
\newcommand{\setRecoc}{\setRecoSym^{\policy}(\varFirstc, \oneCe)}

\newcommand{\weightRO}{{\mathbf{w}^{\oneCe}(\varFirst)}^{\top}}
\newcommand{\weightROc}{{\mathbf{w}^{\oneCe}(\varFirstc)}^{\top}}

\newcommand{\objRO}{\weightRO\varReco}
\newcommand{\objROc}{\weightROc\varReco}

\newcommand{\SingleRO}{\bm{P}(\policy, \bigsCe)}
\newcommand{\SingleROr}{\bm{P}(\policy, \tilde{\bigsCe})}
\newcommand{\objSingleRO}{Z_{P}}
\newcommand{\objSingleROr}{\tilde{Z}_{P}}
\newcommand{\iteTwoStage}{i}

\newcommand{\SSPref}{\bm{Q}(\policy, \varFirstc)}
\newcommand{\objTwoStage}{Z^{\pi}_{Q}(\varFirstc)}

\newcommand{\objTwoStageStar}{Z_{Q}^{\pi,\star}(\varFirstc)}

\newcommand{\polySetxOne}{\matchSet^{\policy}(\varFirst, \oneCe)}

\newcommand{\varRecoZ}{z}

\newcommand{\wTwoStageGralc}{\mathbf{w_{\match}}(\varFirstc)}

\newcommand{\wTwoStagec}{\mathbf{w_{\match}}(\varFirstc)}

\newcommand{\wTwoStageExpanded}{\mathbf{\hat{w}_{\match}}(\varFirstc,\cdteOne)}

\newcommand{\RecoP}{R^{\policy}(\varFirst, \oneCe)}
\newcommand{\Rexp}{R^{\policy}_{\text{exp}}}

\newcommand{\RecoPexpc}{\Rexp(\varFirstc, \cdteOne)}

\newcommand{\varRecoc}{\tilde{\varReco}}

\newcommand{\varCov}{Z}
\newcommand{\UBCovering}{{\bar{\varCov}}^{\pi}_{Q}(\varFirstc)}
\newcommand{\LBCovering}{{\bunderline{\varCov}}^{\pi}_{Q}(\varFirstc)}
\newcommand{\RecoObjVarcdte}{{Z}^{\pi,\star}_{R}(\varFirstc, \cdteOne)}

\newcommand{\RecoExpObjVar}{{Z}^{\pi,\star}_{RE}(\varFirstc, \cdteOne)}

\newcommand{\varRecoOpSndS}{\varReco^{\star\oneCe}}
\newcommand{\varRecoOpSndSc}{\varReco^{\star\cdteOne}}
\newcommand{\VarRecoOpExp}{\psi^{\star\cdteOne}}

\newcommand{\VarRecoOpExpcdte}{\underline{\psi}^{\star\cdteOne}}
\newcommand{\VarRecoOpExpcclt}{\underline{\psi}^{\star\oneCe}}

\newcommand{\incumCe}{\oneCe^{\prime}}
\newcommand{\optCe}{\oneCe^{\star}}

\newcommand{\rhsfeasite}{H^{\cdteOne}}

\newcommand{\cdteOne}{\tilde{\oneCe}}

\newcommand{\RecoPcdte}{R^{\policy}(\varFirstc, \cdteOne)}
\newcommand{\bigsCecR}{\hat{\bigsCe}(\varFirstc)}
\newcommand{\bigsCecdte}{\bigsCe(\varFirst)}
\newcommand{\bigsCec}{\bigsCe(\varFirstc)}

\newcommand{\oneCeNoPrime}{\tilde{\oneCe}}
\newcommand{\oneCePrime}{\tilde{\oneCe}^{\prime}}

\newcommand{\RecoObjVariteEx}{\varCov^{\policy,\star}_{R}(\varFirstc, \oneCeNoPrime)}
\newcommand{\RecoObjVariteExPrime}{\varCov^{\policy,\star}_{R}(\varFirstc,\oneCePrime)}

\newcommand{\cyclechainSetNoPX}{C^{\policy}_{\cycleCap, \chainCap}(\varFirstc, \oneCeNoPrime)}
\newcommand{\cyclechainSetPrimeX}{C^{\policy}_{\cycleCap, \chainCap}(\varFirstc, \oneCePrime)}

\newcommand{\cyclechainSetXvtx}{C^{\policy, \bar{\vertexb}}_{\cycleCap, \chainCap}(\varFirstc)}
\newcommand{\cyclechainSetXarc}{C^{\policy,\bar{a}}_{\cycleCap, \chainCap}(\varFirstc)}

\newcommand{\cover}{\textit{cover}}
\newcommand{\RecoSet}{\mathcal{M}}
\newcommand{\RecoSetexp}{\mathcal{M}_{\text{exp}}}
\newcommand{\HeuAnsTrue}{(\incumCe, \text{\textbf{true}})}
\newcommand{\HeuAns}{(\incumCe, \cover)}

\newcommand{\ZRecoCol}{Z^{\star}_{\text{cg}}}
\newcommand{\ZRecoColUB}{Z^{\text{UB}}_{\text{cg}}}
\newcommand{\TupleZRecoCol}{(\ZRecoCol, \ZRecoColUB, \tilde{\varReco})}

\newcommand{\ZRecoColexp}{Z^{\star}_{\text{Xcg}}}

\newcommand{\ZRecoColUBexp}{Z^{\text{UB}}_{\text{Xcg}}}

\newcommand{\TupleZRecoColexp}{(\ZRecoColexp, \ZRecoColUBexp, \tilde{\varReco})}
\newcommand{\ElSet}{E}
\newcommand{\Elem}{\text{e}}
\newcommand{\Elemstar}{\Elem^\star}
\newcommand{\genmatch}{W}
\newcommand{\Elsta}{s}
\newcommand{\Elwei}{\text{w}}
\newcommand{\matchHeuSet}{\mathcal{E}}
\newcommand{\deltaOutEl}{\delta^{+}(\Elemstar_{n})}
\newcommand{\deltaInEl}{\delta^{-}(\Elemstar_{n})}

\newcommand{\vertxCevar}{t^{\oneCe}_{\vertexb}}
\newcommand{\chainvar}{\delta_{\vertexa\vertexb\ell}}
\newcommand{\chainvarNext}{\delta_{\vertexb\vertexa(\ell + 1)}}
\newcommand{\chainvarNextCe}{\delta^{\oneCe}_{\vertexb\vertexa(\ell + 1)}}
\newcommand{\chainvarCe}{\delta^{\oneCe}_{\vertexa\vertexb\ell}}
\newcommand{\chainvarR}{\delta_{\vertexb\vertexa 1}}
\newcommand{\chainvarRCe}{\delta^{\oneCe}_{\vertexb\vertexa 1}}
\newcommand{\idposset}{\mathcal{L}}
\newcommand{\posset}{\mathcal{L}(\vertexa, \vertexb)}
\newcommand{\cyclesetuv}{\cycleSet^{\vertexa\vertexb}}
\newcommand{\cyclesetv}{\cycleSet^{\vertexb}}

\newcommand{\objPosIdx}{Z}
\newcommand{\arcSetL}{\arcSet_{\ell}}
\newcommand{\cyclevar}{\varRecoZ_\match}
\newcommand{\cyclevarCe}{\varRecoZ_\match^{\oneCe}}

\usepackage{amsmath}
\usepackage{tipx}
\usepackage{graphicx}
\usepackage[colorlinks=true, allcolors=blue]{hyperref}
\usepackage[normalem]{ulem}
\title{A Feasibility-Seeking Approach to Two-stage Robust Optimization in Kidney Exchange}
\author{Lizeth Carolina Riascos-\'{A}lvarez\textsuperscript{\dag}\footnote{Corresponding Author}\\
\texttt{carolina.riascos@mail.utoronto.ca}\\
Merve Bodur\textsuperscript{\dag}\\
\texttt{bodur@mie.utoronto.ca}\\
Dionne M. Aleman\textsuperscript{\dag\ddag}\\ \texttt{aleman@mie.utoronto.ca}\\
\textit{\textsuperscript{\dag}Department of Mechanical \& Industrial Engineering}\\
\textit{\textsuperscript{\ddag}Institute for Health Policy, Management \& Evaluation}\\
\textit{University of Toronto}\\
\textit{5 King’s College Road},
\textit{Toronto, Ontario},
\textit{M5S 3G8, Canada}
}

\begin{document}
\maketitle

\begin{abstract}
Kidney paired donation programs (\KPDPs) match patients with willing but incompatible donors to compatible donors with an assurance that when they donate, their intended recipient receives a kidney in return from a different donor. A patient and donor join a \KPDP\ as a pair, represented as a vertex in a compatibility graph, where arcs represent compatible kidneys flowing from a donor in one pair to a patient in another. A challenge faced in real-world \KPDPs\ is the possibility of a planned match being cancelled, e.g., due to late detection of organ incompatibility or patient-donor dropout. We therefore develop a two-stage robust optimization approach to the kidney exchange problem wherein (1) the first stage determines a kidney matching solution according to the original compatibility graph, and then (2) the second stage repairs the solution after observing transplant cancellations. In addition to considering homogeneous failure, we present the first approach that considers non-homogeneous failure between vertices and arcs. To this end, we develop solution algorithms with a feasibility-seeking master problem and evaluate two types of recourse policies. Our framework outperforms the state-of-the-art kidney exchange algorithm under homogeneous failure on publicly available instances. Moreover, we provide insights on the scalability of our solution algorithms under non-homogeneous failure for two recourse policies and analyze their impact on highly-sensitized patients, patients for whom few kidney donors are available and whose associated exchanges tend to fail at a higher rate than non-sensitized patients.
\end{abstract}

\section{Introduction}
\label{Introduction}

Kidney paired donation programs (\KPDPs) across the globe have facilitated living donor kidney transplants for patients experiencing renal failure who have a willing but incompatible or sub-optimal donor. A patient in need of a kidney transplant registers for a \KPDP\ with their incompatible donor (paired donor) as a pair. The patient then receives a compatible kidney from either the paired donor in another pair, who in turn receives a kidney from another donor, or from a singleton donor who does not expect a kidney in return (called a \emph{non-directed donor}). The transplants are then made possible through the exchange of paired donors between patient-donor pairs. Since first discussed by \cite{Rapaport1986}, and put in practice for the first time in South Korea \citep{Park1999}, kidney exchanges performed through \KPDPs\ have been introduced in several countries around the world, e.g., the United States \citep{Saidman2006}, the United Kingdom \citep{Manlove2015}, Canada \citep{Malik2014} and Australia \citep{AustraliaKPD}, and the underlying matching of patients to donors has been the subject of study from multiple disciplines (see, e.g., \citet{Roth2005,Dickerson2016, Dickerson2019,Carvalho2021, Riascos2020}).

Despite this attention, \KPDPs\ still face challenges from both a practical and theoretical point of view \citep{Ashlagi2021}.
In this work, motivated by the high rate of exchanges that do not proceed to transplant \citep{Bray2015, Dickerson2016, CBS2019}, we provide a robust optimization framework that proposes a set of exchanges for transplant, observes failures in the transplant plan, and then repairs the affected exchanges provided that a recovery rule (i.e., \textit{recourse policy}) is given by \KPDP\ operators. 

Kidney exchange, or the kidney exchange problem as it is known in the literature, can be modeled with a compatibility graph, (i.e., a digraph). Each vertex represents either a pair or a non-directed donor, and each arc indicates that the donor in the starting vertex is blood type and tissue type compatible with the patient in the ending vertex. Arcs may have an associated weight that represents the priority of that transplant. Exchanges then take the form of simple cycles and simple paths (called chains).

Cycles consist of patient-donor pair exchanges only, whereas in a chain the first donation occurs from a non-directed donor to a patient-donor pair, and is followed by a sequence of donations from a paired donor to the patient in the next patient-donor pair. Since a donor/patient can withdraw from a \KPDP\ at any time, or late-detected medical issues can prevent a paired donor from donating in the future, cyclic transplants are performed simultaneously. Unlike cycles, a patient in a chain does not risk exchanging his/her paired donor without first receiving a kidney from another donor. Therefore, transplants in a chain can be executed sequentially, but depending on \KPDP\ regulations, they can also be performed simultaneously \citep{Biro2021}. Due to the simultaneity constraint for cycles, the maximum number of transplants in a cycle is limited by the logistical implications of arranging operating rooms and surgical teams. The maximum number of transplants in a chain can vary depending on the simultaneity constraint and specific regulations. In the United States, the size of chains can be unbounded and theoretically limited only by the number of available donors \citep{Ashlagi2012, Dickerson2012b, Anderson2015,Ding2018, Dickerson2019} by allowing the donor in the last pair of a chain to become a \textit{bridge donor}, i.e., a paired donor that acts as a non-directed donor in a future algorithmic matching. However, in Canada \citep{Malik2014} and in Europe \citep{Biro2021, Carvalho2021}, chains are limited in size since their \KPDPs\ do not use bridge donors. In this case, the donor in the last pair donates to a patient on the deceased donor waiting list, ``closing up'' a chain. 

In constructing kidney exchanges, \KPDP\ operators use ``believed'' information on the compatibility between patients and donors; however, this information may not be accurate as final confirmation of compatibility is confirmed only when a set of transplants has been selected. Thus, there is uncertainty about the existence of vertices and arcs in the compatibility graph. There are multiple reasons a patient-donor pair or a non-directed donor selected for transplant may not be available, e.g., match offer rejection, already transplanted out of the \KPDP, illness, pregnancy, reneging, etc. Thus, even if some of this information is captured prior to the matching process, there is still a chance of subsequent failure. Additionally, tissue type compatibility is not known with certainty when the compatibility graph is built, unlike blood type compatibility. Tissue type compatibility is based on the result of a \textit{virtual crossmatch test}, which typically has lower accuracy than a \textit{physical crossmatch test}. Both tests try to determine if there are antibodies of importance in a patient that can lead to the rejection of a potential donor's kidney. Physical crossmatch tests are challenging to perform, making them impractical, and thus unlikely to be performed between all patients and donors in real life \citep{Carvalho2021}. So, in the first stage, a ``believed'' compatibility graph is built according to the results of the virtual crossmatch test. Once the set of transplants has been proposed, those exchanges undergo a physical crossmatch test to confirm or rule out the viability of the transplants. After confirming infeasible transplants and depending on the \KPDPs' regulations, \KPDP\ operators may attempt to repair the originally planned cycles and chains impacted by the non-existence of a vertex or an arc. We refer to these impacted cycles and chains as \textit{failed} cycles and chains. A cycle fails completely if any of its elements (vertices or arcs) cease to exist, whereas a chain is cut short at the pair preceding the first failed transplant.

Failures in the graph, caused by the disappearance of a vertex or arc, have a significant impact on the number of exchanges that actually proceed to transplant \citep{Dickerson2019,CBS2019} and can even drive \KPDP\ regulations \citep{Carvalho2021}. For instance, \cite{Dickerson2019} reported that for selected transplants from the UNOS program between 2010--2012, 93\% did not proceed to transplant. Of those non-successful transplants, 44\% had a failure reason (e.g., failed physical crossmatch result or patient/donor drop-out), which in turn caused the cancellation of the other 49\% of non-successful transplants. In Canada, between 2009--2018, 62\% of cycles and chains with six transplants failed, among which only 10\% could be repaired \citep{CBS2019}. Half the cycles and chains with three or fewer transplants were successful, and approximately 30\% of the total could not proceed to transplant. 

The set of transplants that is proposed and later repaired by \KPDPs\ does not account for subsequent failures, and thus the number of successful transplants performed is sub-optimal. There is a large body of work concerned with maximizing the expectation of proposed transplants \citep{Awasti2009, Dickerson2014, Dickerson2019, Klimentova2016, Smeulders2022}. From a practical perspective, maximizing expectation could increase the number of planned exchanges that become actual transplants. However, such a policy risks favoring patients that are likely to be compatible with multiple donors over highly-sensitized patients \citep{Carvalho2021}, i.e., patients for whom few kidney donors are available and whose associated exchanges tend to fail at a higher rate than non-sensitized patients. Furthermore, real data is limited and there is not currently a enough understanding of the dynamics between patients and donors to derive a probability distribution of failures that could be generalized to most \KPDPs. We therefore model failure through an uncertainty set that does not target patient sensitization level or probabilistic knowledge, and aims to find a set of transplants that allows the biggest recovery under the worst-case failure scenario in the uncertainty set.

We develop a two-stage robust optimization (RO) approach to the kidney exchange problem wherein (1) the first stage determines a kidney matching solution according to the original compatibility graph, and then (2) the second stage repairs the solution after observing transplant cancellations. We extend the current state-of-the-art RO methodologies \citep{Carvalho2021} by assuming that the failure rate of vertices and arcs is non-homogeneous, since failure reasons, such as a late-detected incompatibility, seem to be independent from patient/donor dropout and vice versa. Moreover, we also consider the impact of scarce match possibilities for highly-sensitized patients. 

The contributions of this work are as follows:
\begin{enumerate}
    \item We develop a novel two-stage RO framework with non-homogeneous failure between vertices and arcs.
    \item We present a novel general solution framework for any recourse policy whose recourse solution set is finite after selection of a first-stage set of transplants.
    \item We are the first to introduce two feasibility-seeking reformulations of the second stage, as opposed to optimality-based formulations \citep{Carvalho2021, Blom2021}, which improves scalability since the number of decision variables in our formulation grows linearly with the number of vertices and arcs. 
    \item We derive dominating scenarios and explore several second-stage solution algorithms to overcome the drawback of a lack of lower bound in our solution framework.
    \item We show that our framework results in significant computational and solution quality improvements compared to state-of-the-art algorithms.
\end{enumerate}

The remainder of the paper is organized as follows. Section \ref{LitReview} presents a collection of related works. Section \ref{Preliminaries} establishes the problem we address. Sections \ref{RobustModels} and \ref{SecondStage} present the first and second-stage formulations, respectively. Section \ref{sec:HSAs} presents the full algorithmic framework. Section \ref{Experiments} shows computational results. Lastly, Section \ref{sec:Conclusion} draws conclusions and states a path for future work.

\section{Related Work}
\label{LitReview}
   
\cite{Abraham2007} and \cite{Roth2007} introduced the first \KPDP\ mixed-integer programming (MIP) formulations for maximizing the number/weighted sum of exchanges. These formulations are the well-known \textit{edge formulation} and \textit{cycle formulation}. The edge formulation
uses arcs in the input graph to index decision variable, whereas
the cycle formulation, which was initially proposed for cycles only, has a decision variable for every feasible cycle and chain in the input graph. Both formulations are of exponential size either in the number of constraints (edge formulation) or in the number of decision variables (cycle formulation). However, the cycle formulation, along with a subsequent formulation \citep{Dickerson2016}, is the MIP formulation with the strongest linear relaxation. Due to its strength and natural adaptability, multiple works have designed branch-and-price algorithms employing the cycle formulation. The branch-and-price algorithm in \citet{Abraham2007} was effective for cycles of size up to three, while that of \citet{Lam2020} solved the problem for long cycles, and \cite{Riascos2020} used decision diagrams to solve large instances with both long cycles and long chains for the first time. More recently, \cite{Omer2022} built on the work in \citet{Riascos2020} by implementing a branch-and-price able to solve remarkably large instances (10,000 pairs and 1,000 non-directed donors), opening the door to large-scale multi-hospital and multi-country efforts. Another trend has focused on new arc-based formulations
   (e.g., \cite{Constantino2013,Dickerson2016}) and arc-and-cycle-based formulations (e.g., \cite{Anderson2015, Dickerson2016}). Between these two approaches, arc-and-cycle-based formulations seem to outperform arc-based formulations \citep{Dickerson2016}, especially for instances with cycles with at most three exchanges.
    
The previously discussed studies do not consider uncertainty in the proposed exchanges. However, the high percentage of planned transplants that end up cancelled suggests a need to plan for uncertainty. There are two sources of uncertainty that have been studied in the literature: weight accuracy (e.g., \cite{Duncan2019}) and vertex/arc existence (e.g., \cite{Dickerson2016, Klimentova2016, Duncan2019, Carvalho2021, Smeulders2022}). Weight accuracy uncertainty assumes that the social benefit (weight) associated with an exchange can vary, e.g., due to changes in a patient's health condition, and from the existence of multiple opinions from policy makers on the priority that should be given to each patient \citep{Duncan2019}. Uncertainty in the existence of a vertex/arc, i.e., whether or not a patient or donor leaves the exchange or compatibility between a patient and donor changes, has received greater attention. There are three main approaches in the literature addressing vertex or arc existence as the source of uncertainty: (1) a maximum expected value approach; (2) an identification of exchanges for which a physical crossmatch test should be performed to maximize the expected number of realized transplants; and (3) a maximization of the number of transplants under the worst-case disruption of vertices and arcs.
    
    The maximum expected value approach is the approach most investigated in the literature. It is concerned with finding the set of transplants with maximum expected value, i.e., a set of transplants that is most likely to yield either the maximum number of exchanges, or the maximum weighted sum of exchanges given some vertex/arc failure probabilities. This approach has mostly been modeled as a deterministic kidney exchange problem (KEP), where the objective function approximates the expected value of a matching using the given probabilities as objective coefficient multipliers of deterministic decisions.
    
    \cite{Awasti2009} considered the failure of vertices in an online setting of the cycle-only version for cycles with at most three exchanges. The authors generate sample trajectories on the arrival of patients/donors and patients survival, then use a REGRETS algorithm as a general framework to approximate the collection of cycles with maximum expectation. \cite{Dickerson2012} proposed a heuristic method to learn the ``potential'' of structural elements (e.g., vertex), that quantifies the future expected usefulness of that element in a changing graph with new patient/donor arrivals and departures.  \cite{Dickerson2013} considered arc failure probabilities and found a matching with maximum expected value, but solution repairs for failures are not considered. This work is extended in \citet{Dickerson2019}. \cite{Klimentova2016} computed the expected number of transplants for the cycle-only version while considering \textit{internal} recourse and \textit{subset} recourse to recover a solution in case of vertex or arc failure. Internal recourse, also known as \textit{back-arcs recourse} (e.g., \cite{Carvalho2021}) allows surviving pairs to match among themselves, whereas subset recourse allows a wider subset of vertices to participate in the repaired solution. To compute the expectation, an enumeration tree is used for all possible failure patterns in a cycle and its extended subset. This subset consists of the additional vertices (for the subset recourse only) such that the pairs in the original cycle can form feasible cycles. To limit the size of the tree, the subset recourse is limited to a small subset of extra vertices and the internal recourse seems to scale for short cycles only. \cite{Alvelos2019} proposed to find the expected value for the cycle-only version while considering internal recourse through a branch-and-price algorithm, finding that the overall run time grew rapidly with the size of the cycles.
    
 To identify exchanges where a physical crossmatch test should be performed, \cite{Blum2013} modeled the KEP in an undirected graph representing pairwise exchanges only. They proposed to perform two physical crossmatch tests per patient-donor pair---one for every arc in a cycle of size two---before exchanges are selected with the goal of maximizing the expected number of transplants. They showed that their algorithm yields near-optimal solutions in polynomial time. Subsequent works \citep{Assadi2019,Blum2020} evaluated adaptive and non-adaptive policies to query edges in the graph. In the same spirit, but for the general kidney exchange problem (with directed cycles and chains), \cite{Smeulders2022} formulated the maximization of the expected number of transplants as a two-stage stochastic integer programming problem with a limited budget on the number of arcs that can be tested in the first stage. Different algorithmic approaches were proposed, but scalability was a challenge.
    
In addressing worst-case vertex/arc disruption, \cite{Duncan2019} found robust solutions with no recourse for budget failure of the number of arcs that fail in the graph. \cite{Carvalho2021} proposed a two-stage robust optimization model that allowed recovery of failed solutions through the back-arcs recourse and full recourse policies. The full recourse policies can be considered a subset recourse policy (e.g., \cite{Klimentova2016}), in which all vertices that were not selected in the matching can be included in the repaired solution. 

Unlike our work, vertex and arc failure probabilities are treated as homogeneous, i.e., both elements fail with the same probability. Since in homogeneous failure there is a worst-case scenario in which all failures are vertex failures, the recourse policies are evaluated under vertex failure only. The back-arcs recourse policy only scales for instances with 20 vertices, whereas the full-recourse policy scales for instances up to 50 vertices. \citet{Blom2021} examined the general robust model for the full-recourse policy studied in \citet{Carvalho2021} and showed its structure to be a defender-attacker-defender model. Two Benders-type approaches were proposed and tested using the same instances from \citet{Carvalho2021}, and showed improved performance over previous branch-and-bound algorithms \citep{Carvalho2021}. This approach, however, is limited to homogeneous failure for the full-recourse policy. 

In our work, we allow for different failure rates between vertices and arcs, and present a solution scheme that can address recourse policies whose recourse solution set corresponds to a subset of the feasible cycles and chains in the compatibility graph. Specifically, we apply our solution scheme to the full recourse policy, as well as to a new recourse policy we introduce later, but our solution scheme can also be used for the back-arcs recourse. Our solution method requires a robust MIP formulation adapted to a specific policy that can be solved iteratively as new failure scenarios are added. The second-stage problem is decomposed into a master problem and a sub-problem. The master problem is formulated as the same feasibility problem regardless of the policy but it is implicit in the constraint set, whereas the sub-problem (i.e., recourse problem) corresponds to a deterministic KEP where only non-failed cycles and chains contribute to the robust objective.
 
\section{Preliminaries}
\label{Preliminaries}

In this section, we describe our two-stage robust model. A full list of notation is in Appendix \ref{sec:notation}. We begin by formally describing the first-stage problem. Specifically, we define the compatibility graph and the feasible set for the first-stage decisions. We similarly define the second-stage problem and then introduce our two-stage robust problem. Finally, we define the uncertainty set and recourse policies.

\paragraph{First-stage compatibility graph.} The KEP can be defined on a directed graph $\graphKEP$, whose vertex set $\vertexSet := \pairSet \cup \nddSet$ represents the set of patient-donor pairs, $\pairSet$, and the set of non-directed donors $\nddSet$. From this point onward, we will refer to patient-donor pairs simply as \emph{pairs}. The arc set $\arcSet \subseteq \vertexSet \times \pairSet$ contains arc $\parc$ if and only if the donor in vertex $\vertexa \in \vertexSet$ is compatible with a patient in vertex $\vertexb \in \pairSet$. 
A matching of donors and patients in the KEP can take the form of simple cycles and simple chains (i.e., paths from the digraph). 
A cycle is feasible if it has no more than $\cycleCap$ arcs, whereas a chain is feasible if it has no more than $\chainCap$ arcs ($\chainCap + 1$ vertices) and starts with a non-directed donor. 
The set of feasible cycles and chains is denoted by $\cycleSet$ and $\chainSet$, respectively. Furthermore, let $\containedVxt$ and $\containedArc$ be the set of vertices and arcs in $(\cdot)$, where $(\cdot)$ refers to a cycle, chain, etc. 

\paragraph{Feasible set of first-stage decisions.} A feasible solution to the KEP corresponds to a collection of vertex-disjoint feasible cycles and chains, referred to as a \emph{matching}\footnote{Note that since each vertex in the KEP compatibility graph corresponds to a patient-donor pair or a non-directed donor, a KEP matching---which is referred to as a matching in the literature and in this paper---is not actually a matching of the digraph. Instead it is a collection of simple cycles and paths in the digraph which leads to an actual underlying 1-1 matching of selected donors and patients.}, i.e., $\matchSet \subseteq \cycleSet \cup \chainSet$ such that $V(\match) \cap V(\match^{\prime}) = \emptyset,$ for all $\match, \match^{\prime} \in \matchSet$ with $\match \ne \match^{\prime}$, noting that $\match$ and $\match^\prime$ are entire cycles or chains in the matching. 
We let $\matchSpace$ denote the set of all KEP matchings in graph $\graph$. Also, we define  
$
\setFirst := \{\charVecFirst: \matchSet \in \matchSpace \}
$
as the set of all binary vectors representing the selection of a feasible matching, where $\charVecFirst$ is the characteristic vector of matching $\matchSet$ in terms of the cycles/chains sets $\cycleSet \cup \chainSet$. That is, $\charVecFirst \in \{0,1\}^{\mid \cycleSet \cup \chainSet \mid}$ with $\varFirstSel_{\matchSet,\match} = 1$  if and only if $\match \in \matchSet$, meaning that a patient in a pair obtains a transplant if it is in some cycle/chain $\match$ selected in matching $\matchSet$.

\paragraph{Transitory compatibility graph.}

Upon selection of a solution $\varFirst \in \setFirst$ but before uncertainty is revealed, there exists a \emph{transitory} compatibility graph $\digraphSndTr = (\vertexSetSndBe, \arcSetSndBe)$ with vertex set $\vertexSetSndBe = \cup_{\match \in \cycleSetSndTr \cup \chainSetSndTr}V(\match)$ and arc set $\arcSetSndBe = \cup_{\match \in \cycleSetSndTr \cup \chainSetSndTr}A(\match)$. The set of feasible cycles $\cycleSetSndTr$ and feasible chains $\chainSetSndTr$  in the transitory compatibility graph (i) are \textit{allowed} under recourse policy $\policy \in \policySet$, and (ii) have \textit{at least one pair} in $\varFirst$. Sections \ref{subsec:uncertaintyset} and \ref{sec:RecoursePolicies} provide details on uncertainty set  $\bigsCe$ and the types of recourse policies $\policySet$, respectively, and provide an explicit definition of condition (i). Condition (ii) enforces that each $\match \in \cycleSetSndTr \cup \chainSetSndTr$ satisfies $V(\match) \cap V(\varFirst) \cap \pairSet \neq \emptyset$, where $V(\varFirst) = \cup_{\match' \in \cycleSet \cup \chainSet:\varFirst_{\match'} = 1} V(\match')$.

\paragraph{Second-stage compatibility graph.} 

Once a failure scenario $\oneCe \in \bigsCe$ is observed in the first-stage compatibility graph, a digraph $\digraphSndS$ for the second-stage problem is \emph{induced} in the transitory graph, thus revealing the cycles and chains with selected pairs from the first stage that were unaffected by $\oneCe$. We refer to $\vertexSetSndBeF$ and  $\arcSetSndBeF$ as the set of vertices and the set of arcs that fail under scenario $\oneCe$ in the first-stage compatibility graph, respectively. We then define a second-stage compatibility graph $\digraphSndS = (\vertexSetSndBe \setminus \{\vertexSetSndBe \cap \vertexSetSndBeF\}, \arcSetSndBe \setminus \{\arcSetSndBe \cap \arcSetSndBeF\})$ as the sub-graph obtained by removing  the vertices and arcs that fail under scenario $\oneCe$ in the first-stage compatibility graph, and are also present in the transitory graph, from the transitory compatibility graph. Thus, we define $\cycleSetSndS$ and $\chainSetSndS$ as the set of \emph{non-failed} cycles and chains in the second-stage compatibility graph, respectively.

\paragraph{Feasible set of second-stage decisions.}
A solution to the second stage is referred to as a \textit{recourse solution} in $\digraphSndS$ under some scenario $\oneCe \in \bigsCe$, leading to an alternative matching where pairs from the first-stage solution $\varFirst \in \setFirst$ are re-arranged into non-failed cycles and chains, among those allowed under policy $\policy \in \policySet$.
We can now define $\polySetxOne:= \{\matchSet \subseteq \cycleSetSndS \cup \chainSetSndS \mid V(\match) \cap V(\match^{\prime}) = \emptyset \text{ for all } \match, \match^{\prime} \in \matchSet \text{; } \match \neq \match^{\prime}\}$ as the set of allowed recovering matchings under policy $\policy$ such that every cycle/chain in  $\polySetxOne$ contains at least one pair in $\varFirst$. In other words, a recourse solution must contain at least one pair from the first-stage solution.
Thus, similar to the first-stage characteristic vector $\charVecFirst$, let
$
\setReco:= \{\charVecSecond:\matchSet \in \polySetxOne \}
$
be the set of all binary vectors representing the selection of a feasible matching in a second-stage compatibility graph with non-failed elements (vertices/arcs), under scenario $\oneCe \in \bigsCe$ and policy $\policy \in \policySet$ that contain at least one  pair in $\varFirst$.

\paragraph{Two-stage RO problem.}

A general two-stage RO problem for the KEP can then be defined as follows: 
\begin{align}
	\label{ROmodel0}
	\max_{\varFirst \in \setFirst} \min_{\oneCe \in \bigsCe} \max_{\varReco \in \setReco} f(\varFirst,\oneCe,\varReco)&
\end{align} 
A set of transplants given by solution $\varFirst \in \setFirst$ (the outer maximization) is selected in the first stage. Then, the uncertainty vector $\oneCe \in \bigsCe$ is observed (the middle minimization) and a recourse solution $\varReco \in \setReco$ is found to repair $\varFirst$ according to recourse policy $\policy \in \policySet$ (the inner maximization). The \emph{second stage}, established by the min-max problem ($\max_{\varFirst \in \setFirst} \min_{\oneCe \in \bigsCe}$), finds a recourse solution by solving the \emph{recourse problem} ($\max_{\varReco \in \setReco}$), whose objective value maximizes $f(\varFirst,\oneCe,\varReco)$ under failure scenario $\oneCe \in \bigsCe$, but it is the lowest objective value among all failure scenarios. The scenario optimizing the second-stage problem is then referred to as the \emph{worst-case scenario} for solution $\varFirst \in \setFirst$. 

The \emph{recourse objective function}, $f(\varFirst,\oneCe,\varReco)$, assigns weights to the cycles and chains of a recovered matching associated with a recourse solution $\varReco$ under failure scenario $\oneCe$, based on the number of pairs matched in the first stage solution $\varFirst$. 
Thus, we define the recourse objective as 
\begin{equation}
\label{eq:RecourseCostFunction}
f(\varFirst,\oneCe,\varReco) := \sum_{\match \in \cycleSetSndS \cup \chainSetSndS} \mathbf{w}_{\match} (\varFirst) \varReco_{\match} := \weightRO\varReco,
\end{equation}
where $\mathbf{w}_{\match}(\varFirst) := |V(\match) \cap  \SelVerX \cap \pairSet|$ is the weight of a cycle/chain $\match \in \cycleSetSndS \cup \chainSetSndS$ corresponding to the \textit{number of pairs} that were matched in the first stage by solution $\varFirst \in \setFirst$, and can now also be matched in the second stage by recourse solution $\varReco \in \setReco$, after failures are observed. As a result, the weight of every cycle and chain selected by $\varReco$ is the number of pairs present from the first-stage solution $\varFirst$. Accordingly, an optimal KEP robust solution is a matching in the first stage which, under the worst-case scenario, has the fewest pairs disrupted by the recovery in the second stage.
For the sake of a compact representation, and with a slight abuse of notation, we represent the recourse cost function as the inner product given on the right-hand side of Equation \eqref{eq:RecourseCostFunction}, where the superscript $\oneCe$ on the weight vector determines its dimension and in turn makes it consistent with that of the recourse decision vector. 
Thus, our two-stage RO problem for the KEP is defined as follows:
\begin{align}
	\label{ROmodel}
	\max_{\varFirst \in \setFirst} \min_{\oneCe \in \bigsCe} \max_{\varReco \in \setReco} \objRO&
\end{align}

\subsection{Uncertainty set} 
\label{subsec:uncertaintyset}
Failures in a planned matching are observed after a final checking process, leading to the removal of the affected vertices and arcs from the first-stage compatibility graph. The failure of a vertex/arc causes the failure of the entire cycle to which that element belongs and the shortening of the chain at the last vertex before the first failure.

Unlike previous studies, we consider \textit{non-homogeneous failure} by allowing vertices and arcs to fail at different rates. We define two failure budgets, one for vertices and one for arcs; in other words, we assume that there exist two unknown probability distributions causing vertices and arcs to fail independently from one another, rather than assuming that both vertices and arcs follow the same failure probability distribution. We note that this approach can still model homogeneous failure, which is a special case of non-homogeneous failure.

Traditionally, the uncertainty set in robust optimization does not depend on the first-stage decision. While we consider vertex and arc failures as exogenous random variables, we define the uncertainty set $\bm{\bigsCe}$ in terms of  all the uncertainty sets $\bm{\bigsCe}(\varFirst)$. Failure scenarios in $\bm{\bigsCe}(\varFirst)$ lead to a second-stage compatibility graph $\digraphSndS$ denoting the dependency between the selection of a first-stage decision and the ability to observe a failure. Recall that the checking process leading to the discovery of failures is based on the matching proposed for transplant. In other words, failures that involve vertices and arcs not present in the transitory graph cannot be discovered even if they occur since the set of alternative cycles and chains in the transitory graph that can be used to repair a first-stage solution $\varFirst$ are defined by policy $\policy$ regardless of the failure scenario $\oneCe$. Therefore, we define the uncertainty set $\bm\bigsCe$ as follows:
\begin{subequations}
\label{UDef}
\begin{align}
\bm{\bigsCe}(\varFirst) &:= \left(\bm{\UbigsCe}(\varFirst), \bm{\AbigsCe}(\varFirst)\right) \text{ where,}\\
\bm{\UbigsCe}(\varFirst)  &:=\{\bm{\smallsCeVx} \in \{0,1\} ^{\mid\vertexSet \mid }\mid  \lvert \vertexSetSndBe  \cap \vertexSetSndBeF \rvert \le \sum_{\vertex \in \vertexSet} \smallsCeVx_{\vertex} \le \vBudget\}\\
\bm{\AbigsCe}(\varFirst)  &:= \{\bm{\smallsCeArc} \in \{0,1\} ^{\mid \arcSet \mid}\mid \lvert \arcSetSndBe \cap \arcSetSndBeF \rvert \le \sum_{\parc \in \arcSet} \smallsCeArc_{\narc} \le \aBudget  \}\\
\bm{\bigsCe} &:= \bigcup\limits_{\varFirst \in \setFirst} \bm{\bigsCe}(\varFirst)
\end{align}
\end{subequations}

A failure scenario $\oneCe = (\smallsCeVx, \smallsCeArc)$ is represented by binary vectors $\smallsCeVx$ and $\smallsCeArc$, where $\smallsCeVx$ refers to vertex failures and $\smallsCeArc$ refers to arc failures. A vertex $\vertex \in \vertexSet$ and arc $\parc \in \arcSet$ from the first-stage compatibility graph fail under a realized scenario if $\smallsCeVx_{\vertex} = 1$ or $\smallsCeArc_{\narc} = 1$, respectively. The total number of vertex and arc failures in the first-stage compatibility graph is controlled by user-defined parameters $\vBudget$ and $\aBudget$, respectively. Therefore, the number of vertex failures in the transitory compatibility graph $\digraphSndTr$ leading to a second-stage compatibility graph  $\digraphSndS$ cannot exceed $\vBudget$. Likewise, the number of arc failures in $\digraphSndTr$ cannot exceed $\aBudget$. Thus, the uncertainty set $\bm\bigsCe$ is the union over all failure scenarios leading to a second-stage compatibility graph when a set of transplants in the first stage has been proposed for transplantation. 
Note that this uncertainty set definition only distinguishes failures by element type (vertex or arc), and does not distinguish between sensitized and non-sensitized patients.

\subsection{Recourse policies}
\label{sec:RecoursePolicies}
An important consideration made by KPDPs is the guideline according to which a selected matching is allowed to be repaired. Although re-optimizing the deterministic KEP with non-failed vertices/arcs is an option \citep{Dutch2005}, some KPDPs opt for recovery strategies when failures in the matching given by the deterministic model are observed \citep{CBS2019, Manlove2015}. Thus, it is reasonable to use those strategies as recourse policies when uncertainty is considered. We consider the \textit{full-recourse policy} studied in \citet{Blom2021, Carvalho2021} and introduce a natural extension of this policy, which we refer to as \textit{the first-stage-only recourse policy}. Examples of each are shown in Figure \ref{fig:RecourseExamples}.

\subsubsection{Full recourse policy}
\begin{figure}[tbp]
	\centering	
	\begin{adjustbox}{minipage=\linewidth,scale=1}
		\begin{subfigure}[t]{.48\textwidth}
			\centering
			\tikzstyle{place}=[circle,draw=blue!50,fill=blue!20,thick,
			inner sep=0pt,minimum size=4mm]
			\tikzstyle{transition}=[circle,draw=black!50,fill=black!20,thick,
			inner sep=0pt,minimum size=4mm]
			\tikzstyle{dummy}=[rectangle,draw=white,fill=white,thick,
			inner sep=0pt,minimum size=4mm]
			\begin{tikzpicture}
				
				\node[transition] (cuatro) at (0,2) {$4$}; 
				\node[transition] (tres) at (1,1) {$3$};
				\node[transition] (dos) at (-1,1) {$2$};
				\node[transition] (nueve) at (0,1) {$9$};
				\node[transition] (diez) at (-2,1) {$10$};
				
				\node[transition] (cinco) at (2,1) {$5$};
				\node[transition] (uno) at (4,1)  {$1$}; 
				\node[transition] (ocho) at (4,2.5)  {$8$}; 
				\node[transition] (seis) at (3,2)  {$6$}; 
				\node[transition] (siete) at (3,1)  {$7$}; 
				
				\draw [-{latex},thick] (dos) to [bend left=45] (nueve);
				\draw [-{latex},thick] (nueve) to [bend left=45] (diez);
				\draw [-{latex},thick] (diez) to [bend left=45] (dos);
				
				\draw [-{latex}, dashed] (cuatro) to [bend right=45] (dos);
				\draw [-{latex}, dashed] (dos) to [bend right=45] (tres);
				\draw [-{latex},thick] (cuatro) to [bend right=45] (tres);
				\draw [-{latex},thick] (tres) to [bend right=45] (cuatro);
				
				\draw [-{latex},thick] (seis) to [bend left=45] (uno);
				\draw [-{latex}, thick] (uno) to [bend left=45] (cinco);
				\draw [-{latex},dashed] (siete) to [bend left=45] (uno);
				\draw [-{latex},dashed] (cinco) to [bend left=45] (siete);
				\draw [-{latex},dashed] (ocho) to [bend left=45] (uno);
				\draw [-{Classical TikZ Rightarrow},thick] (cinco) to [bend left=45] (seis);
				
				\draw [-{latex},dashed] (tres) to [bend left=45] (seis);
				\draw [-{latex},dashed] (seis) to [bend right=35] (cuatro);
			\end{tikzpicture}
			\caption{First-stage compatibility graph $\graphKEP$ and first-stage solution with eight transplants depicted in bold arcs}
			\label{fig:CompGraph}
		\end{subfigure}
            \hspace{0.5cm}
		\begin{subfigure}[t]{.48\textwidth}
			\centering
    	    \tikzstyle{place}=[circle,draw=blue!50,fill=blue!20,thick,
			inner sep=0pt,minimum size=4mm]
			\tikzstyle{transition}=[circle,draw=black!50,fill=black!20,thick,
			inner sep=0pt,minimum size=4mm]
			\tikzstyle{dummy}=[rectangle,draw=white,fill=white,thick,
			inner sep=0pt,minimum size=4mm]
			\begin{tikzpicture}
				
				\node[transition] (cuatro) at (0,2) {$4$}; 
				\node[transition] (tres) at (1,1) {$3$};
				\node[transition, draw=red] (dos) at (-1,1) {$2$};
				\node[transition] (nueve) at (0,1) {$9$};
				\node[transition] (diez) at (-2,1) {$10$};
				
				\node[transition] (cinco) at (3,1) {$5$};
				\node[transition] (uno) at (5,1)  {$1$}; 
				\node[transition] (ocho) at (5,2.5)  {$8$}; 
				\node[transition] (seis) at (4,2)  {$6$}; 
				\node[transition] (siete) at (4,1)  {$7$}; 
				
				\draw [-{latex},thick] (dos) to [bend left=45] (nueve);
				\draw [-{latex},thick] (nueve) to [bend left=45] (diez);
				\draw [-{latex},thick] (diez) to [bend left=45] (dos);
				
				\draw [-{latex}, thick] (cuatro) to [bend right=45] (dos);
				\draw [-{latex}, thick] (dos) to [bend right=45] (tres);
				\draw [-{latex}, thick] (cuatro) to [bend right=45] (tres);
				\draw [-{latex},thick] (tres) to [bend right=45] (cuatro);
				
				\draw [-{latex},thick] (seis) to [bend left=45] (uno);
				\draw [-{latex}, thick] (uno) to [bend left=45] (cinco);
				\draw [-{latex},dashed] (ocho) to [bend left=45] (uno);
				\draw [-{latex},thick, draw=red] (cinco) to [bend left=45] (seis);
				\draw [-{latex},dashed] (siete) to [bend left=45] (uno);
				\draw [-{latex},dashed] (cinco) to [bend left=45] (siete);
				
				\draw [-{latex}, dashed] (tres) to [bend left=45] (seis);
				\draw [-{latex}, dashed] (seis) to [bend right=35] (cuatro);
				
				\draw[dashed, -latex,preaction = {draw,blue!30,-,double=blue!30,double distance = 5\pgflinewidth}] (ocho) to [bend left=45] (uno);
				\draw[thick, -latex,preaction = {draw,blue!30,-,double=blue!30,double distance = 3\pgflinewidth}] (uno) to [bend left=45] (cinco);
				\draw[thick, -latex,preaction = {draw,blue!30,-,double=blue!30,double distance = 3\pgflinewidth}] (cuatro) to [bend right=45] (tres);
				\draw[dashed, -latex,preaction = {draw,blue!30,-,double=blue!30,double distance = 4\pgflinewidth}] (seis) to [bend right=35] (cuatro);
				\draw[dashed, thick, -latex,preaction = {draw,blue!30,-,double=blue!30,double distance = 3\pgflinewidth}] (tres) to [bend left=45] (seis);
				
			\end{tikzpicture}
			\caption{Full recourse solution in shaded arcs, with five transplants including pairs from the first-stage solution}
			\label{fig:FullRecourseGraph}
		\end{subfigure}
  
  \centering
  \begin{subfigure}[t]{.52\textwidth}
			\centering
    	    \tikzstyle{place}=[circle,draw=blue!50,fill=blue!20,thick,
			inner sep=0pt,minimum size=4mm]
			\tikzstyle{transition}=[circle,draw=black!50,fill=black!20,thick,
			inner sep=0pt,minimum size=4mm]
			\tikzstyle{dummy}=[rectangle,draw=white,fill=white,thick,
			inner sep=0pt,minimum size=4mm]
			\begin{tikzpicture}
				
				\node[transition] (cuatro) at (0,2) {$4$}; 
				\node[transition] (tres) at (1,1) {$3$};
				\node[transition, draw=red] (dos) at (-1,1) {$2$};
				\node[transition] (nueve) at (0,1) {$9$};
				\node[transition] (diez) at (-2,1) {$10$};
				
				\node[transition] (cinco) at (3,1) {$5$};
				\node[transition] (uno) at (5,1)  {$1$}; 
				\node[transition] (ocho) at (5,2.5)  {$8$}; 
				\node[transition] (seis) at (4,2)  {$6$}; 
				\node[transition] (siete) at (4,1)  {$7$}; 
				
				\draw [-{latex},thick] (dos) to [bend left=45] (nueve);
				\draw [-{latex},thick] (nueve) to [bend left=45] (diez);
				\draw [-{latex},thick] (diez) to [bend left=45] (dos);
				
				\draw [-{latex}, thick] (cuatro) to [bend right=45] (dos);
				\draw [-{latex}, thick] (dos) to [bend right=45] (tres);
				\draw [-{latex}, thick] (cuatro) to [bend right=45] (tres);
				\draw [-{latex},thick] (tres) to [bend right=45] (cuatro);
				
				\draw [-{latex},thick] (seis) to [bend left=45] (uno);
				\draw [-{latex}, thick] (uno) to [bend left=45] (cinco);
				\draw [-{latex},dashed] (ocho) to [bend left=45] (uno);
				\draw [-{latex},thick, draw=red] (cinco) to [bend left=45] (seis);
				\draw [-{latex},dashed] (siete) to [bend left=45] (uno);
				\draw [-{latex},dashed] (cinco) to [bend left=45] (siete);
				
				\draw [-{latex}, dashed] (tres) to [bend left=45] (seis);
				\draw [-{latex}, dashed] (seis) to [bend right=35] (cuatro);
				
				\draw[thick, -latex,preaction = {draw,blue!30,-,double=blue!30,double distance = 3\pgflinewidth}] (cuatro) to [bend right=45] (tres);
				\draw[dashed, -latex,preaction = {draw,blue!30,-,double=blue!30,double distance = 4\pgflinewidth}] (seis) to [bend right=35] (cuatro);
				\draw[dashed, thick, -latex,preaction = {draw,blue!30,-,double=blue!30,double distance = 3\pgflinewidth}] (tres) to [bend left=45] (seis);
				
			\end{tikzpicture}
			\caption{First-stage-only recourse solution in shaded arcs, with three transplants including pairs from the first-stage solution}
			\label{fig:FirstStageOnlyRecourseGraph}
		\end{subfigure}
	\end{adjustbox}
	\caption{Recourse policies when vertex 2 and arc (5,6) in first-stage compatibility graph (a) fail.}
	\label{fig:RecourseExamples}
\end{figure}
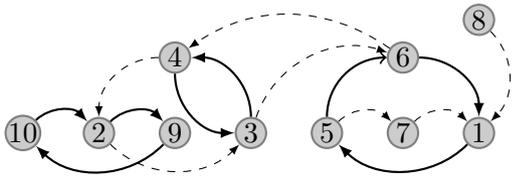
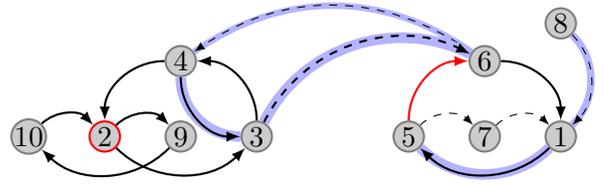
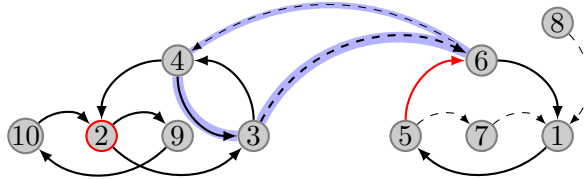

Under the full-recourse policy, pairs that were selected in the first stage belonging to failed components are allowed to be re-arranged in the second stage within non-failed cycles and chains that may involve any other vertex (pair or non-directed donor), regardless of whether that vertex was selected or not in the first stage. Figure \ref{fig:FullRecourseGraph} shows an example of the full-recourse policy arising from the compatibility graph in Figure \ref{fig:CompGraph}. The first-stage solution depicted with bold arcs has a total weight of eight, since there are eight exchanges. Suppose there is a scenario in which $\smallsCeVx_{2} = 1$ and $\smallsCeArc_{56} = 1$. Assuming that $\cycleCap = \chainCap = 4$, then the best recovery plan under this scenario, depicted by the recourse solution with shaded arcs, is to re-arrange vertices 3, 4 and 6 by bringing them together into a new cycle and include vertices 1 and 5 in a chain started by non-directed donor 8 (which was not selected in the first stage). Alternatively, vertices 1 and 5 could be selected along with vertex 7 to form a cycle. In both cases, the recourse solution involves only five pairs from the first stage.

\subsubsection{First-stage-only recourse policy}

We refer to the first-stage-only as a recourse policy in which only vertices selected in the first stage can be used to repair a first-stage solution, i.e., the new non-failed cycles and chains selected in the second stage must include vertices from that first-stage solution only. The recourse solution set of the first-stage-only recourse policy corresponds to a subset of the full recourse policy. Although more conservative, KPDPs can opt for the first-stage-only since in the full recourse there can be vertices selected in the second stage that were not selected in the first stage, and thus have not yet been checked by KPDP operators, adding additional uncertainty about the actual presence of vertices/arcs. The back-arcs recourse policy studied in \citet{Carvalho2021} also allows recovery of a first-stage solution with pairs that were selected in that solution only, but such a policy only allows recovery of a cycle (chain) if there exists other cycles (chains) nested within it, making the back-arcs recourse more conservative than the first-stage-only recourse. In Figure \ref{fig:FirstStageOnlyRecourseGraph}, the recourse solution under the first-stage-only recourse includes vertices 3, 4 and 6 in a cycle just as in the full-recourse policy, but chains started by vertex 8 and the cycle to which vertex 7 belongs would not be within the feasible recourse set. Thus, the recourse objective value corresponds to only three pairs from the first stage as opposed to five pairs.

\section{Robust model and first stage}
\label{RobustModels}

For a finite (yet possibly large) uncertainty set, the second optimization problem in Model \eqref{ROmodel} ($\min_{\oneCe \in \bigsCe}$) can be removed as an optimization level and included as a set of scenario constraints instead \citep{Carvalho2021}. Model \eqref{ROmodel} can be expressed as the following MIP robust formulation: 
\begin{subequations}
	\label{SingleROmodel0}
    \begin{align}
	\SingleRO: \quad \max_{\varFirst, \objSingleRO} \quad& \objSingleRO \label{eq:SingleROobj} \\
	& \objSingleRO \le  \max_{\varReco \in \setReco} \objRO & \oneCe \in \bigsCe  \label{eq:scenarios}\\
	&\varFirst \in \setFirst
\end{align} 
\end{subequations}
The optimization problem in \eqref{eq:scenarios} is the recourse problem, which given a first-stage solution $\varFirst \in \setFirst$ and policy $\policy \in \policySet$, it finds a recourse solution in the set of binary vectors $\setReco$ whose cycles and chains under scenario $\oneCe \in \bigsCe$ have the largest number of pairs from the first stage. Having the recourse problem in the constraints is challenging since it must be solved as many times as the number of failure scenarios $\oneCe \in \bigsCe$ for any fixed decision $\varFirst$, which can be prohibitively large. To reach a more tractable form, first observe that the directions of the external (Objective \eqref{eq:SingleROobj}) and internal (Constraint \eqref{eq:scenarios}) objectives in Model \eqref{SingleROmodel0} are maximization. We can then define a binary decision vector $\varReco^{\oneCe}$ for every scenario $\oneCe \in \bigsCe$ whose feasible space corresponds to a recourse solution in $\setReco$. Thus, an equivalent model for the first-stage formulation (FSF) has the inner optimization problem removed and the set $\setReco$ included as part of the constraints set for every failure scenario. The resulting FSF model is as follows:
\begin{subequations}
    \begin{align}
	\tag{FSF} \SingleRO=  \max_{\varFirst, \varReco^{\oneCe}, \objSingleRO} \quad& \objSingleRO \label{SingleROmodel}\\
	& \objSingleRO \le  \objRO & \oneCe \in \bigsCe  \label{eq:scenariosBound}\\
 	&\varReco^{\oneCe} \in \setReco & \oneCe \in \bigsCe  \label{eq:scenariosReco}\\
 	&\varFirst \in \setFirst
\end{align} 
\end{subequations}
The first-stage solution and its corresponding recourse solution can be found in a single-level MIP formulation if the set of scenarios $\bigsCe$ can be enumerated. Observe that once a scenario $\oneCe \in \bigsCe$ is fixed, $\setReco$ can be expressed by a set of linear constraints and new decision variables, e.g., the position-indexed formulation for the robust KEP under full recourse proposed in \citet{Carvalho2021}. This formulation follows the structure of Model \eqref{SingleROmodel}. We present this formulation (Appendix \ref{FullFirstSFormulation}) and a variant of it that satisfies the first-stage-only recourse policy (Appendix \ref{sec:fisrtOnlyROFormulation}) in the Online Supplement.

To overcome the large size of $\bigsCe$, we solve the robust KEP in Model \eqref{SingleROmodel} iteratively by finding a new failure scenario whose associated optimal recourse solution $\varReco^{\oneCe} \in \setReco$ recovers the maximum number of pairs from the first stage under that scenario in Model \eqref{SingleROmodel} (Algorithm \ref{Alg:ScenarioGeneration}). We note that adding new failure scenarios worsens the optimal objective value in Model \eqref{SingleROmodel}. The goal is to find a scenario that yields a recourse solution with an objective value lower than the current objective value of the robust KEP model, if such scenario exists. If a scenario is found, then a new set of recourse decision variables along with Constraints \eqref{eq:scenariosBound} and \eqref{eq:scenariosReco} are added to Model \eqref{SingleROmodel}.

Algorithm \ref{Alg:ScenarioGeneration} starts with an empty restricted set of scenarios $\tilde{\bigsCe}$. First, $\SingleROr$ is solved and its incumbent solution $\varFirstc$ is retrieved (Step \ref{alg:R0-solve}). Then, the second-stage problem is solved to optimality for an incumbent solution $\varFirstc$ under policy $\policy$ (Step \ref{alg:RO-2ndstage}). The solution to the second-stage problem yields a failure scenario $\optCe \in \bigsCe$ with the lowest optimal recourse objective value among all possible failure scenarios. We refer to that scenario as the worst-case scenario for a first-stage solution $\varFirstc \in \setFirst$. If such recourse objective value is lower than $\objSingleROr$, then $\optCe$ is added to the restricted set of scenarios $\tilde{\bigsCe}$. Then, the corresponding Constraints \eqref{eq:scenariosBound} and \eqref{eq:scenariosReco} are added to Model \eqref{SingleROmodel} along with decision vector $\varReco^{\optCe}$, and Model \eqref{SingleROmodel} is re-optimized. We refer to Model \eqref{SingleROmodel} with a restricted set of scenarios $\tilde{\bigsCe}$ as the first-stage problem or first-stage formulation to imply that the search for the optimal robust solution corresponding to some $\varFirstc \in \setFirst$ continues. In principle, as long as a failure scenario $\oneCe \in \bigsCe$ leads to an optimal recourse solution with an objective value lower than $\objSingleROr$, that scenario could be added to Model \eqref{SingleROmodel}. However, for benchmark purposes, at Step \ref{alg:R0-solve} the worst-case scenario is found for a given first-stage solution $\varFirstc$ under policy $\policy$. If the optimal objective value of the second-stage problem equals $\objSingleROr$, then an optimal robust solution is returned.

\begin{algorithm}[tbp]

	\algorithmicrequire { A policy $\policy \in \policySet$ and restricted set of scenarios $\tilde{\bigsCe}$, $\tilde{\bigsCe}:= \emptyset$}\\
	\algorithmicensure { Optimal robust solution $\varFirst^{\star}$}
	\begin{algorithmic}[1]
    \REPEAT
	\STATE \label{alg:R0-solve} Solve $\SingleROr$ and obtain optimal solution $\varFirstc$  \\ 
	\IF{$\bm{\min}_{\oneCe \in \bigsCe}\bm{\max}_{\varReco \in \setRecoc} \weightROc\varReco < \objSingleROr$} \label{alg:RO-2ndstage}
	\STATE $\tilde{\bigsCe} \gets \tilde{\bigsCe} \cup \{\optCe\}$, add recourse decision vector $\varReco^{\optCe}$, Constraints \eqref{eq:scenariosBound}, \eqref{eq:scenariosReco} and go to Step 1
    \ENDIF
    \UNTIL{$\bm{\min}_{\oneCe \in \bigsCe}\bm{\max}_{\varReco \in \setRecoc} \weightROc\varReco \ge \objSingleROr$}
	\STATE $\varFirst^{\star} \gets \varFirst$
    \RETURN $\varFirst^{\star}$
	\caption{Solving the robust KEP in Model \eqref{SingleROmodel}}
	\label{Alg:ScenarioGeneration}
\end{algorithmic}
\end{algorithm}

Algorithm \ref{Alg:ScenarioGeneration} converges in a finite number of iterations due to the finiteness of $\setFirst$ and $\bigsCe$. Due to the large set of scenarios, Step \ref{alg:R0-solve} is critical to efficiently solve the robust problem. In subsequent sections, we decompose the second-stage problem into a master problem yielding a failure scenario and a sub-problem (the recourse problem) finding an alternative matching with the maximum number of pairs from the first stage.

\section{New second-stage decompositions}
\label{SecondStage}

 In this section, we present two new decompositions of the second-stage problem that is solved at Step \ref{alg:RO-2ndstage} in Algorithm \ref{Alg:ScenarioGeneration} as part of the iterative solution to the robust KEP. Both consist of a feasibility-seeking master problem that finds a failure scenario and a sub-problem that finds an alternative matching under that scenario. Each decomposition solves a recourse problem in either the second-stage compatibility graph (Appendix \ref{sec:RecourseFormls}, Model \eqref{RecoPModel}) or the transitory graph (Appendix \ref{sec:RecourseFormls}, Model \eqref{RecoPModelexp}). Such models use cycle-chain decision variables as proposed in \citet{Blom2021}. Although the optimal solutions provided by both formulations are identical, we use the structure of each solution set to formulate the master problem of each decomposition accordingly.

\subsection{An initial feasibility-seeking formulation for the second stage}
\label{sec:basicfeas}
In the initial feasibility-seeking formulation, we formulate the second-stage problem as a linear binary program whose optimal solution can be obtained through another linear binary program that finds a feasible failure scenario. There exists one constraint in the feasibility-seeking formulation for every constraint in its optimality-seeking counterpart. Each constraint is associated with an optimal recourse solution under some failure scenario. However, the number of constraints may be exponential as the formulation requires knowing an optimal recourse solution under all possible failure scenarios in advance. To efficiently solve such formulations (optimality-seeking and the feasibility-seeking), we decompose them it into a master problem and sub-problem (recourse problem). The master problems finds a failure scenario which is then used by the recourse problem to obtain an optimal recourse solution, which is then added to the master problem and a new failure scenario is found.

We show that any feasible failure scenario leads to an optimal recourse solution whose objective value is an upper bound on the optimal objective value of the second-stage problem. Based on this upper bound, we also show that finding the worst-case scenario in the optimality-seeking formulation implies the failure of some of the cycles/chains in \emph{every} optimal recourse solution in that formulation. The feasibility-seeking master problem tries to find a failure scenario satisfying that at least one cycle/chain fails in every recourse solution of the optimality-seeking master problem. When a new failure scenario leads to an optimal recourse solution with the lowest objective value among the ones found in previous iterations, we update the constraints in the feasibility-seeking master problem to increase in the minimum number of cycles/chains that should fail in every recourse solution of the optimality-seeking master problem. We refer to this update  as a \emph{strengthening} procedure to indicate that the feasible space of failure scenarios in our feasibility-seeking formulation can be reduced based on a new upper bound on the optimal objective value of the second-stage problem.

Our solution framework is purely a cutting plane since no new decision variables need to be created when a new constraint is added to the feasibility-seeking master problem.  From this point onward, our feasibility-seeking formulations are shown (and referred to) as master problems, since we assume that at the start of the cutting plane algorithm we have an empty set of failure scenarios $\bigsCecR \subseteq \bigsCec$ for a given $\varFirstc \in \setFirst$ rather than the full set. We show that once the master problem becomes infeasible, an optimal solution to the second-stage problem has been found.

In this formulation, a recourse solution is found on the second-stage compatibility graph. Recall that a second-stage compatibility graph does not include cycles/chains that fail under the scenario that induced the graph in the transitory graph. Therefore, we refer to the master problem in the initial feasibility-seeking formulation as MasterSecond (MS). We first present the recourse problem formulation on the second-stage problem, then we formally define MasterSecond and prove its validity. Finally, we present non-valid inequalities to strengthen the right-hand side of constraints in MasterSecond.

\subsubsection{The recourse problem on the second-stage compatibility graph}
Given a first-stage solution $\varFirstc \in \setFirst$, let $\bigsCecR \subseteq \bigsCec$ be a restricted set of failure scenarios such that $\cdteOne \in \bigsCecR$ induces a second-stage compatibility graph $\digraphSndSc$ as defined in Section \ref{Preliminaries}. The recourse problem consists of finding a matching of maximum weight in $\digraphSndSc$, i.e., a matching with the largest number of pairs selected in the first stage. We define $\RecoPcdte$ as the MIP formulation (Appendix \ref{sec:RecourseFormls}) for the recourse problem in the second-stage compatibility graph $\digraphSndSc$, and define $\varRecoOpSndSc$ as an optimal recourse solution in $\setRecocdte$ with objective value $\RecoObjVarcdte$.
\subsubsection{The master problem: MasterSecond}
We continue to use $\smallsCeVx$ and $\smallsCeArc$ as binary vectors representing the failure of a vertex ($\smallsCeVx_{\vertex} = 1 \ \forall \vertexa \in \vertexSet$) and arc ($\smallsCeArc_{\narc} = 1 \ \forall \parc \in \arcSet$), respectively. Thus, the master problem for the second-stage problem can be formulated as follows:
\begin{subequations}
\label{mo:basicfeas}
    \begin{align}
        \tag{MS} \text{MasterSecond}(\varFirstc) \text{: }  &&& \qquad \qquad \text{Find } \oneCe \label{MasterBasic}\\
        && \sum_{\match \in \cycleSet \cup \chainSet: \varFirstc_{\match} = 1} &\left(\sum_{\vertexa \in V(\match)} \smallsCeVx_{\vertexa} +  \sum_{\parc \in A(\match)} \smallsCeArc_{\narc} \right) \ge 1   &&    \label{eq:AtLeastOneX}
        \\
         && \sum_{\match \in \cychSetSndScdte: \hat{\varReco}_{\match} = 1} &\left(\sum_{\vertexa \in V(\match)} \smallsCeVx_{\vertexa} +  \sum_{\parc \in A(\match)} \smallsCeArc_{\narc} \right) \ge 1   && \cdteOne \in \bigsCecR; \varRecoOpSndSc   \label{eq:AtLeastOne}\\
        &&& \qquad \qquad \oneCe \in \bigsCecR&& \label{basicGamma}
    \end{align}
\end{subequations}
where $\cychSetSndScdte = \cycleSetSndSc \cup \chainSetSndSc$. Because we start with an empty set of scenarios, Constraints \eqref{eq:AtLeastOneX}-\eqref{eq:AtLeastOne} correspond to covering constraints in a set cover problem, which is NP-complete \citep{Karp1972}.

At least one vertex or arc in \textit{every} recourse solution $\varRecoOpSndSc \in \setRecocdte$ \textit{must} fail before finding the optimal worst-case scenario $\oneCeStar \in \bigsCe$ for a first-stage solution $\varFirstc \in \setFirst$.  Observe that at the start, when $\bigsCecR = \emptyset$ the second-stage compatibility graph is equivalent to the transitory graph, and so the first recourse solution known is $\varFirstc$ itself. Intuitively, we can see that we can find a worse failure scenario if at least one element (vertex/arc) in $\varFirstc$ fails (Constraint \eqref{eq:AtLeastOneX}). Otherwise, the optimal objective value of the second stage would equal $\objSingleROr$ and we would have found the optimal solution to the robust KEP in Algorithm \ref{Alg:ScenarioGeneration}. Suppose we ``generate'' a failure scenario $\cdteOne$ affecting solution $\varFirstc$ which we use to induce a new second-stage compatibility graph and include that scenario in our restricted set $\bigsCecR$. We then need to find an optimal recourse solution $\varRecoOpSndSc$ under that scenario. However, unless we cause a failure in this new solution (Constraint \eqref{eq:AtLeastOne}) we will not succeed at finding a scenario that maximally decreases the value of $\objSingleROr$ in Model \eqref{SingleROmodel}. We repeat the same procedure until we realize that we can no longer generate a failure affecting the last recourse solution added to MasterSecond without violating the failure vertex budget $\vBudget$ and arc failure budget $\aBudget$ implicit in Constraint \eqref{basicGamma}. Therefore, when MasterSecond becomes infeasible, the worst-case scenario must be the one  that led to the optimal recourse solution with the fewest number of pairs from the first stage.

Algorithm \ref{Alg:BasicCovering} generates new failure scenarios $\cdteOne \in \bigsCecR$ until \ref{MasterBasic} becomes infeasible, at which point the worst-case scenario, $\oneCeStar$, and its associated optimal recourse solution value, i.e., the optimal objective value of the second-stage problem, $\objTwoStageStar$, have already been found. Algorithm  \ref{Alg:BasicCovering} starts with an empty subset of scenarios $\bigsCecR$, and with an upper bound on the objective value of the second stage $\UBCovering = \objSingleROr$. \ref{MasterBasic} is solved for the first time to obtain a failure scenario $\cdteOne$, which is then added to $\bigsCecR$. At Step \ref{Alg:SolveRecourse}, the recourse formulation $\RecoPcdte$ is solved and an optimal recourse solution $\varRecoOpSndS$ is obtained with objective value $\RecoObjVarcdte$. The recourse solution is then used to create Constraint \eqref{eq:AtLeastOne} in \ref{MasterBasic}. If the objective value of the recourse solution $\RecoObjVarcdte$ is less than the current upper bound $\UBCovering$, then the upper bound is updated, along with the incumbent failure scenario $\incumCe$. At Step \ref{Alg:ResolveMS}, an attempt is made to find a feasible failure scenario in \ref{MasterBasic}. If such a scenario exists, the process is repeated, otherwise the algorithm returns an optimal solution ($\optCe, \objTwoStageStar$).

\begin{algorithm}[tbp]
	\algorithmicrequire { A recourse policy $\policy \in \policySet$, first-stage solution $\varFirstc \in \setFirst$ and current KEP robust objective $\objSingleROr$}\\
	\algorithmicensure { Optimal recovery plan value $\objTwoStageStar$ and worst-case scenario $\optCe \in \bigsCec$}
    \begin{algorithmic}[1]
    	\STATE  $\bigsCecR = \emptyset$; $\UBCovering \gets \objSingleROr$
        \STATE Solve \ref{MasterBasic} with Constraint \eqref{eq:AtLeastOneX} to obtain scenario $\cdteOne$  \\ 
        \WHILE{MasterSecond$(\varFirstc)$ is feasible}
        \STATE $\bigsCecR \gets \bigsCecR \cup \{\cdteOne\}$
    	\STATE  \label{Alg:SolveRecourse} Solve $\RecoPcdte$ to obtain objective value $\RecoObjVarcdte$ and recourse solution $\varRecoOpSndSc$; create Constraint \eqref{eq:AtLeastOne}\\
    	\IF {$\RecoObjVarcdte < \UBCovering$} 
            \STATE \label{Alg:UpdateUB} $\UBCovering  \gets \RecoObjVarcdte$, $\incumCe \gets \cdteOne$
        \ENDIF
        \STATE \label{Alg:ResolveMS} Attempt to solve MasterSecond$(\varFirstc)$ to get a new candidate scenario $\cdteOne$\\
        \ENDWHILE
    	\STATE $\optCe \gets \incumCe$, $\objTwoStageStar \gets \UBCovering$
        \RETURN $\optCe$, $\objTwoStageStar$
    	\caption{Solving the second-stage problem (Model \eqref{ROmodel}, Step \ref{alg:RO-2ndstage} in Algorithm \ref{Alg:ScenarioGeneration}) using \ref{MasterBasic}}
    	\label{Alg:BasicCovering}
	\end{algorithmic}
\end{algorithm}

To prove the validity of Algorithm \ref{Alg:BasicCovering}, we state the following:
\begin{proposition}
\label{prop:BCproof}
Algorithm \ref{Alg:BasicCovering} returns the optimal objective value of the second stage $\objTwoStageStar$ and worst-case scenario $\oneCeStar \in \bigsCec$ for a first-stage decision $\varFirstc \in \setFirst$.
\end{proposition}

\proof In the first part of the proof, we show that any optimal objective value of the recourse problem $\RecoObjVarcdte$, regardless of the scenario $\cdteOne \in \bigsCec$ is an upper bound on the optimal objective value of the second stage for a given $\varFirstc \in \setFirst$. In the second part, we show that the second-stage problem can be decomposed into an optimality-seeking non-linear binary program which can be linearized by the iterative addition of recourse solutions. We then show that constraints in the binary program have a one-to-one correspondence with those in \ref{MasterBasic}.

Part I. 
Observe that given a candidate solution $\varFirstc \in \setFirst$ and policy $\policy \in \policySet$
\begin{align*}
\max_{\varReco \in \setRecocdte} \objROc &\ge \min_{\oneCe \in \bigsCec} \max_{\varReco \in \setReco} \objROc &  \cdteOne &\in \bigsCec
\\
\sum_{\match \in \cychSetSndSc: \varRecoOpSndSc_{\match} = 1} \wTwoStageGralc &\ge \objTwoStageStar & \cdteOne &\in \bigsCec; \varRecoOpSndSc 
\end{align*}
That is, the optimal objective value of a recourse solution, $\sum_{\match \in \cychSetSndSc: \varRecoOpSndSc_{\match} = 1} \wTwoStageGralc  = \RecoObjVarcdte$, regardless of the scenario $\cdteOne \in \bigsCec$ is at least as big as the smallest value that $\objTwoStageStar$ can reach.

Part II.
In what follows $\mathbbm{1}_{\match,\oneCe}$ is an indicating variable that takes on value one if cycle/chain $\match \in \cychSetSndS$ fails under scenario $\oneCe \in \bigsCec$. Then, we define $\alpha_{\match}$ as a binary decision variable for a cycle/chain $\match \in \cycleSetSndTrc \cup \chainSetSndTrc$ in the transitory graph that takes on value one if such cycle/chain fails, and zero otherwise. We define $\setRecoTrcdte$ as the set of binary vectors representing a recourse solution in the transitory graph and denote $\varRecoc$ to be  a feasible recourse solution, as opposed to its decision vector counterpart $\varReco$. In Appendix \ref{Alg:CyChrecourse} we present a procedure to find $\cycleSetSndTrc$ and $\chainSetSndTrc$ given a first-stage solution $\varFirstc \in \setFirst$. Note that since all feasible chains from length 1 to $\chainCap$ are found, the shortening of a chain when a failure occurs is represented by some $\alpha_{\match}$ taking on value zero and another one taking on value one.
Then, the second-stage problem (SSP) can be reformulated as follows:
\begin{subequations}
\begin{align}
    &\min_{\oneCe \in \bigsCec} \max_{\varReco \in \setRecoc} \sum_{\match \in \cychSetSndS} \wTwoStageGralc\varReco_{\match} \tag{SSP} \label{eq:SSP}\\
    \text{\ref{eq:SSP}} =\min_{\oneCe \in \bigsCecdte} \quad &\objTwoStage \label{refSndSvOne}\\
	& \objTwoStage \ge  \max_{\varReco \in \setRecoc} \sum_{\match \in \cychSetSndS} \left(\wTwoStageGralc - \wTwoStageGralc \mathbbm{1}_{\match,\oneCe} \right) \varReco_{\match} \\
	\text{\ref{eq:SSP}} =\min_{\oneCe \in \bigsCecdte} \quad &\objTwoStage \label{refSndSvTwo}\\
	& \objTwoStage \ge  \sum_{\match \in \cychSetSndS: \varRecoc_{\match} = 1} \left(\wTwoStageGralc - \wTwoStageGralc \mathbbm{1}_{\match,\oneCe} \right) \varRecoc_{\match} & \varRecoc \in \setRecoc
 \end{align} 
    \setcounter{storesubequations}{\value{equation}}%
\end{subequations}

Note that all feasible failure scenarios $\oneCe \in \bigsCec$ exist as a combination of vertex and arc failures in cycles/chains of the transitory graph, i.e., $\match \in \cycleSetSndTrc \cup \chainSetSndTrc$. Therefore, we can model all possible values of the indicating variable $\mathbbm{1}_{\match,\oneCe}$ by replacing it with a decision variable $\alpha_\match$ that takes on value one if cycle/chain $ \match \in \cycleSetSndTrc \cup \chainSetSndTrc$ fails and zero otherwise. The idea is to link this new decision variable to a set of linear constraints such that whenever a vertex/arc in a cycle/chain fails so does the latter. To this end, we introduce the binary vectors $(\smallsCeVx, \smallsCeArc) \in \bigsCec$ as defined in Section \ref{subsec:uncertaintyset}. Thus, we include the minimization over the scenarios within the constraint set as follows:
\addtocounter{equation}{-1}
\begin{subequations}\setcounter{equation}{\value{storesubequations}}
\begin{align}
	\text{\ref{eq:SSP}} =\min \quad &\objTwoStage \label{objMPOpt}\\
	& \objTwoStage \ge  \sum_{\match \in \cycleSetSndTrc \cup \chainSetSndTrc: \varRecoc_{\match} = 1} \left(\wTwoStageGralc - \wTwoStageGralc \alpha_{\match} \right) \varRecoc_{\match}  & \varRecoc \in \setRecoTrcdte \label{eq:solsTwo}\\
	& \alpha_{\match} \le \sum_{\vertexa \in V(\match)} \smallsCeVx_{\vertexa} + \sum_{\parc \in A(\match)} \smallsCeArc_{\narc}
    &\match \in \cycleSetSndTrc \cup \chainSetSndTrc \label{eq:Inter}\\
    &\oneCe \in \bigsCec \label{eq:budget}\\
    &\alpha_\match \in \{0,1\} & \match \in \cycleSetSndTrc \cup \chainSetSndTrc \label{eq:Bound}
\end{align}
    \setcounter{storesubequations}{\value{equation}}%
\end{subequations}

Therefore, Constraints \eqref{eq:Inter} mandate that when a cycle/chain fails, $\alpha_\match$ takes on value one due to the minimization objective. Otherwise, $\alpha_\match$ becomes zero and thus the weight of such cycle/chain is considered in Constraints \eqref{eq:solsTwo}. Observe that once $\alpha_\match$ is set to a value (either zero or one), Model \eqref{objMPOpt} is equivalent to solving Models \eqref{refSndSvOne} or \eqref{refSndSvTwo} when  $\oneCe = \cdteOne$ and the indicating variables $\mathbbm{1}_{\match,\oneCe}$ take their values accordingly for all $\match \in \cychSetSndS$. Lastly, Constraints \eqref{eq:budget} require that the number of vertices and arcs that fail is not exceeded. Since there can exist an exponential number of Constraints \eqref{eq:solsTwo}, a natural way to solve Model \eqref{objMPOpt} is by assuming that the first recourse solution is $\varFirstc$, which can be obtained in the transitory graph.  Then, use $\varRecoc = \varFirstc$ to create the first constraint of Constraints \eqref{eq:solsTwo}. Afterwards, a feasible solution $\tilde{\alpha}$ mapping to some $\cdteOne \in \bigsCec$ can be found, allowing us to solve a recourse problem with optimal solution $\varRecoOpSndSc$. Then, a new Constraint \eqref{eq:solsTwo} can be created using $\varRecoOpSndSc$ and a new failure scenario $\cdteOne$ can be found. This process continues until $\objTwoStage = \RecoObjVarcdte$. Thus, a formulation for the second-stage problem can be defined by the following second-stage formulation (\ref{objMPOptQ}):
\addtocounter{equation}{-1}
\begin{subequations}\setcounter{equation}{\value{storesubequations}}
\begin{align}
	\SSPref = \min \quad &\objTwoStage \tag{SSF} \label{objMPOptQ}\\
    & \objTwoStage \ge  \objSingleROr - \sum_{\match \in \cycleSetSndTrc \cup \chainSetSndTrc: \varFirstc_{\match} = 1}  \wTwoStageGralc \alpha_{\match} \varFirstc_{\match}  & \label{eq:solsTwoX}\\
	& \objTwoStage \ge  \RecoObjVarcdte - \sum_{\match \in \cycleSetSndSc \cup \chainSetSndSc: \varRecoOpSndSc_{\match} = 1}  \wTwoStageGralc \alpha_{\match} \varRecoOpSndSc_{\match}   &\cdteOne \in \bigsCecR; \varRecoOpSndSc  \label{eq:solsTwoQ}\\
	& \alpha_{\match} \le \sum_{\vertexa \in V(\match)} \smallsCeVx_{\vertexa} + \sum_{\parc \in A(\match)} \smallsCeArc_{\narc} 
    &\match \in \cycleSetSndTrc \cup \chainSetSndTrc \label{eq:InterQ}\\
    &\oneCe \in \bigsCec \label{eq:budgetQ}\\
    &\alpha_\match \in \{0,1\} & \match \in \cycleSetSndTrc \cup \chainSetSndTrc \label{eq:BoundQ}
\end{align} 
\end{subequations}

Since $\objSingleROr \ge \RecoObjVarcdte \ge \objTwoStage$, for $\objTwoStage$ to be as small as possible, at least one cycle or chain---and thus, at least one vertex or arc---must fail for some cycle/chain in \emph{every} matching associated with an optimal recourse solution (Constraints \eqref{eq:solsTwoX} and \eqref{eq:solsTwoQ}). If it was not for the failure budget limiting the maximum number of failed vertices and arcs (Constraint \eqref{eq:budgetQ}), $\objTwoStage$ could reach zero. Note that if Model \eqref{objMPOptQ} is unable to cause the failure of a new optimal recourse solution $\varRecoOpSndSc$, it is because Constraint \eqref{eq:budgetQ} would be violated. Thus, the worst-case scenario $\optCe$ is found when there exists a Constraint \eqref{eq:solsTwoQ} associated to a recourse solution $\varRecoOpSndSc$ with non-failed cycles/chains or equivalently when $\objTwoStage = \RecoObjVarcdte$ for some $\cdteOne \in \bigsCecR$. Thus, there is a one-to-one correspondence between (i) Constraint \eqref{eq:AtLeastOneX} in Model \eqref{objMPOptQ} and Constraint \eqref{eq:solsTwoX} in \ref{MasterBasic}, and (ii) Constraints \eqref{eq:solsTwoQ} in Model \eqref{objMPOptQ} and Constraints \eqref{eq:AtLeastOne} in \ref{MasterBasic}. Therefore, the optimal value to the second stage is the smallest value of $\RecoObjVarcdte$ among all scenarios $\cdteOne \in \bigsCecR$ before \ref{MasterBasic} becomes infeasible. \hfill$\square$

It is worth noting that a cycle-chain-based formulation similar to \ref{objMPOptQ} was proposed by \citet{Blom2021} for the full-recourse policy under homogeneous failure, using the set of feasible cycles and chains in what we refer to as the first-stage compatibility graph. However, with our two-stage optimization formulation, Model \eqref{objMPOptQ} can address multiple policies under homogeneous and non-homogeneous failure. 

\subsubsection{Strengthening constraints in \ref{MasterBasic}}
\label{StrengthMS}

We now derive a novel family of non-valid inequalities to narrow the search of the worst-case scenario in \ref{MasterBasic}. Unlike valid inequalities, the so-called \textit{non-valid} inequalities cut off feasible solutions \citep{Atamturk2000, Hooker1994}, and therefore are invalid in the standard sense, though optimal solutions are preserved.

Observe that every time the value of $\UBCovering$ is updated in Algorithm \ref{Alg:BasicCovering}, we know from Part II of Preposition \ref{prop:BCproof} that $\UBCovering$ could be the optimal objective value of the second-stage problem and if so, \ref{objMPOptQ} would need to ``cause" the failure of the previously added constraints, i.e., Constraints \eqref{eq:solsTwoX} and \eqref{eq:solsTwoQ}, in such a way that $\objTwoStage$ is at most $\UBCovering$. Thus, there exists a minimum number of cycles/chains that fail in Constraints \eqref{eq:solsTwoX} and \eqref{eq:solsTwoQ} to achieve this goal, and because those constraints are also represented in the feasibility-seeking master problem, we aim to solve the feasibility-seeking master problem based on our knowledge of the optimality-seeking master problem.  
Thus, we strengthen the right-hand side of Constraints \eqref{eq:AtLeastOneX} and \eqref{eq:AtLeastOne} by updating the minimum number of vertices/arcs that should fail in each of those constraints whenever a smaller value of $\UBCovering$ is found at Step \ref{Alg:UpdateUB} of Algorithm \ref{Alg:BasicCovering}.

For a recourse solution $\varRecoOpSndSc$ and assuming that the failure of a vertex/arc completely causes the failure of a cycle or chain, we sort cycle and chain weights  $\wTwoStagec$ with $\match \in \cychSetSndSc \ \forall \varRecoOpSndSc_\match = 1$ in non-increasing order so that $\wTwoStagec_1 \ge \wTwoStagec_2 ... \ge \wTwoStagec_{\mid \cychSetSndSc \mid}$. Thus, we can now state the following:

\begin{proposition}
\label{prop:RHSi}
Constraints \eqref{eq:AtLeastOneX} and \eqref{eq:AtLeastOne} with right-hand side value $t$ is a non-valid inequality such that $t$ is the smallest index for which the following condition is true: $\RecoObjVarcdte - \sum_{t = 1}^{\mid \cychSetSndSc \mid} \wTwoStagec_t < \UBCovering$.
\end{proposition}
\proof 
Observe that unless optimal, $\objTwoStage < \UBCovering \forall \cdteOne \in \bigsCecR$ due to Proposition \ref{prop:BCproof} and Step \ref{Alg:UpdateUB} of Algorithm \ref{Alg:BasicCovering}. Thus, some cycles/chains must fail in each of the analogous constraints of \eqref{eq:AtLeastOneX} and \eqref{eq:AtLeastOne} in Model \eqref{objMPOptQ}, i.e., Constraints \eqref{eq:solsTwoX} and \eqref{eq:solsTwoQ}. Also, observe that Constraints \eqref{eq:InterQ} imply that whenever at least one vertex/arc fails, so does its associated cycle/chain. Therefore, finding the minimum number of cycles/chains that should fail in Constraints \eqref{eq:solsTwoX} and \eqref{eq:solsTwoQ} to satisfy $\objTwoStage < \UBCovering$ implies finding the minimum number of vertices/arcs that should fail in Constraints \eqref{eq:AtLeastOneX} and \eqref{eq:AtLeastOne}. Thus, sorting the cycle/chain weight of every cycle/chain in non-increasing order, and then subtracting it in that order from $\RecoObjVarcdte$ until the subtraction is strictly less than $\UBCovering$, yields a valid lower bound on the number of cycles/chains that must fail.
\hfill$\square$

\subsection{An expanded feasibility-seeking decomposition}
\label{sec:ExtFeasForm}
Next, we introduce the second novel decomposition for the second-stage problem. This decomposition ``expands" the recourse solutions to allow failed cycles and chains with the goal of generating failure scenarios that speed up the convergence of our master problem to infeasibility and thus, to optimality. To this end, we make use of the transitory graph to find optimal recourse solutions as opposed to the second-stage compatibility graph.

\subsubsection{The recourse problem on the transitory graph}
So far, we have solved the recourse problem using the second-stage compatibility graph $\digraphSndSc$. However, an optimal solution to the recourse problem $\RecoPcdte$ can also be found in the transitory graph $\digraphSndTrc$ by fitting failed cycles/chains into the solution given by $\RecoPcdte$. Fitting failed cycles/chains into recourse solutions with non-failed cycles/chains was considered by \cite{Blom2021} to lift the constraints in their optimality-seeking second-stage formulation. Such formulation is referred to as ``the attack generation problem'' in their work. Our approach borrows a modified objective function to find recourse solutions in \ref{RecoP:Objexp}. Unlike their work, our master problems generalize to homogeneous/non-homogeneous failure and multiple recourse policies.

Although an optimal recourse solution in the transitory graph has the same objective value as an optimal recourse solution in the second-stage compatibility graph (failed cycles/chains do not increase the recourse objective value), as we will show, it can prevent some dominated scenarios from being explored, and therefore help reduce the number of times the recourse problem is resolved in Algorithm \ref{Alg:BasicCovering}. A failure scenario $\oneCePrime \in \bigsCec$ dominates another failure scenario $\oneCeNoPrime \in \bigsCec$ if $\varCov^{\star}_{R}(\varFirst,\oneCePrime) \le \varCov^{\star}_{R}(\varFirst, \oneCeNoPrime)$. We define to $\RecoPexpc$ as the recourse problem solved in the transitory graph whose optimal recourse solutions are also optimal with respect to $\RecoPcdte$ but may allow some failed cycles/chains in the solution. We refer the reader to Appendix \ref{sec:RecourseFormls} for proof that optimal recourse solutions to $\RecoPexpc$ are also optimal to $\RecoPcdte$, which is adapted from \citet{Blom2021}. Thus,  we let $\VarRecoOpExp \in \setRecoTrcdte$ be an optimal recourse solution to $\RecoPexpc$ in the transitory graph that is also optimal to $\RecoPcdte$ under scenario $\cdteOne \in \bigsCecR$.

\subsubsection{The master problem: MasterTransitory}
In the basic reformulation, we assumed the recourse solutions correspond to matchings with non-failed components only. In our expanded formulation, we assume that recourse solutions can be expanded to fit some failed cycles/chains by solving the recourse problem in the transitory graph $\digraphSndTrc$. We let $\VarRecoOpExpcdte \subseteq \VarRecoOpExp$ be the subset of the optimal recourse solution $\VarRecoOpExp$ that has no failed cycles/chains and thus corresponds to a feasible solution  in $\setRecocdte$. Therefore, the expanded feasibility-seeking reformulation, which we call MasterTransitory, is expressed as follows:
\begin{subequations}
\label{mo:extendedfeas}
    \begin{align}
        \tag{MT} \text{MasterTransitory}(\varFirstc) \text{: }  &&& \qquad \qquad \text{Find } \oneCe \label{MasterExp}\\
        && \sum_{\match \in \cycleSet \cup \chainSet: \varFirstc_{\match} = 1} &\left(\sum_{\vertexa \in V(\match)} \smallsCeVx_{\vertexa} +  \sum_{\parc \in A(\match)} \smallsCeArc_{\narc} \right) \ge 1   &&   \label{eq:AtLeastOneXExp}\\
         && \sum_{\match \in \cychSetSndScdte: \varRecoOpSndSc_{\match} = 1} &\left(\sum_{\vertexa \in V(\match)} \smallsCeVx_{\vertexa} +  \sum_{\parc \in A(\match)} \smallsCeArc_{\narc} \right) \ge 1 &&  \cdteOne \in \bigsCecR; \VarRecoOpExpcdte   \label{eq:AtLeastOneExp}\\
        && \sum_{\match \in \cycleSetSndTrc \cup \chainSetSndTrc: \varRecoOpSndSc_{\match} = 1} &\left(\sum_{\vertexa \in V(\match)} \smallsCeVx_{\vertexa}+ \sum_{\parc \in A(\match)} \smallsCeArc_{\narc} \right) \ge \rhsfeasite  && \ \cdteOne \in \bigsCecR; \VarRecoOpExp \label{eq:AtLeastOneRHSExp}\\
        &&& \qquad \qquad \oneCe \in \bigsCe&& \label{GammaExp}
    \end{align}
\end{subequations}

Constraints \eqref{eq:AtLeastOneXExp} and \eqref{eq:AtLeastOneExp} are equivalent to Constraints \eqref{eq:AtLeastOneX} and \eqref{eq:AtLeastOne}. Constraints \eqref{eq:AtLeastOneRHSExp} require that a combination of at least $\rhsfeasite$ vertices and arcs fail in an expanded recourse solution $\VarRecoOpExp \in \setRecoExpc$.  Proposition \ref{prop:RHSExtdi} defines the value of $\rhsfeasite$. The goal of \ref{MasterExp} is to enforce the failure of two different recourse solutions with identical objective values under scenario $\cdteOne \in \bigsCecR$, i.e. a recourse solution found in the second-stage compatibility graph $\digraphSndSc$ and a solution found in the transitory graph $\digraphSndTrc$. 

For a recourse solution $\VarRecoOpExp \in \setRecoExpc$ again assuming that the failure of a vertex/arc completely causes the failure of a cycle or chain, we sort cycle and chain weights  $\wTwoStagec$  $\forall \varReco_\match = 1$ with $\match \in \cycleSetSndTrc \cup \chainSetSndTrc$, in non-increasing order so that $\wTwoStagec_1 \ge \wTwoStagec_2 ... \ge \wTwoStagec_{\mid \cycleSetSndTrc \cup \chainSetSndTrc \mid}$. Thus, we can now state:

\begin{proposition}
\label{prop:RHSExtdi}
$\rhsfeasite = t$ is a non-valid inequality such that $t$ is the smallest index for which the following condition is true $\RecoObjVarcdte - \sum_{t = 1}^{\mid \cycleSetSndTrc \cup \chainSetSndTrc \mid} \wTwoStagec_t < \UBCovering$.
\end{proposition}
\proof
Follows by the same arguments given for Proposition 
\ref{prop:RHSi}.\hfill$\square$

An example demonstrating that \ref{MasterExp} can yield failure scenarios that dominate those in \ref{MasterBasic} can be found in Appendix \ref{Example}.

\subsubsection{Dominating scenarios}
\label{DominatingScenarios}

Next, we derive another novel family of non-valid inequalities based on dominating scenarios. We say that a failure scenario $\oneCePrime \in \bigsCec$ dominates another failure scenario $\oneCeNoPrime \in \bigsCec$ when $\RecoObjVariteExPrime \le \RecoObjVariteEx$. Let $ I(\oneCePrime)$  and $I(\oneCeNoPrime)$ be the number of failed vertices and arcs under their corresponding scenario. Moreover, let $\cyclechainSetPrimeX \subseteq \cycleSetSndTrc \cup \chainSetSndTrc$ and $\cyclechainSetNoPX \subseteq \cycleSetSndTrc \cup \chainSetSndTrc$ be the set of feasible cycles and chains that fail in the transitory graph $\digraphSndTrc$ under scenarios $\oneCePrime$ and $\oneCeNoPrime$, respectively. We then have the following result:

\begin{proposition}
\label{prop:dominance}
If $\cyclechainSetNoPX \subseteq \cyclechainSetPrimeX$, $\oneCePrime$ dominates $\oneCeNoPrime$ and the following dominance inequality is valid for the second-stage problem:
\begin{align}
    \label{eq:dominance}
       \sum_{\vertexa \in \vertexSet: \oneCeNoPrime_{\vertexa} = 1} \smallsCeVx_{\vertexa} + \sum_{\parc \in \arcSet: \oneCeNoPrime_{\narc} = 1} \smallsCeArc_{\narc} \le I(\oneCeNoPrime) \left(I(\oneCePrime) - 
            \sum_{\vertexa: \oneCePrime_{\vertexa} = 1} \smallsCeVx_{\vertexa} - \sum_{\parc: \oneCePrime_{\narc} = 1} \smallsCeArc_{\narc}
       \right)
\end{align}
\end{proposition}

\proof By definition, $\cyclechainSetNoPX \subseteq \cyclechainSetPrimeX$, thus, all cycles and chains that fail under $\oneCeNoPrime$ also fail under $\oneCePrime$, which means that $\RecoObjVariteExPrime \le \RecoObjVariteEx$. Therefore, it follows that $\oneCePrime$ dominates scenario $\oneCeNoPrime$.\hfill$\square$

Finding all $\oneCeNoPrime, \oneCePrime \in \bigsCec$ satisfying Proposition \ref{prop:dominance} is not straightforward and there may be too many constraints of type \eqref{eq:dominance} to feed into our master problem formulations. Therefore, we explore two alternatives to separate these dominating-scenario cuts. The first consists of identifying a subset of scenarios that satisfy constraint \eqref{eq:dominance} a priori, which we call adjacent-failure separation. The second consists of attempting to solve \ref{MasterBasic} and \ref{MasterExp} via a heuristic and then ``discovering'' dominating scenarios on the fly, which we call single-vertex-arc separation. We refer to this approach as a heuristic because it may not always find a feasible failure scenario satisfying the constraints in \ref{MasterBasic} or \ref{MasterExp}.

\paragraph{Adjacent-failure separation} 
A vertex failure in the transitory graph $\digraphSndTr$ is equivalent to the failure of that vertex and the failure of any arc adjacent to that vertex, since in both cases the failed cycles and chains are the same. Thus, for every vertex $\vertex \in \vertexSet$, we can build a scenario where the only non-zero value in vector $\oneCeNoPrime^{\prime\text{v}}$ is $\oneCeNoPrime^{\prime\text{v}}_{\vertex} = 1$, which dominates the scenario where all arcs either leave $\vertex$, i.e., $\oneCeNoPrime^{\text{a}}_{\narc} = 1$  or point towards it, i.e., $\oneCeNoPrime^{\text{a}}_{\narcrev} = 1$. Note that the following constraints satisfy Proposition \ref{prop:dominance}:
\begin{align}
    \label{eq:adjacent}
    \sum_{\parc \in \arcSet} \smallsCeArc_{\narc} + \sum_{\parcrev \in \arcSet} \smallsCeArc_{\narcrev} \le \aBudget \left(1 -  \smallsCeVx_{\vertexa} \right)&& \vertexa \in \vertexSet
\end{align}
That is, if a vertex fails, the arcs adjacent to it get disconnected from the graph. Therefore, the objective value of a recourse problem with a failure scenario where an arc adjacent to a failed vertex also fails is equivalent to the objective value obtained when we consider the scenario with only the vertex failure. Thus, Constraints \eqref{eq:adjacent} can be added to \ref{MasterBasic} and  \ref{MasterExp} before the start of Algorithm \ref{Alg:BasicCovering}.

\paragraph{Single-vertex-arc separation} 
\label{vxtarcsep}
Suppose we have a set of proposed failed arcs $A(\oneCeNoPrime^{\text{a}})$ such that $\smallsCeArc_{\narc} = 1 \ \forall {\parc} \in A(\oneCeNoPrime^{\text{a}})$, and a set of proposed failed vertices, $V(\oneCeNoPrime^{\text{v}})$, such that $\oneCeNoPrime^{\text{v}}_{\vertexa} = 1 \ \forall \vertexa \in V(\oneCeNoPrime^{\text{v}})$. The second-stage compatibility graph corresponds to $\digraphSndSpr = (\vertexSetSndBe \setminus V(\oneCeNoPrime^{\text{v}}), \arcSetSndBe \setminus A(\oneCeNoPrime^{\text{a}}))$. Then, the following two cases also satisfy Proposition \ref{prop:dominance}:

\begin{enumerate}
    \item \label{CaseOne} Suppose there is a candidate pair $\bar{\vertexb} \in \vertexSetSndBe \setminus V(\oneCeNoPrime^{\text{v}})$ and $\cyclechainSetXvtx$ is the set of feasible cycles and chains in $\digraphSndSpr$ that include vertex $\bar{\vertexb}$. Then, if  $\cyclechainSetXvtx = \emptyset$, scenario $\oneCeNoPrime$ dominates scenario $\oneCeNoPrime \cup \{\bar{\vertexb}\}$ and thus for $\UBCovering$ to decrease, another vertex $\bar{\vertexb}$ should be proposed to fail.
    
    \item \label{CaseTwo} Suppose there is a candidate arc $\bar{a} \in \arcSetSndBe \setminus A(\oneCeNoPrime^{\text{a}})$ and $\cyclechainSetXarc$ is the set of feasible cycles and chains in $\digraphSndSpr$ that include arc $\bar{a}$. Then, if  $\cyclechainSetXarc = \emptyset$, scenario $\oneCeNoPrime$ dominates scenario $\oneCeNoPrime \cup \{\bar{a}\}$ and thus another arc $\bar{a}$ should be proposed to fail if $\UBCovering$ is to be decreased.
\end{enumerate}

 We perform the single-vertex-arc separation within a heuristic while attempting to find a feasible failure scenario satisfying the constraints in \ref{MasterBasic} and  \ref{MasterExp}, as detailed in the next section.
\section{Solution algorithms for the second-stage problem}
\label{sec:HSAs}

This section presents solution algorithms to solve the second stage.
We refer to \textit{hybrid} solution algorithms as the combination of a linear optimization solver and a heuristic attempting to solve  the second-stage problem. Our goal is to specialize the steps of Algorithm \ref{Alg:BasicCovering} to improve the efficiency of our final solution approaches.

\subsection{A feasibility-based solution algorithm for MasterBasic}
\label{FBSA-MB}

Algorithm \ref{Alg:BasicHybrid} is the first of our two feasibility-based solution algorithms (FBSAs), which has \ref{MasterBasic} as master problem, and it is referred to as FBSA\_MB. FBSA\_MB requires a policy $\policy \in \policySet$ and a first-stage solution $\varFirstc$ found at Step \ref{alg:R0-solve} of Algorithm \ref{Alg:ScenarioGeneration}. We refer to $\textbf{ToMng}(\cdot)$ as a function that ``extracts'' the matching from a recourse solution. This matching is then added to a set $\RecoSet$ containing all matchings associated with the recourse solutions found up to some iteration. The procedure \textbf{Heuristic}($\RecoSet$)---whose algorithmic details we present shortly---attempts to find a failure scenario satisfying the constraints in \ref{MasterBasic}. The heuristic returns a tuple $\HeuAns$ with two outputs: a candidate failure scenario $\incumCe$ and a boolean variable $\cover$ that indicates whether $\incumCe$ is feasible to \ref{MasterBasic}. If $\cover = \text{\textbf{true}}$, MasterSecond$(\varFirstc)$ is feasible. 

Algorithm \ref{Alg:BasicHybrid} starts with an empty subset of scenarios $\bigsCecR$, and with an upper bound on the objective value of the second stage $\UBCovering = \objSingleROr$. At Step \ref{HeuTrue}, because there is only one matching in $\RecoSet$, it is trivial to heuristically choose a failure scenario that satisfies Constraint \eqref{eq:AtLeastOneX} and thus the result of the boolean variable $\cover$ is \textbf{true}. This failure scenario is then added to $\bigsCecR$. At Step \ref{ColGenAlg3}, $\textbf{ColGen}(\RecoPcdte)$ finds an upper bound $\ZRecoColUB$ on the optimal value of the recourse problem $\RecoPcdte$ through Column Generation (CG). CG is a technique for solving linear programs with a large number of variables or columns \citep{Barnhart1998}. For $\textbf{ColGen}(\RecoPcdte)$, the columns correspond to cycles and chains in the set $\cycleSetSndSc \cup \chainSetSndSc$. CG considers two problems: a master problem and subproblem(s). In Algorithm \ref{Alg:BasicHybrid}, the master problem is the linear programming relaxation of \ref{RecoP:Obj}. At the start, the master problem starts with no columns. The subproblem consists of finding columns that have the potential to increase the objective value of the master problem, i.e., columns with positive reduced cost (we refer the reader to \citep{Riascos2020} for further details). We add up to 10 cycle columns and up to 10 chain columns at every iteration of the CG algorithm, unless $\chainCap \ge 4$, case in which we add up to 5 chain columns. Moreover, if $\chainCap \ge 4$, we split the set of chains $\chainSetSndSc$ into two equally sized sets. After having unsuccessfully searched for chain columns with positive reduced cost in the first set, we then continue the search for such chain columns over all chains in $\chainSetSndSc$.

Once no more columns with positive reduced cost are found, and thus the optimality of the CG master problem is proven, an upper bound $\ZRecoColUB$ on the optimal value of $\RecoPcdte$ is returned. We then take the decision variables from the optimal base of the master problem and turn them into binary decision variables to obtain a feasible recourse solution $\tilde{\varReco}$ with objective value $\ZRecoCol$. If $\ZRecoCol$ equals $\ZRecoColUB$, then $\tilde{\varReco}$ is optimal to formulation $\RecoPcdte$. Thus, at Step \ref{SolCGisOPT_Alg3}, $\tilde{\varReco}$ becomes the optimal recourse solution $\varRecoOpSndS$ and $\RecoObjVarcdte$ is updated. Otherwise, we incur in the cost of solving \ref{RecoP:Obj} from scratch as a MIP instance to obtain $\varRecoOpSndS$ and $\RecoObjVarcdte$. The matching associated with $\varRecoOpSndS$ is then added to $\RecoSet$. If the new recourse value $\RecoObjVarcdte$ is smaller than the upper bound on the objective value of the second stage, $\UBCovering$, then we update both  $\UBCovering$ and the incumbent failure scenario $\optCe$ that led to that value. Since a new lower bound has been found, the right-hand side of Constraints \eqref{eq:AtLeastOneX}-\eqref{eq:AtLeastOne} is updated following Proposition \ref{prop:RHSi}. At Step \ref{newIteAlg3}, the heuristic tries to find a new failure scenario feasible to \ref{MasterBasic}+Constraints\eqref{eq:adjacent}. If unsuccessful, \ref{MasterBasic}+Constraints\eqref{eq:adjacent} is solved as a MIP instance. If a feasible failure scenario is found, then the algorithm continues. Otherwise, the algorithm ends and returns the optimal recovery plan ($\objTwoStageStar$, $\optCe$).

\begin{algorithm}[tbp]
	\algorithmicrequire { A recourse policy $\policy \in \policySet$, first-stage solution $\varFirstc \in \setFirst$ and current KEP robust objective $\objSingleROr$} \\
	\algorithmicensure { Optimal recovery plan value $\objTwoStageStar$ and worst-case scenario $\optCe \in \bigsCec$}
    \begin{algorithmic}[1]
      \STATE $\bigsCecR = \emptyset$; $\UBCovering \gets \objSingleROr$
      \STATE Create Constraint \eqref{eq:AtLeastOneX} in \ref{MasterBasic}
      \STATE $\RecoSet \gets  \textbf{ToMng}(\varFirstc)$
	   \STATE \label{HeuTrue} $\HeuAnsTrue \gets$ \text{\textbf{Heuristic}}($\RecoSet$); $\cdteOne \gets \incumCe$ 
    \WHILE{MasterSecond$(\varFirstc)$ is feasible}
        \STATE $\bigsCecR \gets \bigsCecR \cup \{\cdteOne\}$
	    \STATE \label{ColGenAlg3} $\TupleZRecoCol \gets \textbf{ColGen}\left(\RecoPcdte \right)$
	\IF{$\ZRecoCol = \ZRecoColUB$}
	    \STATE   $\RecoObjVarcdte \gets \ZRecoCol$ 
        \STATE \label{SolCGisOPT_Alg3} $\varRecoOpSndSc \gets \tilde{\varReco}$
	\ELSE
	\STATE Solve $\RecoPcdte$ to obtain solution $\varRecoOpSndSc \in \setRecocdte$ with objective value $\RecoObjVarcdte$
	\ENDIF
    \STATE $\RecoSet \gets  \RecoSet \cup \textbf{ToMng}(\varRecoOpSndSc)$
	\IF{$\RecoObjVarcdte < \UBCovering$}
	      \STATE $\UBCovering \gets \RecoObjVarcdte$; 
        \STATE $\optCe \gets \cdteOne$
	    \STATE Update right-hand side of constraints in \ref{MasterBasic}, \ $\forall \cdteOne \in \bigsCecR$ following Proposition \ref{prop:RHSi}
	\ENDIF
    \STATE \label{newIteAlg3} $\HeuAns \gets \text{\textbf{Heuristic}}(\RecoSet)$
    \IF{$\cover = \text{\textbf{true}}$}
        \STATE $\cdteOne \gets \incumCe$
    \ELSE
    \STATE Create Constraint \eqref{eq:AtLeastOne} for $\varRecoOpSndSc \in \setRecocdte$ in \ref{MasterBasic}
    \STATE Attempt to solve \ref{MasterBasic}+Consts.\eqref{eq:adjacent} to get new failure scenario $\cdteOne$
    \ENDIF
    \ENDWHILE
    \STATE $\objTwoStageStar \gets \UBCovering$; Return $\optCe$ and $\objTwoStageStar$
	\caption{FBSA\_MB: A Feasibility-based Solution Algorithm for \ref{MasterBasic}}
	\label{Alg:BasicHybrid}
	\end{algorithmic}
\end{algorithm}

\subsection{A feasibility-based solution algorithm for MasterExp}

Algorithm \ref{Alg:EnhancedHybrid} is our second feasibility-based solution algorithm and has \ref{MasterExp} as the master problem of the second stage. We refer to Algorithm \ref{Alg:EnhancedHybrid} as FBSA\_ME. 
FBSA\_ME requires a recourse policy $\policy \in \policySet$ and a first-stage solution $\varFirstc \in \setFirst$. In addition to the notation defined for FBSA\_MB, we introduce new one. We refer to $\RecoSetexp$ as the set of matchings associated to the optimal recourse solutions in $\RecoPexpc$. The CG algorithm, $\textbf{XColGen}\left(\RecoPexpc \right)$, returns a tuple with three values: $\ZRecoColUBexp$, $\ZRecoColexp$ and $\tilde{\varReco}$. $\ZRecoColUBexp$  corresponds to the optimal objective value returned by the master problem of the CG algorithm; $\ZRecoColexp$ is the feasible solution to the KEP after turning the decision variables in the optimal basis of the CG master problem into binaries. Lastly, $\tilde{\varReco}$ is the feasible recourse solution found by the CG algorithm with objective value $\ZRecoColexp$. For $\textbf{XColGen}\left(\RecoPexpc \right)$, the master problem is the linear programming relaxation of \ref{RecoP:Obj}. The subproblem also splits $\chainSetSndTrc$ in half when $\chainCap \ge 4$ and cycle/chain columns are searched and added in the same way as described in Section \ref{FBSA-MB}.

At the start of Algorithm \ref{Alg:EnhancedHybrid}, \textbf{Heuristic}($\RecoSet$) returns the first failure scenario that is added to $\bigsCecR$. At Step \ref{XColAlg4}, the CG algorithm in $\textbf{XColGen}\left(\RecoPexpc \right)$ is solved. If $\ZRecoColexp = \ZRecoColUBexp$ then a function $\textbf{TrueVal}(\tilde{\varReco})$ computes the objective value of $\tilde{\varReco}$ in terms of the original recourse objective function, i.e., \ref{RecoP:Obj}. In case $\ZRecoColexp$ is lower than $\ZRecoColUBexp$, the original recourse problem $\RecoPcdte$ is solved as a MIP instance. The reason for this decision is that  $\RecoPexpc$ includes a larger set of cycle-and-chain decision variables, possibly leading to scalability issues when the recourse problem is solved as a MIP instance. At Step \ref{UpdateNonFMngSet}, the matching associated with the non-failed cycles and chains in the optimal recourse solution is included in $\RecoSet$. If the upper bound $\UBCovering$ is updated, so is the right-hand side of constraints in \ref{MasterExp}. At Step \ref{newIteAlg4}, the heuristic tries to find a new failure scenario feasible to \ref{MasterBasic}+Constraints\eqref{eq:adjacent}. If unsuccessful, \ref{MasterBasic}+Constraints\eqref{eq:adjacent} is solved as a MIP instance. If a feasible failure scenario is found, then the algorithm continues. Otherwise, the algorithm ends and returns the optimal recovery plan ($\objTwoStageStar$, $\optCe$).

\begin{algorithm}[tbp]
	\algorithmicrequire { A recourse policy $\policy \in \policySet$, a first-stage solution $\varFirstc \in \setFirst$ and current KEP robust objective $\objSingleROr$}\\
	\algorithmicensure { Optimal recovery plan value $\objTwoStageStar$ and worst-case scenario $\oneCe^{\star} \in \bigsCecdte$}
    \begin{algorithmic}[1]
    \STATE $\UBCovering \gets \objSingleROr$; $\bigsCecR = \emptyset$
    \STATE Create Constraint \eqref{eq:AtLeastOneXExp} in \ref{MasterExp}
    \STATE $\RecoSet \gets \textbf{ToMng}(\varFirstc)$
    \STATE  $\HeuAnsTrue \gets$ \text{\textbf{Heuristic}}($\RecoSet$); $\cdteOne \gets \incumCe$
    \WHILE{MasterTransitory$(\varFirstc)$ is feasible}
    \STATE \label{AddNewCAlg4} $\bigsCecR \gets \bigsCecR \cup \{\cdteOne\}$
	\STATE \label{XColAlg4} $\TupleZRecoColexp \gets \textbf{XColGen}\left(\RecoPexpc \right)$
	\IF{$\ZRecoColexp = \ZRecoColUBexp$}
	\STATE $\RecoObjVarcdte \gets \textbf{TrueVal}(\tilde{\varReco})$
    \STATE $\VarRecoOpExp \gets \tilde{\varReco}$
    \STATE $\VarRecoOpExpcdte \subseteq \VarRecoOpExp$, s.t. $\VarRecoOpExpcdte \in \setRecocdte$
	\STATE $\RecoSetexp \gets \RecoSetexp \cup \textbf{ToMng}(\VarRecoOpExp)$
	\ELSE
	\STATE Solve $\RecoPcdte$ to obtain recourse objective value $\RecoObjVarcdte$ and recourse solution $\varRecoOpSndSc$
    \STATE $\VarRecoOpExpcdte = \varRecoOpSndSc$;  $\VarRecoOpExp = \emptyset$
	\ENDIF
    \STATE \label{UpdateNonFMngSet} $\RecoSet \gets \RecoSet \cup \textbf{ToMng}(\VarRecoOpExpcdte)$
	\IF{$\RecoObjVarcdte < \UBCovering$}
	\STATE $\UBCovering \gets \RecoObjVarcdte$ 
    \STATE $\optCe \gets \cdteOne$
	\STATE Update right-hand side of \eqref{eq:AtLeastOneXExp} and \eqref{eq:AtLeastOneExp} following Proposition  \ref{prop:RHSi} and \eqref{eq:AtLeastOneRHSExp} following Proposition \ref{prop:RHSExtdi} \ $\forall \cdteOne \in \bigsCecR$
	\ENDIF \label{EndMTUpdateAlg4}
    \STATE \label{newIteAlg4} $\HeuAns \gets \text{\textbf{Heuristic}}(\RecoSet \cup \RecoSetexp)$
    \IF{$\cover = \text{\textbf{true}}$}
    \STATE $\cdteOne \gets \incumCe$
    \ELSE
    \STATE Create Const.\eqref{eq:AtLeastOneExp}  for $\VarRecoOpExpcdte \in \setRecocdte$ and Const.\eqref{eq:AtLeastOneRHSExp} for $\VarRecoOpExp \in \setRecoTrcdte$ in \ref{MasterExp}
    \STATE Attempt to solve \ref{MasterExp}+Consts.\eqref{eq:adjacent} to get new failure scenario $\cdteOne$ 
    \ENDIF \label{EndnewIteAlg4}
    \ENDWHILE
    \STATE $\objTwoStageStar \gets \UBCovering$; Return $\optCe$ and $\objTwoStageStar$
	\caption{FBSA\_ME: A Feasibility-based Solution Algorithm for \ref{MasterExp}}
	\label{Alg:EnhancedHybrid}
	\end{algorithmic}
\end{algorithm}

\subsection{Hybrid solution algorithms}
The main difference between the hybrid algorithms and the feasibility-based ones is the possibility of transitioning from the feasibility-seeking master problems to the optimality-seeking problem \ref{objMPOptQ} after TR iterations. The goal is to obtain a lower bound on the objective value of the second stage, $\LBCovering$, that can be compared to $\UBCovering$ to prove optimality. We refer to these hybrid algorithms as HSA\_ME (Algorithm \ref{Alg:HybridHSA-ME}) and HSA\_MB. The former has \ref{MasterExp} as master problem, whereas the latter has \ref{MasterBasic}. For both hybrid algorithms, we solve the recourse problem $\RecoPexpc$, but in the case of HSA\_MB, whenever an attempt to find a failure scenario occurs, it is done with respect to \ref{MasterBasic} + Constraints \eqref{eq:adjacent}. The reason for this approach is to accumulate the recourse solutions given by $\RecoPexpc$ and use them to solve \ref{objMPOptQ} when the number of iterations exceeds TR. The HSA\_MB algorithm is the HSA\_ME algorithm (Algorithm \ref{Alg:HybridHSA-ME}) with the following modifications:
\begin{itemize}
    \item Step 2: Create Constraint \eqref{eq:AtLeastOneX} in \ref{MasterBasic}
    \item Step 5: Change MasterTransitory$(\varFirstc)$ to MasterSecond$(\varFirstc)$ in the \texttt{while} loop
    \item Step 21:  Update right-hand side of constraints in \ref{MasterBasic} $\forall \cdteOne \in \bigsCecR$ following Proposition \ref{prop:RHSi}
    \item Step 27:  Create Constraint \eqref{eq:AtLeastOne} for $\VarRecoOpExpcdte \in \setRecocdte$ in \ref{MasterBasic}
    \item Step 28: Attempt to solve \ref{MasterBasic}+Constraints\eqref{eq:adjacent} to get new failure scenario $\cdteOne$
\end{itemize}

\begin{algorithm}[tbp]
	\algorithmicrequire{ Parameter TR indicating iteration from which \ref{objMPOptQ} is solved and those parameters in Algorithm \ref{Alg:EnhancedHybrid}}\\
 	\algorithmicensure { Optimal recovery plan value $\objTwoStageStar$ and worst-case scenario $\oneCe^{\star} \in \bigsCecdte$}
    \begin{algorithmic}[1]
    \STATE $\UBCovering \gets \objSingleROr$; $\LBCovering = 0$; $\bigsCecR = \emptyset$
    \STATE Create Constraint \eqref{eq:AtLeastOneXExp} in \ref{MasterExp}
    \STATE $\RecoSet \gets \textbf{ToMng}(\varFirstc)$
    \STATE  $\HeuAnsTrue \gets$ \text{\textbf{Heuristic}}($\RecoSet$); $\cdteOne \gets \incumCe$
    \WHILE{MasterTransitory$(\varFirstc)$ is feasible \OR $\LBCovering < \UBCovering$}
    \STATE Step \ref{AddNewCAlg4} to \ref{EndMTUpdateAlg4} in Algorithm \ref{Alg:EnhancedHybrid}
    \IF{$\iteTwoStage \ge$ TR}
        \STATE Create a Const.\eqref{eq:solsTwoQ} per matching in $\RecoSetexp$ and solve \ref{objMPOptQ} to obtain $\LBCovering$ and failure scenario $\cdteOne$
        \IF{$\LBCovering = \UBCovering$}
            \STATE $\optCe \gets \cdteOne$
        \ENDIF
    \ELSE
        \STATE Steps \ref{newIteAlg4} to \ref{EndnewIteAlg4} in Algorithm \ref{Alg:EnhancedHybrid}
    \ENDIF
    \ENDWHILE
    \STATE $\objTwoStageStar \gets \UBCovering$; Return $\optCe$ and $\objTwoStageStar$
	\caption{HSA\_ME}
	\label{Alg:HybridHSA-ME}
	\end{algorithmic}
\end{algorithm}

\subsection{Heuristic}
We now discuss the heuristic (Algorithm \ref{Alg:Heuristics}), used in the feasibility-based and hybrid algorithms. Algorithm \ref{Alg:Heuristics} requires a set of matchings $\genmatch$ as input. The heuristic returns a failure scenario $\incumCe$ and a boolean variable \cover. If \cover = \textbf{true}, $\incumCe$ satisfies either \ref{MasterBasic}+Constraints\eqref{eq:adjacent} or \ref{MasterExp}+Constraints\eqref{eq:adjacent}, accordingly. Although the heuristic cannot ensure that $\incumCe$ is always feasible with respect to the feasibility-seeking master problems, it does ensure (i) Constraints \eqref{eq:adjacent} are satisfied and (ii) $\incumCe$ is not dominated according to the two single-vertex-arc separation cases presented in Section \ref{DominatingScenarios}.

We first introduce notation specific to the heuristic, and then detail the steps of the heuristic. A function \textbf{UniqueElms}($\genmatch$) returns the set of unique vertices and arcs among all matchings in $\genmatch$, which are then kept in $\ElSet$.  With a slight abuse of notation, we say that each vertex/arc $\Elem_{n}  \in \ElSet$ with $n = 1,...\lvert \ElSet \rvert$ has an associated boolean variable $\Elsta_{n}$ and a calculated weight $\Elwei_{n}$. The boolean variable $\Elsta_{n}$ becomes \textbf{true} if element $\Elem_{n}$ has been checked, i.e., either it has been proposed for failure in $\incumCe$ or it has been proven that adding $\Elem_{n}$ to the current failure scenario will result in a dominated scenario according to Section \ref{DominatingScenarios}. Function $\text{\textbf{Weight}}(\Elem_{n},\genmatch)$ returns the weight $\Elwei_{n}$ corresponding to the number of times element $\Elem_{n}$ is repeated in the set of matchings $\genmatch$. Moreover, a function $\text{\textbf{IsNDD}}(\Elemstar_{n})$ determines whether element $\Elemstar_{n}$ is a non-directed donor. 

At Step \ref{whileHeu}, while it is true, a list of unique and non-checked elements $\matchHeuSet$ is created as long as the vertex and arc budgets, $\vBudget$ and $\aBudget$ are not exceeded, respectively. If $\matchHeuSet$ turns out to have no elements, the algorithm ends at Step \ref{firstEndHeu}. Otherwise, the elements in $\matchHeuSet$ are sorted in non-increasing order of their weights, and among the ones with the highest weight, an element $\Elemstar_{n}$ is selected randomly. 

At Step \ref{SepHeu}, the single-vertex-arc separation described in Section \ref{vxtarcsep} is performed if (i) $\Elemstar_{n}$ is a pair, and (ii) there is at least one arc or at least two vertices proposed for failure already. The reason we check that $\Elemstar_{n}$ is a pair is because a non-directed donor does not belong to a cycle by definition. Moreover, when a non-directed donor is removed from the transitory graph or second-stage compatibility graph, its removal trivially causes the failure of all chains triggered by it. Thus, a boolean variable \textit{ans} remains \textbf{true} if $\Elemstar_{n}$ is a non-directed donor or a vertex/arc does not satisfy neither Case \ref{CaseOne} nor Case \ref{CaseTwo}. At Step \ref{checked}, $\Elemstar_{n}$ is labeled as checked by setting $\Elemstar_{n} = \text{\textbf{true}}$. Then, if  \textit{ans} = \textbf{true}, it is checked whether $\Elemstar_{n}$ is a vertex or an arc, and the proposed failure scenario $\HatoneCe$ is updated accordingly. If element $\Elemstar_{n}$ is a vertex, at Step \ref{adjArcstHeu} the set of arcs adjacent to $\Elemstar_{n}$ is labeled as checked, so that $\HatoneCe$ satisfies Constraints \eqref{eq:adjacent}.

We say that scenario $\HatoneCe$ covers all matchings in $\genmatch$, if for every matching in that set, the number of vertices/arcs proposed for failure satisfies Proposition \ref{prop:RHSi} (for Constraints \eqref{eq:AtLeastOneX}-\eqref{eq:AtLeastOne} or Constraints \eqref{eq:AtLeastOneXExp}-\eqref{eq:AtLeastOneExp}) and Proposition \ref{prop:RHSExtdi} (for Constraints \eqref{eq:AtLeastOneRHSExp}), accordingly. If $\HatoneCe$ covers all matchings in $\genmatch$, then \cover = \textbf{true}. The heuristic is greedy in the sense that even when \cover = \textbf{true}, it will attempt to use up the vertex and arc budgets as checked at Step \ref{budgetHeu}. If another vertex/arc can still be proposed for failure, then a new iteration in the While loop starts. Thus, the heuristic ends at either Steps \ref{firstEndHeu} or \ref{SndEndHeu}.

\begin{algorithm}[tbp]
	\algorithmicrequire{ A set of matchings $\genmatch$}\\
	\algorithmicensure{ Failure scenario $\HatoneCe$ and a boolean variable \cover\ indicating whether $\HatoneCe$ covers matchings in $\genmatch$}
	
    \begin{algorithmic}[1]
    \STATE \cover $=$ \textbf{false}; $\ElSet \gets \textbf{UniqueElms}(\genmatch)$
        \WHILE{\textit{true}} \label{whileHeu}
        \STATE $\matchHeuSet = \emptyset$
        \FOR{$n = 1,..., \lvert \ElSet \rvert$}
            \IF{$\Elem_{n} \in \arcSet$ and $\Elsta_{n} =$ \textbf{false} and $\sum_{\vertex \in \vertexSet} \HatoneCe_{\vertex}^{\text{v}}  < \vBudget$}
            \STATE $\Elwei_{n} = \text{\textbf{Weight}}(\Elem_{n},\genmatch)$
            ; $\matchHeuSet \gets \matchHeuSet \cup \{\Elem_{n}\}$
            \ELSIF{$\Elem_{n} \in \vertexSet$ and $\Elsta_{n} =$ \textbf{false} and $\sum_{\parc \in \arcSet} \HatoneCe_{\narc}^{\text{a}} < \aBudget$}
            \STATE $\Elwei_{n} = \text{\textbf{Weight}}(\Elem_{n},\genmatch)$
            ; $\matchHeuSet \gets \matchHeuSet \cup \{\Elem_{n}\}$
            \ENDIF
        \ENDFOR
        \IF{$\matchHeuSet = \emptyset$}
        \STATE Return $\HatoneCe$ and \cover\label{firstEndHeu}
        \ENDIF
        \STATE Sort $\matchHeuSet$ in non-increasing order of values s.t. $\Elwei_{1} \ge \Elwei_{2} ...\ge \Elwei_{\lvert \ElSet \rvert}$ 
        \STATE Select $\Elemstar_{n} \in \matchHeuSet$ randomly among elements whose weight equals $\Elwei_{1}$
        \STATE \textit{ans} = \textbf{true}
        \STATEnonum \textbf{---Start single-vertex-arc separation \ref{vxtarcsep}---}
        \IF{$(\sum_{\vertex \in \vertexSet} \HatoneCe_{\vertex}^{\text{v}} + \sum_{\parc \in \arcSet} \HatoneCe_{\narc}^{\text{a}}) \ge 1$ and $\text{\textbf{IsNDD}}(\Elemstar_{n}) = \textbf{false}$ and $(\sum_{\parc \in \arcSet} \HatoneCe_{\narc}^{\text{a}} \ge 1 \text{ or } \sum_{\vertex \in \vertexSet} \HatoneCe_{\vertex}^{\text{v}} \ge 2)$} \label{SepHeu}
        \STATE $V(\oneCeNoPrime^{\text{v}}) \gets V(\HatoneCe^{\text{v}})$; $A(\oneCeNoPrime^{\text{a}}) \gets A(\HatoneCe^{\text{a}})$
        \STATE \textit{ans} = \textbf{false}
            \IF{$\Elemstar_{n} \in \vertexSet$}
                \STATE $\bar{\vertexb} \gets \Elemstar_{n}$; $\cyclechainSetXvtx \gets$ Check Case \ref{CaseOne}
                \IF{$\cyclechainSetXvtx \neq \emptyset$}
                    \STATE \textit{ans} = \textbf{true}
                \ENDIF
            \ELSE 
                \STATE $\bar{a} \gets \Elemstar_{n}$; $\cyclechainSetXarc \gets $ Check Case \ref{CaseTwo}
                \IF{$\cyclechainSetXarc \neq \emptyset$}
                    \STATE \textit{ans} = \textbf{true}
                \ENDIF
            \ENDIF
        \ENDIF
        \STATEnonum \textbf{---End separation---}
        \STATE $\Elemstar_{n} = \text{\textbf{true}}$ \label{checked}
        \IF{\textit{ans} = \textbf{true}}
            \IF{$\Elemstar_{n} \in \vertexSet $}
                \STATE  $\HatoneCe^{\text{v}}_{\Elemstar_{n}} = 1$; $\Elsta_{n} =$ \textbf{true}  $\forall \ \Elem_{n} \in \deltaOutEl \cup \deltaInEl$ \label{adjArcstHeu}
            \ELSE
                \STATE  $\HatoneCe^{\text{a}}_{\Elemstar_{n}} = 1$
            \ENDIF
        \ENDIF
        \IF{$\HatoneCe \text{ covers all matchings in }  \genmatch$}
            \STATE \cover = \textbf{true}
        \ENDIF
        \IF{\cover = \textbf{true}}
            \IF{$\sum_{\vertex \in \vertexSet} \HatoneCe_{\vertex}^{\text{v}} = \vBudget$ and $\sum_{\parc \in \arcSet} \HatoneCe_{\narc}^{\text{a}} = \aBudget$} \label{budgetHeu}
             \STATE Return $\HatoneCe$ and \cover \label{SndEndHeu}
            \ENDIF
        \ENDIF
        \ENDWHILE
	\caption{Heuristic}
	\label{Alg:Heuristics}
	\end{algorithmic}
\end{algorithm}

\section{Computational Experiments}
\label{Experiments}

    We test our framework on the same instances of 20, 50, and 100 vertices tested in \citet{Blom2021, Carvalho2021}. There are 30 instances for each network size. For homogeneous failure, we focus on the 100-vertex sets, which are more computationally challenging. For non-homogeneous failure, we focus on the 50- and 100-vertex sets.  In the first part of this section, we compare the efficiency of our solution approaches to the state-of-the-art algorithm addressing the full-recourse policy under homogeneous failure, Benders-PICEF, proposed by \citet{Blom2021}. Next, we analyze the performance and practical impacts of our approaches under non-homogeneous failure for the full-recourse and first-stage-only recourse policies presented in Section \ref{sec:RecoursePolicies}. Our implementations are coded in C++ using CPLEX 12.10, including the state-of-the-art algorithm we compare against, on a machine with Debian GNU/Linux as the operating system and a 3.60GHz processor Intel(R) Core(TM). A time limit of one hour is given to every run. The TR value for the HSA\_MB and HSA\_ME algorithms was 150 iterations for each run.

\subsection{Homogeneous failure analysis}

\begin{figure}[tbp]
  \centering
  \includegraphics[width=0.8\linewidth]{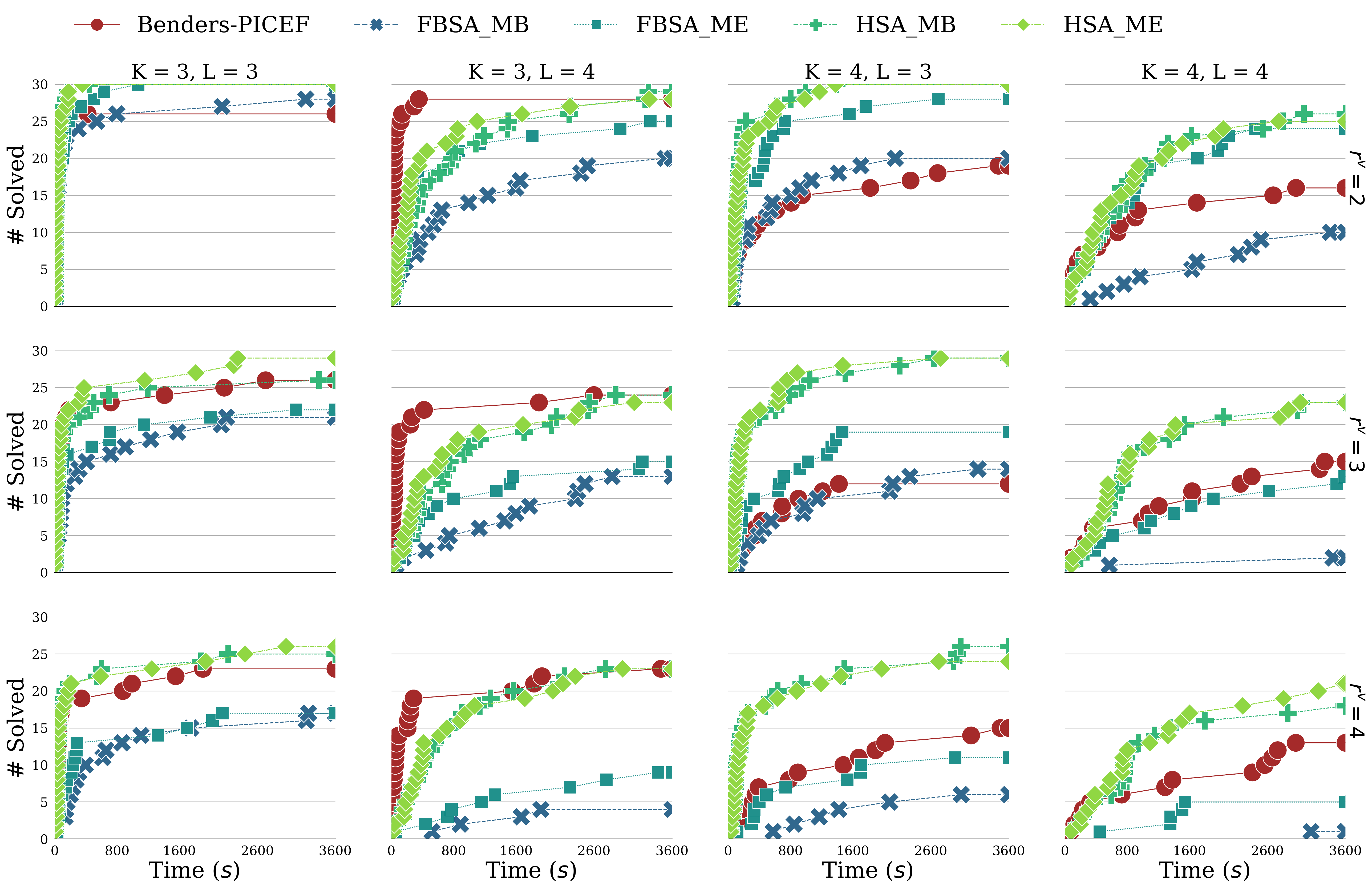}
  \caption{Performance profile for multiple $\vBudget$ values and $\aBudget = 0$ under full recourse}
  \label{fig:Performance}
\end{figure}

Figure \ref{fig:Performance} shows the performance profile for five algorithms: our lifted-constraint version of Benders-PICEF, FBSA\_MB, FBSA\_ME, HSA\_MB and HSA\_ME. Recall that FBSA\_MB and FBSA\_ME are feasibility-based algorithms without the optimization step. Figure \ref{fig:Performance} shows that solving the recourse problem in the transitory graph pays off for FBSA\_ME compared to FBSA\_MB. Although Benders-PICEF solves more instances in general than FBSA\_MB and FBSA\_ME, the performance of FBSA\_ME noticeably outperforms that of Benders-PICEF when the maximum length of cycles is four and that of chains is three and up to three vertices are allowed to fail. In the same settings, FBSA\_MB is comparable to Benders-PICEF. However, as soon as the feasibility-seeking master problems incorporate the optimization step (Algorithm \ref{Alg:HybridHSA-ME}), the performance of HSA\_MB and HSA\_ME is consistently ahead of all other algorithms. As we show shortly, in most cases HSA\_MB and HSA\_ME need a small percentage of iterations in the optimization version of the second-stage master problem to converge. Benders-PICEF is fast when cycles and chains of size up to three and four are considered, respectively; however, increased cycle length and budget failure result in worse performance of Benders-PICEF compared to the hybrid approaches.

Table \ref{tab:performance} summarizes the computational performance details of FBSA\_ME, HSA\_ME and Benders-PICEF. On average, column 2SndS divided by column 1stS indicates the number of iterations that were needed per first-stage iteration to solve the second-stage problem. The average number of iterations per first-stage iteration can also be obtained by adding up columns Alg-\ref{Alg:Heuristics}-true, \ref{MasterExp} and \ref{objMPOptQ} and then dividing the sum by column 1stS. The first observation is that the CG algorithm successfully finds optimal recourse solutions in most cases and it is responsible for most of the average total time for FBSA\_ME and HSA\_ME. 
On the other hand, the average total number of iterations needed by FBSA\_ME to converge was significantly higher than that needed by HSA\_ME and yet the average total time per iteration of the second stage with respect to the total time (2ndS/total) is between 0.16 and 1.84 seconds. 

FBSA\_ME clearly generates scenarios that \ref{objMPOptQ} would not explore due to its cycle-and-chain decision variables and minimization objective, both properties that FBSA\_MB and FBSA\_ME lack. However, also for the full recourse policy, \citet{Blom2021} tested a master problem analogous to \ref{objMPOptQ} and showed that due to the large number of cycles and chains in the formulation, its scalability is limited as the size of cycles and chains grows. We find that solving the feasibility problem is much more efficient than solving \ref{objMPOptQ}. In fact, (i) the heuristic Alg\ref{Alg:Heuristics}, even when close to 1700 iterations, only accounts for about 26\% of the CPU time, and (ii) even when the heuristic  fails and \ref{MasterExp} is solved as a MIP, around 1000 iterations are needed to account for about 39\% of the total time, while for \ref{objMPOptQ}, $<$200 iterations are needed to account for about the same percentage.

Observe that in some cases the average total number of second-stage iterations spent by HSA\_ME is a third of those spent by FBSA\_ME, indicating that \ref{objMPOptQ} helps convergence. However, the number of \ref{objMPOptQ} iterations per first-stage solution (\ref{objMPOptQ}/1stS) is on average 40.5 while the same statistics for the feasibility-seeking iterations are on average 132.5, thus, highlighting that \ref{MasterExp} helps reduce the number of iterations spent by \ref{objMPOptQ}. As a result, the hybrid algorithms are able to converge quickly in all the tested settings. Thus, the feasibility-seeking master problems are able to find near-optimal failure scenarios within the first hundreds of iterations, but in the absence of a lower bound, they may require thousands of iterations to reach infeasibility and thus prove optimality.

\begin{table}[tbp]
    \caption{Computational performance details of FBSA\_ME, HSA\_ME and Benders-PICEF under homogeneous failure and full recourse for 30 instances with $\lvert\vertexSet \rvert = 100$. MPBP and RPBP are the master problem and recourse problem presented in \citet{Blom2021}, respectively. A column is left empty if it does not apply. Data includes both optimal and sub-optimal runs.}
    \label{tab:performance}
    \begin{adjustbox}{width=1\textwidth}
        \centering
        \begin{small}
    \begin{tabular}{rrrrrrrrrrrrrrrrr}
        \toprule
        &  &  &  \multicolumn{14}{c}{FBSA\_ME}\\
        \cmidrule(lr){4-17}
         &  &  &  &   \multicolumn{6}{c}{Avg.\ time} & \multicolumn{7}{c}{Avg. total \# iterations}\\
        \cmidrule(lr){5-10} \cmidrule(lr){11-17}
        $\cycleCap$   & $\chainCap$   & $\vBudget$ &   opt &   total (s) &   Alg\ref{Alg:Heuristics} (\%) &   \ref{MasterExp} (\%) &   \ref{RecoP:Objexp} (\%) &   CG (\%) &\ref{objMPOptQ} (\%)  & 1stS &  2ndS     & Alg\ref{Alg:Heuristics}-true &   \ref{MasterExp} &   \ref{RecoP:Objexp} &   CG-true & \ref{objMPOptQ}\\
    \hline
       3 &   3 &              2 &                           30 &                          130.86 &                             3.39 &            10.77&      2.56             &62.18 & - &                               1.50 &                               241.23 &                                   67.90  &                            173.33 &                               4.23 &                                 237.00    &  -                         \\
       3 &   3 &              3 &                           22 &                         1252.26 &                            14.93 &            31.89&      1.90             &37.68  & - &                             1.33  &                              1190.77 &                                  251.67 &                            939.10  &                              26.23 &                                1164.53  &    -                      \\
       3 &   3 &              4 &                           17 &                         1804.25 &                            26.32 &            39.55&      1.52             &24.91 & - &                              1.31 &                              1698.69 &                                  434.34 &                           1264.34 &                              34.07 &                                1664.62  &   -                       \\
       3 &   4 &              2 &                           25 &                         1076.27 &                             0.44 &            2.59&       4.80             &70.77  & - &                               1.23 &                                229.83 &                                   70.60  &                            159.23 &                               3.43 &                                 226.40 &   -                        \\
       3 &   4 &              3 &                           15 &                         2226.59 &                             3.71 &            7.72&      4.79             &65.04  & - &                              1.21 &                                768.66 &                                  186.38 &                            582.28 &                              12.48 &                                 756.17 &    -                       \\
       3 &   4 &              4 &                            9 &                         2952.83 &                             5.40  &            9.80&      3.51              &63.95 & - &                             1.13 &                               977.13 &                                  249.80  &                            727.33 &                              11.50  &                                 965.63   &   -                      \\
       4 &   3 &              2 &                           28 &                          638.01 &                             1.38 &            5.11&      1.24              &74.94 & - &                              1.47 &                               336.40 &                                   85.00    &                            251.40  &                               3.87 &                                 332.53  &   -                       \\
       4 &   3 &              3 &                           19 &                         1681.14 &                             8.01 &            20.28&     0.84             &62.28  & - &                              1.40 &                              1097.70 &                                  154.70  &                            943.00    &                              19.60  &                                1078.10   &    -                      \\
       4 &   3 &              4 &                           11 &                         2632.47 &                            17.62 &            23.54&     0.70             & 53.31  & - &                             1.23 &                              1722.87 &                                  445.77 &                           1277.10  &                              26.27 &                                1696.60    &    -                     \\
       4 &   4 &              2 &                           24 &                         1385.25 &                             0.29 &            1.56&      2.49              &81.29 & - &                               1.27 &                               219.27 &                                   67.17 &                            152.10  &                               3.03 &                                 216.23 &      -                      \\
       4 &   4 &              3 &                           13 &                         2658.12 &                             1.59 &            4.16&      2.15              &79.13 & - &                               1.23 &                               581.60 &                                  149.27 &                            432.33 &                               7.33 &                                 574.27 &       -                    \\
       4 &   4 &              4 &                            5 &                         3197.12 &                             4.37 &            5.91&      2.11              &75.08 & - &                             1.10 &                               856.62 &                                  284.28 &                            572.34 &                              10.72 &                                 845.90    &     -                     \\
       \cmidrule(lr){1-17}
        &  &  &  \multicolumn{14}{c}{HSA\_ME}\\
       \cmidrule(lr){4-17}
        &  &  &  &  \multicolumn{6}{c}{Avg.\ time} & \multicolumn{7}{c}{Avg. total \# iterations}\\
        \cmidrule(lr){5-10} \cmidrule(lr){11-17}
       $\cycleCap$   & $\chainCap$   & $\vBudget$ &   opt &   total (s) &   Alg\ref{Alg:Heuristics} (\%) &   \ref{MasterExp} (\%) &   \ref{RecoP:Objexp} (\%)  &   CG(\%) &\ref{objMPOptQ} (\%)  &1stS &  2ndS &Alg\ref{Alg:Heuristics}-true &   \ref{MasterExp} &   \ref{RecoP:Objexp} &  CG-true & \ref{objMPOptQ}\\
    \hline
       3 &   3 &   2 &               30 &                          48.91 &                            1.71 &                               7.62 &                                2.94 &                              57.36 & 4.84 &              1.47&  156.03&               62.57 &                            84.87 &                              1.87 &                                154.17  & 8.60  \\
       3 &   3 &   3 &               29 &                         436.87 &                            2.39 &                              10.61 &                                3.01 &                              44.25 & 23.18 &              1.67&  322.83&               123.30  &                           106.27 &                              3.70  &                                319.13 & 93.27 \\
       3 &   3 &   4 &               26 &                         825.24 &                            2.60  &                               9.07 &                                2.32 &                              35.9  & 36.81 &             1.70&  379.80&              175.80  &                            78.03  &                              5.53 &                                374.27  &125.97\\
       3 &   4 &   2 &               28 &                         719.42 &                            0.27 &                               1.85 &                                5.44 &                              69.56 & 0.75 &              1.20&  130.03&                56.37 &                           66.47  &                              1.20  &                                128.83 &7.20\\
       3 &   4 &   3 &               23 &                        1409.03 &                            0.40  &                               1.73 &                                3.57 &                              63.41 & 11.11 &              1.67&  318.10&               132.17 &                           98.03 &                              2.97 &                                315.13  &87.90\\
       3 &   4 &   4 &               23 &                        1469.34 &                            0.47 &                               1.63 &                                4.51 &                              59.50  & 14.07 &              1.43&  295.10&               155.50  &                           56.50  &                              3.77 &                                291.33 &83.10\\
       4 &   3 &   2 &               30 &                         264.44 &                            0.69 &                               2.94 &                                1.27 &                              73.33 & 2.05 &              1.57&  202.13&                81.40  &                           107.50 &                              1.97 &                                200.17 &13.23\\
       4 &   3 &   3 &               29 &                         473.65 &                            1.81 &                               2.97 &                                1.08 &                              68.36 & 11.97 &              1.50&  280.87&               136.00    &                           79.60 &                              3.07 &                                277.80   &65.27\\
       4 &   3 &   4 &               24 &                        1099.37 &                            2.29 &                               1.66 &                                0.65 &                              55.46 & 28.64 &              1.67&  426.70&               197.63 &                           52.17 &                              5.43 &                                421.27  &176.90\\
       4 &   4 &   2 &               25 &                        1243.41 &                            0.21 &                               1.15 &                                2.14 &                              80.57 & 0.79 &              1.50&  194.57&                81.77 &                           97.87  &                              2.40  &                                192.17 &14.93\\
       4 &   4 &   3 &               23 &                        1554.62 &                            0.46 &                               0.89 &                                1.79 &                              78.23 & 4.04 &              1.50&  245.47&               127.73 &                           81.17 &                              4.10  &                                241.37  &36.57\\
       4 &   4 &   4 &               21 &                        1900.28 &                            0.58 &                               0.54 &                                1.84 &                              72.29 & 11.25 &              1.60&  291.63&               190.50  &                           35.73 &                              4.40  &                                287.23  &65.40\\
    \cmidrule(lr){1-17}
        &  &  &  \multicolumn{14}{c}{Benders-PICEF}\\
       \cmidrule(lr){4-17}
        &  &  &  &   \multicolumn{6}{c}{Avg.\ time} & \multicolumn{7}{c}{Avg. total \# iterations}\\
        \cmidrule(lr){5-10} \cmidrule(lr){11-17}
       $\cycleCap$   & $\chainCap$   & $\vBudget$ &   opt &   total (s) &   Alg\ref{Alg:Heuristics} (\%) &   MPBP (\%) &   RPBP (\%)  &   CG(\%) &\ref{objMPOptQ} (\%)  &1stS &  2ndS &Alg\ref{Alg:Heuristics}-true &   MPBP &   RPBP &  CG-true & \ref{objMPOptQ}\\
       \hline
       3 &     3 &                2 &      26 &    271.82 & - &   6.64 &  78.39 & - & - &        8.14 &      167.29 & - &         167.29 &          167.29&  - & - \\
       3 &     3 &                3 &      26 &    722.85 & - &  25.65 &  66.11 & - & - &       11.83 &      278.70  & - &         278.70  &          278.70& - & - \\
       3 &     3 &                4 &      23 &    854.58 & - &  36.04 &  60.33 & - & - &       8.00    &      516.86 & - &         516.86 &          516.86&  - & - \\
       3 &     4 &                2 &      28 &    135.46 & - &   6.70  &  80.54 & - & - &        3.69 &      126.93 & - &         126.93 &          126.93&  - & - \\
       3 &     4 &                3 &      24 &    867.87 & - &  18.64 &  75.15 & - & - &        5.57 &      235.80  & - &         235.80  &          235.80& -  & -  \\
       3 &     4 &                4 &      23 &   1122.50  & - &  24.25 &  71.05 & - & - &        8.23 &      386.63 & - &         386.63 &          386.63& -  & - \\
       4 &     3 &                2 &      19 &   1801.11 & - &   0.73 &  92.42 & - & - &        1.21 &       22.45 & - &          22.45 &           22.45&  -  & - \\
       4 &     3 &                3 &      12 &   2434.24 & - &   2.26 &  91.52 & - & - &        5.07 &      134.13 & - &         134.13 &          134.13&  - & - \\
       4 &     3 &                4 &      15 &   2418.72 & - &   3.73 &  91.97 & - & - &        1.17 &       44.13 & - &          44.13 &           44.13& -  & - \\
       4 &     4 &                2 &      16 &   2099.70  & - &   0.63 &  91.13 & - & - &        3.27 &       34.87 & - &          34.87 &           34.87&  - & - \\
       4 &     4 &                3 &      15 &   2486.05 & - &   2.31 &  93.35 & - & - &        2.60  &      126.00    & - &         126.00    &          126.00&  - & -    \\
       4 &     4 &                4 &      13 &   2659.81 & - &   2.81 &  93.09 & - & - &        3.40  &      136.90  & - &         136.90  &          136.90&   - & - \\
    \hline
    \end{tabular}     
\end{small}
    \end{adjustbox}

{ \footnotesize
\begin{description}
\item $(\%)$:  Percentage of total time spent in that part of the optimization
\item \ref{MasterExp}:  Time solving MT as a MIP when the heuristic failed to find a feasible solution
\item \ref{RecoP:Objexp}:  Time solving RE as a MIP when the CG algorithm failed to find an optimal recourse solution
\item \ref{objMPOptQ}: Total \# of iterations SSF was solved as a MIP when more than TR = 150 iterations passed without MT becoming infeasible
\item 1stS:  Total \# of first-stage decisions required before finding the robust solution
\item 2SndS:  Total \# of iterations to solve the second-stage problem for all the first-stage decisions
\item Alg\ref{Alg:Heuristics}-true: Total \# of iterations the heuristic found a feasible failure scenario
\item CG-true:  Total \# of iterations the CG algorithm found an optimal recourse solution
\end{description}
}
\end{table}

\subsection{Non-homogeneous failure analysis}
In this part of the analysis, we focus on understanding both the scalability of the solution algorithms under non-homogeneous failure and the impact of the two recourse policies on the total number of pairs that can be re-arranged into new cycles and chains. We additionally examine the percent of re-arranged pairs that correspond to highly-sensitized patients, who are more impacted by non-homogenous failure rates due to their real-life higher failure rates. Pairs in the instances published by \cite{Carvalho2021} have an associated panel reactive antibody (PRA), which determines how likely a patient may reject a potential donor. The higher this value, the more sensitized a patient becomes and thus the patient becomes less likely to get a compatible donor. Typically, a patient with a PRA greater or equal to 90\% is considered sensitized. 

Table \ref{tab:policiesk3} presents computational details for HSA\_MB with a maximum cycle length of three. A similar table for a maximum cycle length of four is presented in Appendix \ref{CycleLengthFourNonHomoFail}. The average total number of dominated scenarios for instances in Table \ref{tab:performance} was negligible and thus, not presented. The non-homogeneous case is more difficult to solve based on the average total time, the average total number of first-stage iterations, and the average percentage of total time to solve the optimality problem \ref{objMPOptQ} after 150 iterations. The average number of dominated scenarios is particularly high for the first-stage-only recourse. This behaviour may be explained because when a failure occurs, the second-stage compatibility graph under this policy is smaller and thus more susceptible to having vertices and arcs that can no longer be part of cycles and chains. In terms of the robust objective, all instances that were able to be solved optimally under both policies, obtained the same objective. This interesting fact indicates that for the tested instances, there exists a set of transplants pairs that can be recovered amongst themselves, and there are that many transplants that can be recovered by the full-recourse policy. 

\begin{table}[tbp]
    \caption{Computational performance details of HSA\_MB under non-homogeneous failure for 30 instances with $\cycleCap = 3$. Columns are mostly the same as in Table \ref{tab:performance}, with new columns defined below. Data includes both
optimal and sub-optimal runs.$^*$}
    \label{tab:policiesk3}
    \begin{adjustbox}{width=\textwidth}
        \centering
        \begin{tabular}{rrrrrrrrrrrrrrrrrrr}
    \toprule
    &  &  &  \multicolumn{16}{c}{Full recourse}\\
    \cmidrule(lr){4-19}
     &  &  &  &  & \multicolumn{6}{c}{Time} & \multicolumn{7}{c}{Avg. total \#}\\
    \cmidrule(lr){6-11} \cmidrule(lr){12-19}
    $\chainCap$   & $\vBudget$   &  $\aBudget$  &      opt &      HSP (\%) &   total (s)     & Alg\ref{Alg:Heuristics} (\%) &   \ref{MasterExp}(\%) &   \ref{RecoP:Objexp} (\%) &   CG (\%)         &\ref{objMPOptQ} (\%)  & 1stS      &  2ndS         &   Alg\ref{Alg:Heuristics}-true  &   \ref{MasterExp} &   \ref{RecoP:Objexp} &      CGtrue &   \ref{objMPOptQ} &   DomS           \\
       \cmidrule(lr){1-19}
          $\lvert\vertexSet\rvert = 50$&             &             &         &           &          &              &                 &                  &                 &               &              &              &              &                &                 &                &              &          \\
               3 &           1 &        1.80  &      30 &   34.21 &      20.13 &         2.69 &           35.18 &             2.06 &           22.30  &         18.08 &         3.73 &       395.13 &       107.73 &         264.83 &            8.87 &         386.27 &        22.57 &         2.13 \\
               3 &           2 &        1.80  &      29 &   25.36 &     220.27 &         2.94 &           26.95 &             0.81 &           17.18 &         42.41 &         3.30  &       481.67 &       178.90  &         249.83 &           10.47 &         471.20  &        52.93 &         1.70  \\
               3 &           1 &        2.97 &      30 &   28.54 &     160.76 &         2.59 &           30.82 &             1.05 &           17.35 &         35.13 &         6.30  &       784.70  &       185.93 &         546.10  &           13.33 &         771.37 &        52.67 &         3.97 \\
               3 &           2 &        2.97 &      25 &   25.12 &     699.03 &         1.99 &           23.75 &             0.57 &           11.47 &         56.42 &         3.03 &       545.23 &       198.63 &         229.10  &           13.50  &         531.73 &       117.50  &         4.77 \\
               4 &           1 &        1.80  &      30 &   34.33 &      27.07 &         1.71 &           22.58 &             5.02 &           34.48 &          8.87 &         2.87 &       338.70  &       101.97 &         223.90  &            7.73 &         330.97 &        12.83 &         1.93 \\
               4 &           2 &        1.80  &      28 &   25.49 &     316.27 &         1.76 &           17.42 &             2.73 &           27.20  &         32.20  &         2.87 &       439.03 &       173.80  &         222.93 &            5.33 &         433.70  &        42.30  &         1.33 \\
               4 &           1 &        2.97    &      29 &   28.20 &     327.59 &         1.43 &           16.87 &             3.22 &           24.20  &         34.20  &         6.23 &       803.80  &       219.20  &         515.90  &           12.30  &         791.50  &        68.70  &         3.47 \\
               4 &           2 &        2.97    &      26 &   22.82  &     663.25 &         1.53 &           15.36 &             2.71 &           20.22 &         46.67 &         2.90  &       478.27 &       203.40  &         203.00    &            8.80  &         469.47 &        71.87 &         3.40  \\
       \cmidrule(lr){1-19}
       $\lvert\vertexSet \rvert =  100$&                &           &          &           &          &              &                 &                  &                 &               &              &              &              &                &                 &                &              &          \\
              3 &           1 &        3.27 &      16 &   48.88 &    2182.03 &         0.40  &            3.28 &             1.28 &            8.72 &         81.56 &         3.13 &       633.47 &       193.83 &         265.00    &           22.53 &         610.93 &       174.63 &         3.83 \\
              3 &           2 &        3.27 &       4 &   33.11 &    3205.20  &         0.33 &            2.35 &             0.39 &            4.17 &         89.23 &         1.93 &       376.23 &       226.60  &          62.53 &           10.47 &         365.77 &        87.10  &         1.30  \\
              3 &           1 &        5.90  &       2 &   35.43 &    3407.89 &         0.05 &            0.34 &             0.24 &            1.17 &         95.16 &         1.77 &       396.83 &       129.93 &         133.80  &           13.23 &         383.60  &       133.10  &        11.17 \\
              3 &           2 &        5.90  &       3 &   26.11  &    3295.06 &         0.04 &            0.14 &             0.15 &            0.95 &         93.96 &         1.33 &       250.63 &       169.33 &          29.77 &            9.97 &         240.67 &        51.53 &         6.23 \\
              4 &           1 &        3.27 &       7 &   41.48  &    2925.14 &         0.14 &            1.20  &             7.80  &           34.43 &         44.43 &         2.37 &       475.50  &       214.67 &         133.80  &           11.80  &         463.70  &       127.03 &         1.60  \\
              4 &           2 &        3.27 &       4 &   37.62 &    3144.74 &         0.09 &            0.41 &             3.81 &           23.47 &         63.66 &         1.47 &       287.70  &       184.43 &          32.27 &            8.53 &         279.17 &        71.00    &         0.50  \\
              4 &           1 &        5.90  &       2 &   20.25 &    3457.09 &         0.05 &            0.22 &             4.08 &           18.12 &         70.91 &         1.57 &       320.27 &       163.17 &          68.47 &            8.73 &         311.53 &        88.63 &         2.20  \\
              4 &           2 &        5.90  &       1 &   27.27 &    3597.92 &         0.02 &            0.08 &             2.70  &           14.77 &         74.21 &         1.30  &       240.43 &       172.37 &          19.27 &            8.03 &         232.40  &        48.80  &         1.00    \\
       \cmidrule(lr){1-19}
    &  &  &  \multicolumn{16}{c}{First-stage-only recourse}\\
    \cmidrule(lr){4-19}
     &  &  &  &  & \multicolumn{6}{c}{Time} & \multicolumn{7}{c}{Avg. total \#}\\
    \cmidrule(lr){6-11} \cmidrule(lr){12-19}
    $\chainCap$   & $\vBudget$   &  $\aBudget$  &      opt &     HSP (\%)  & total (s)       & Alg\ref{Alg:Heuristics} (\%) &   \ref{MasterExp}(\%) &   \ref{RecoP:Objexp} (\%) &   CG (\%)         &\ref{objMPOptQ} (\%)  & 1stS      &  2ndS         &   Alg\ref{Alg:Heuristics}-true  &   \ref{MasterExp} &   \ref{RecoP:Objexp} &      CGtrue &   \ref{objMPOptQ} &   DomS           \\
       \cmidrule(lr){1-19}
       $\lvert\vertexSet\rvert = 50$&             &             &         &           &          &              &                 &                  &                 &               &              &              &              &                &                 &                &              &          \\
       3 &           1 &        1.80  &      21 &   36.12 &    1087.77 &         2.72 &           27.77 &             0.89 &           19.14 &          9.34 &        53.70  &      4193.23 &      1562.20  &        2535.73 &           82.80  &        4110.43 &        95.30  &       170.80  \\
       3 &           2 &        1.80  &      18 &   29.96 &    1478.79 &         3.43 &           19.63 &             0.37 &           13.05 &         36.94 &        90.53 &      9182.00   &      3180.70  &        5645.50  &           99.10  &        9082.90  &       355.80  &      1338.67 \\
       3 &           1 &        2.97 &      16 &   28.23 &    1752.26 &         3.09 &           20.81 &             0.60  &           11.51 &         28.90  &        85.10  &      8549.13 &      2214.97 &        6030.57 &           80.33 &        8468.80  &       303.60  &      1573.20  \\
       3 &           2 &        2.97 &      18 &   27.39 &    1680.72 &         2.94 &           17.91 &             0.38 &           11.37 &         50.09 &        60.23 &      7778.93 &      2425.60  &        4997.17 &          134.30  &        7644.63 &       356.17 &       920.67 \\
       4 &           1 &        1.80  &      27 &   32.84 &     390.91 &         2.64 &           21.45 &             2.09 &           40.54 &          5.68 &        23.37 &      1718.87 &       559.00    &        1143.70  &           20.67 &        1698.20  &        16.17 &        43.43 \\
       4 &           2 &        1.80  &      19 &   28.36 &    1344.28 &         2.66 &           16.44 &             1.23 &           27.06 &         27.99 &        61.30  &      6597.37 &      2517.10  &        3770.47 &           60.53 &        6536.83 &       309.80  &       235.07 \\
       4 &           1 &        2.97    &      16 &   28.44 &    1736.95 &         2.56 &           14.78 &             1.17 &           27.97 &         21.99 &        63.13 &      7111.33 &      2058.10  &        4629.93 &          173.60  &        6937.73 &       423.30  &       785.73 \\
       4 &           2 &        2.97    &      16 &   22.98 &    1725.39 &         2.12 &           12.97 &             0.73 &           21.23 &         41.89 &        75.13 &      8909.87 &      2798.9  &        5779.70  &          159.60  &        8750.27 &       331.27 &      2502.30  \\
       \cmidrule(lr){1-19}
       $\lvert\vertexSet \rvert =  100$&             &             &         &           &          &              &                 &                  &                 &               &              &              &              &                &                 &                &              &          \\
       3 &           1 &        3.27 &      13 &   45.10 &    2405.01 &         0.52 &            5.62 &             1.74 &            9.33 &         74.43 &        11.40  &      1884.80  &       492.00    &        1136.20  &          128.00    &        1756.80  &       256.60  &        36.07 \\
       3 &           2 &        3.27 &       3 &   51.89 &    3326.33 &         0.19 &            0.77 &             0.61 &            3.59 &         90.97 &         7.53 &      1322.60  &       648.83 &         480.00    &           62.67 &        1259.93 &       193.77 &        74.47 \\
       3 &           1 &        5.90  &       1 &   12.50 &    3504.05 &         0.11 &            0.53 &             0.16 &            1.35 &         95.12 &         5.80  &      1048.97 &       375.50  &         487.57 &           35.43 &        1013.53 &       185.90  &       196.97 \\
       3 &           2 &        5.90  &       1 &   32.00 &    3519.36 &         0.13 &            0.46 &             0.21 &            1.12 &         90.28 &         6.10  &      1002.27 &       497.10  &         416.97 &           97.80  &         904.47 &        88.20  &        98.97 \\
       4 &           1 &        3.27 &       9 &   47.28 &    2842.43 &         0.27 &            2.02 &             3.71 &           43.35 &         40.41 &         5.60  &      1026.67 &       409.57 &         414.57 &           19.20  &        1007.47 &       202.53 &        44.47 \\
       4 &           2 &        3.27 &       5 &   40.91 &    3228.62 &         0.09 &            0.29 &             1.56 &           19.30  &         72.94 &         2.70  &       503.83 &       313.03 &          85.83 &           13.43 &         490.40  &       104.97 &         6.57 \\
       4 &           1 &        5.90  &       2 &   20.25 &    3393.66 &         0.09 &            0.44 &             1.91 &           16.86 &         76.00    &         3.23 &       614.73 &       274.70  &         204.60  &           37.47 &         577.27 &       135.43 &        26.83 \\
       4 &           2 &        5.90  &       0 &   -     &    3600.55 &         0.05 &            0.13 &             1.12 &           13.44 &         79.01 &         3.30  &       575.37 &       371.07 &         121.10  &           27.30  &         548.07 &        83.20  &        36.33 \\
    \bottomrule
 \end{tabular}
    \end{adjustbox}
{ \footnotesize
\begin{description}
\item $\aBudget$: Average number of arcs that fail, corresponding to either 5\% or 10\% of the arcs in the deterministic solution of an instance
\item $^*$HSP $(\%)$: Average percentage of highly-sensitized patients re-matched in the second stage. Includes only instances solved to optimality.
\item DomS: Average total number of dominated scenarios that are found by the single-vertex-arc separation procedure (Section \ref{vxtarcsep})
\end{description}
}
\end{table}

\begin{figure}[tbp]
  \centering
  \includegraphics[width=0.7\linewidth]{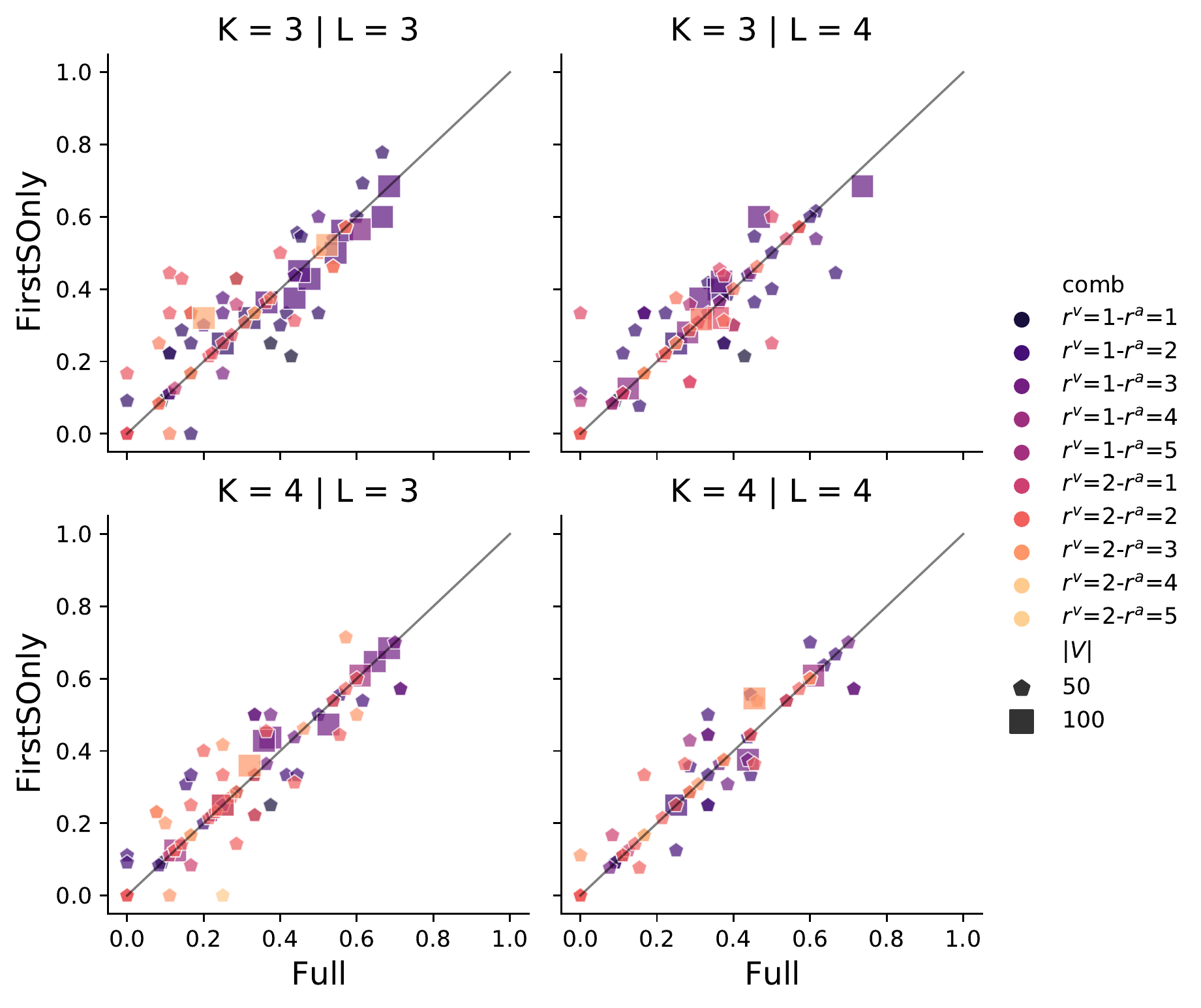}
  \caption{Percentage of highly-sensitized patients selected in both first and second stage with respect to their total number in an instance. Markers above the line indicate a higher percentage for the first-stage-only recourse policy. Every combination (``comb'' in the legend) of $\vBudget$ (number of vertex failures) and $\aBudget$ (number of arc failures) is indicated. The values for $\aBudget$ result from considering 5\% and 10\% of the arcs failing in the deterministic solution of an instance. Only instances that were optimal under both policies are shown.}
  \label{fig:hsp-policies}
\end{figure}

The percentage of highly-sensitized patients that are selected in the first stage and re-matched in the second stage under both full-recourse and first-stage-only recourse policies is shown in Figure \ref{fig:hsp-policies}. On average, around 30\% of highly-sensitized patients are re-matched under both policies, but since this percentage comes from a worst-case scenario, Figure \ref{fig:hsp-policies} is a lower bound on the percentage of highly-sensitized patients that can be transplanted in practice under both policies. While there is variability in all cases, instances with the smallest failure budgets account for most instances that can re-match more than 50\% of the highly-sensitized patients.

\section{Conclusion}      
\label{sec:Conclusion}

We presented a general solution framework for KPDPs to identify robust matchings that can be repaired after failure through a predefined recourse policy. Our framework consists of a two-stage RO problem whose second-stage master problem seeks feasibility rather than optimality, and yet is able to model multiple policies. The most significant limitation of this approach is the inability to obtain a lower bound on the optimal objective value of the second-stage problem, which we overcame by (i) separating dominated scenarios, and (ii) solving an optimality-seeking master problem when needed.

Under homogeneous failure, our purely feasibility-based solution algorithms outperform the state-of-the-art in some settings, and our hybrid methods consistently outperform the state-of-the-art. Under non-homogeneous failure, we solved all instances to optimality under both recourse policies, finding that while the percentage of highly-sensitized patients that can be recovered varies, 30\% seems to be the lower bound. 

Future directions include simulations to assess matching paradigm effectiveness over time, as opposed to just in a single matching, and expanding uncertainty set possibilities. For example, the uncertainty set could model a policy where an intervention is performed over some arcs/vertices prior to proposing a matching (e.g., performing immunotherapy to patients, running a questionnaire on pairs, etc.), thus reducing the likelihood of an exchange failing if resources are allocated for that intervention. Or, failure could be seen as a change in the quality of a matching rather than the non-existence of it. Moreover, adaptations to our framework could include non-unitary weights for vertices/arcs in our definition of the uncertainty set and similarly for pairs in the robust objective.

\section*{Acknowledgements}      
\label{sec:Acknowledgements}
This research was supported by NSERC Discovery Grant, RGPIN-2021-02609

\setlength\bibsep{0.2pt}
\bibliographystyle{apalike}
{
\footnotesize
\bibliography{kep}
}
\clearpage
\appendix

\section{Notation}
\label{sec:notation}

The following is a list of the notation divided by section. Notation is defined in the section where it was first introduced.

\paragraph{Preliminaries}

\begin{itemize}
    \item[] $\cycleCap$: Maximum allowed length of cycles, user defined. 
    \item[] $\chainCap$: Maximum allowed length of chains in terms of arcs, user defined. 
    \item[] $\pairSet$: Set of patient-donor pairs or simply \emph{set of pairs}.
    \item[] $\nddSet$: Set of non-directed donors.
    \item[] $\vertexSet := \pairSet \cup \nddSet$ represents the set of patient-donor pairs, $\pairSet$, and the set of non-directed donors $\nddSet$.
    \item[] $\arcSet \subseteq \vertexSet \times \pairSet$: Arc set containing arc $\parc$ if and only if the donor in vertex $\vertexa \in \vertexSet$ is compatible with patient in vertex $\vertexb \in \pairSet$.
    \item[] $\graphKEP$: First-stage compatibility graph.
    \item[] $\cycleSet$: Set of feasible cycles (up to size $\cycleCap$) in the first-stage compatibility graph.
    \item[] $\chainSet$: Set of feasible chains (up to size $\chainCap$) in the first-stage compatibility graph.
    \item[] $\containedVxt$ and $\containedArc$ be the set of vertices and arcs in $(\cdot)$, where $(\cdot)$ refers to a cycle, chain, etc.
    \item[] $\matchSet \subseteq \cycleSet \cup \chainSet$: A feasible solution to the KEP, i.e., a collection of vertex-disjoint cycles and chains (referred to as a \emph{matching}) if $V(\match) \cap V(\match^{\prime}) = \emptyset,$ for all $\match, \match^{\prime} \in \matchSet$ with $\match \ne \match^{\prime}$, noting that $\match$ and $\match^\prime$ are entire cycles or chains in the matching.
    \item[] $\matchSpace$: Set of all KEP matchings in the first-stage compatibility graph $\graph$.
    \item[] $\setFirst := \{\charVecFirst: \matchSet \in \matchSpace \}$ is the set of all binary vectors representing the selection of a feasible matching in the first-stage compatibility graph, where $\charVecFirst$ is the characteristic vector of matching $\matchSet$ in terms of the cycles/chains sets $\cycleSet \cup \chainSet$.
    \item[] $\varFirst \in \setFirst$: A first-stage solution representing a matching in the first-stage compatibility graph $\graph$.
    \item[] $\policy \in \policySet$: Recourse policy (or simply \emph{policy}) deciding on the cycles and chains that can be used to repair a first-stage decision $\varFirst \in \setFirst$ after observing vertex/arc failures.
    \item[] $\cycleSetSndTr$: Set of all cycles (up to size $\cycleCap$) in the first-stage compatibility graph satisfying policy $\pi \in \policySet$ and including at least one pair in $\varFirst$.
    \item[] $\chainSetSndTr$: Set of all chains (up to size $\chainCap$) in the first-stage compatibility graph satisfying policy $\pi \in \policySet$ and including at least one pair in $\varFirst$.
    \item[] $\digraphSndTr = (\vertexSetSndBe, \arcSetSndBe)$: Transitory graph with vertex set $\vertexSetSndBe = \cup_{\match \in \cycleSetSndTr \cup \chainSetSndTr}V(\match)$ and arc set $\arcSetSndBe = \cup_{\match \in \cycleSetSndTr \cup \chainSetSndTr}A(\match)$.
    \item[] $V(\varFirst) = \cup_{\match' \in \cycleSet \cup \chainSet:\varFirst_{\match'} = 1} V(\match')$: Vertices (pairs and non-directed donors) selected by a first-stage solution $\varFirst \in \setFirst$.
    \item[] $\oneCe \in \bigsCe$: A failure scenario represented as a binary vector in the uncertainty set $\bigsCe$.
    \item[] $\vertexSetSndBeF$ and  $\arcSetSndBeF$: Set of vertices and set of arcs that fail under scenario $\oneCe \in \bigsCe$ in the first-stage compatibility graph, respectively.
    \item[] $\digraphSndS = (\vertexSetSndBe \setminus \{\vertexSetSndBe \cap \vertexSetSndBeF\}, \arcSetSndBe \setminus \{\arcSetSndBe \cap \arcSetSndBeF\})$: A second-stage compatibility graph induced in the transitory graph $\digraphSndTr$ by a failure scenario $\oneCe \in \bigsCe$.
    \item[] $\cycleSetSndS$: Set of existing cycles (up to size $\cycleCap$) in a second-stage compatibility graph  as they do not fail under scenario $\oneCe \in \bigsCe$ and are allowed by policy $\policy \in \policySet$.
    \item[] $\chainSetSndS$: Set of existing chains (up to size $\chainCap$) in a second-stage compatibility graph as they do not fail under scenario $\oneCe \in \bigsCe$ and are allowed by policy $\policy \in \policySet$.
    \item[] $\polySetxOne:= \{\matchSet \subseteq \cycleSetSndS \cup \chainSetSndS \mid V(\match) \cap V(\match^{\prime}) = \emptyset \text{ for all } \match, \match^{\prime} \in \matchSet \text{; } \match \neq \match^{\prime}\}$: Set of allowed recovering matchings under policy $\policy$ such that every cycle/chain in  $\polySetxOne$ contains at least one pair in $\varFirst$.
    \item[] $\setReco:= \{\charVecSecond:\matchSet \in \polySetxOne\}$: Set of all binary vectors representing the selection of a feasible matching in a second-stage compatibility graph $\digraphSndS$ with non-failed elements (vertices/arcs) under scenario $\oneCe \in \bigsCe$ and policy $\policy \in \policySet$ that contain at least one  pair in $\varFirst$. Here, $\charVecSecond$ is the characteristic vector of matching $\matchSet$ in terms of $\cycleSetSndS \cup \chainSetSndS$.
    \item[] $\varReco \in \setReco$: A second-stage solution (a.k.a \emph{recourse solution}) found in a second-stage compatibility graph. 
    \item[] $f(\varFirst,\oneCe,\varReco)$: General robust objective function (a.k.a \emph{recourse objective function}) that assigns weights to the cycles and chains of a recovered matching associated to a recourse solution $\varReco$ under failure scenario $\oneCe$, based on the number of pairs matched in the first stage solution $\varFirst$. 
    \item[] $\mathbf{w}_{\match}(\varFirst) := |V(\match) \cap  \SelVerX \cap \pairSet|$: weight of a cycle/chain $\match \in \cycleSetSndS \cup \chainSetSndS$ corresponding to the \textit{number of pairs} that having been matched in the first stage by solution $\varFirst \in \setFirst$, can also be matched in the second stage by recourse solution $\varReco \in \setReco$, after failures are observed.
    \item[] $\bigsCe(\varFirst) \in \bigsCe$: Observable uncertainty set in the transitory graph when the first-stage solution is $\varFirst$.
    \item[] $\bm{\smallsCeVx} \in \{0,1\}^{\mid\vertexSet \mid}$: Binary vector representing a vertex failure scenario. If the entry corresponding to vertex $\vertexa \in \vertexSet$ in $\bm{\smallsCeVx}$ takes on value one then $\vertexa$ fails, and on value zero otherwise.
    \item[] $\bm{\smallsCeArc} \in \{0,1\} ^{\mid \arcSet \mid}$: Binary vector representing an arc failure scenario. If the entry corresponding to arc $\parc \in \arcSet$ in $\bm{\smallsCeArc}$ takes on value one then $\parc$ fails, and on value zero otherwise.
    \item[] $\vBudget$: User-defined number of vertices that can occur in the first-stage compatibility graph.
    \item[] $\aBudget$: User-defined number of arcs that can occur in the first-stage compatibility graph.
\end{itemize}

\paragraph{Robust model and first stage}

\begin{itemize}
    \item[] $\tilde{\bigsCe}$: Restricted set of scenarios for the first-stage problem \ref{SingleROmodel}.
    \item[] $\varFirstc \in \setFirst$: Optimal first-stage solution to \ref{SingleROmodel} in Algorithm \ref{Alg:ScenarioGeneration} with restricted set of scenarios $\tilde{\bigsCe}$.
    \item[] $\objSingleROr$: Objective value of solution $\varFirstc \in \setFirst$.
\end{itemize}

\paragraph{New second-stage decompositions}
\begin{itemize}
    \item[] $\bigsCecR \subseteq \bigsCec$: Restricted set of failure scenarios observable in the transitory graph.
    \item[] $\cdteOne \in \bigsCec$: A failure scenario that induces the second-stage compatibility graph $\digraphSndSc$.
    \item[] $\RecoPcdte$: MIP formulation (\ref{RecoP:Obj}) for the recourse problem when solved in the second-stage compatibility graph $\digraphSndSc$.
    \item[] $\cychSetSndScdte = \cycleSetSndSc \cup \chainSetSndSc$.
    \item[] $\varRecoOpSndS$: Optimal recourse solution to the recourse problem $\RecoPcdte$ with objective value $\RecoObjVarcdte$.
    \item[] $\objTwoStageStar$: Optimal objective of the second-stage problem \ref{eq:SSP}.
    \item[] $\objTwoStage$: Decision variable indicating the objective value of the second-stage problem \ref{eq:SSP}.  
    \item[] $\UBCovering$: Upper bound on the objective value of the second-stage problem \ref{eq:SSP}.
    \item[] $\mathbbm{1}_{\match,\oneCe}$: Indicating variable that takes on value one if cycle/chain $\match \in \cychSetSndS$ fails under scenario $\oneCe \in \bigsCec$.
    \item[] $\setRecoTrcdte$: Set of binary vectors representing a recourse solution in the transitory graph.
    \item[] $\RecoPexpc$: MIP formulation (\ref{RecoP:Objexp}) for the recourse problem solved in the transitory graph , whose optimal recourse solutions are also optimal to $\RecoPcdte$ but may allow some failed cycles/chains in the solution.
    \item[] $\VarRecoOpExp \in \setRecoTrcdte$: Optimal recourse solution to $\RecoPexpc$ in the transitory graph that is also optimal to $\RecoPcdte$ under scenario $\cdteOne \in \bigsCecR$.
    \item[] $\VarRecoOpExpcdte \subseteq \VarRecoOpExp$: Subset of the optimal recourse solution $\VarRecoOpExp$ that has no failed cycles/chains and thus corresponds to a feasible solution in $\setRecocdte$.
    \item[] $ I(\oneCePrime)$: Number of failed vertices and arcs under scenario $\oneCePrime$.
    \item[] $I(\oneCeNoPrime)$: Number of failed vertices and arcs under scenario $\oneCeNoPrime$.
    \item[] $\cyclechainSetPrimeX \subseteq \cycleSetSndTrc \cup \chainSetSndTrc$: Set of feasible cycles and chains that fail in the transitory graph $\digraphSndTrc$ under scenario $\oneCePrime$.
    \item[] $\cyclechainSetNoPX \subseteq \cycleSetSndTrc \cup \chainSetSndTrc$: Set of feasible cycles and chains that fail in the transitory graph $\digraphSndTrc$ under scenario $\oneCeNoPrime$.
    \item[] $\cyclechainSetXvtx$: Set of feasible cycles and chains in $\digraphSndSpr$ that include vertex $\bar{\vertexb}$.
    \item[] $\cyclechainSetXarc$: Set of feasible cycles and chains in $\digraphSndSpr$ that include arc $\bar{a}$.
\end{itemize}

\paragraph{Solution algorithms for the second-stage problem}
\begin{itemize}
    \item[] FBSA\_MB: Algorithm \ref{Alg:BasicHybrid}.
    \item[] FBSA\_ME: Algorithm \ref{Alg:EnhancedHybrid}.
    \item[] $\LBCovering$: Lower bound on the objective value of the second-stage problem \ref{eq:SSP}.
    \item[] HSA\_MB and HSA\_ME: Algorithms with additional steps needed for FBSA\_MB and FBSA\_ME to find $\LBCovering$, respectively. Such steps are given by Algorithm \ref{Alg:HybridHSA-ME}. 
    \item[] $\RecoSet$,$\RecoSetexp$, $\genmatch$: Set of matchings
    \item[] $\textbf{ToMng}(\cdot)$: Function that ``extracts'' the matching from a recourse solution $(\cdot)$.
    \item[] \textbf{Heuristic}($\cdot$):  Algorithm \ref{Alg:Heuristics} which takes as input a set of matchings $(\cdot)$ and returns a tuple $\HeuAns$.
    \item[] $\textbf{ColGen}(\RecoPcdte)$: Column generation algorithm used to solve the linear relaxation of $\RecoPcdte$ with optimal value $\ZRecoColUB$.
    \item[] $\ZRecoCol$: Objective value of the feasible KEP solution obtained when the columns obtained by  $\textbf{ColGen}(\RecoPcdte)$ are turned to binary.
    \item[] $\textbf{XColGen}\left(\RecoPexpc \right)$: Column generation algorithm used to solve the linear relaxation of $\RecoPexpc$ with objective value $\ZRecoColUBexp$.
    \item[] $\ZRecoColexp$: Objective value of the feasible KEP solution obtained when the columns obtained by  $\textbf{XColGen}(\RecoPexpc)$ are turned to binary.
    \item[] $\textbf{TrueVal}(\ZRecoColexp)$: A function that returns  $\mathbf{w}_{\match}(\varFirst) := |V(\match) \cap  \SelVerX \cap \pairSet|$ for all cycles and chains $\match$ in the feasible solution with objective value $\ZRecoColexp$.
    \item[] \textbf{UniqueElms}($\genmatch$): Function that returns the set of unique vertices and arcs among all matchings in $\genmatch$, referred to as $\ElSet$.
    \item[] $\Elem_{n} \in \ElSet$: Element, either vertex or arc in the set $\ElSet$.
    \item[] $\text{\textbf{Weight}}(\Elem_{n}, \genmatch)$: Function that returns the weight $\Elwei_{n}$ of an element $\Elem_{n}$ corresponding to the number of times $\Elem_{n}$ is repeated in the set of matchings $\genmatch$.
    \item[] $\text{\textbf{IsNDD}}(\Elemstar_{n})$: Function that returns true if element $\Elemstar_{n}$ is a non-directed donor, and false otherwise.
\end{itemize}

\section{Robust formulation for full recourse}
\label{sec:fullROFormulation}
Our solution approach supports the use of any MIP formulation that has the structure shown in Section \ref{RobustModels}. For the results presented in this paper, we adapt a variant presented by \cite{Carvalho2021} of the position-indexed cycle edge formulation (PICEF) proposed by \cite{Dickerson2016}. Cycles are modeled through cycle variables $\varRecoZ_\match \ \forall \match \in \cycleSet$ for first-stage decisions and through $\cyclevarCe \ \forall \match \in \cycleSet$  for recourse cycles under scenario $\oneCe \in \bigsCe$.  However, chains are modeled through first-stage decision variables $\chainvar$, indexed by arc $\parc \in \arcSet$ and the feasible position $\ell \in \posset$ of that arc within a chain. The set $\posset \subseteq \idposset = \{1,...,\chainCap\}$ corresponds to the set of positions for which that arc is reached from some non-directed donor in a simple path with $\ell \le \chainCap$ arcs. For vertices $\vertexa \in \nddSet$, the set of possible arc positions becomes $\posset  = \{1\}$, since non-directed donors always start a chain. To identify $\posset$ for the other arcs, a shortest-path based search can be performed \citep{Dickerson2016}. Thus, we define $\arcSetL = \{\parc \in \arcSet \mid \ell \in \posset\}$. Likewise, recourse chain decision variables for every scenario $\oneCe \in \bigsCe$ are denoted by $\chainvarCe$. A binary decision variable $\vertxCevar$ is also defined for every pair $\vertexb \in \pairSet$ and scenario $\oneCe \in \bigsCe$ to identify the pairs that are selected in both the first stage and in the second stage under some scenario $\oneCe \in \bigsCe$. Moreover, we denote by $\cyclesetv$ and $\cyclesetuv$ the set of feasible cycles including vertex $\vertexb \in \pairSet$ and the set of feasible cycles including arc $\parc \in \arcSet$, respectively.
\label{FullFirstSFormulation}
\begin{subequations}
\small
\begin{align}
    \max \qquad \objPosIdx \label{objmaxRO}\\
    \objPosIdx - \sum_{\vertexb \in \pairSet} \vertxCevar \le 0 && \oneCe \in \bigsCe \label{selidxF}\\
    \vertxCevar - \sum_{\match \in \cyclesetv} \varRecoZ_\match - \sum_{\ell \in \idposset} \sum_{\parc \in \arcSetL} \chainvar \le 0 && \oneCe \in \bigsCe, \vertexb \in \pairSet \label{FirstSSol}\\
    \vertxCevar - \sum_{\match \in \cyclesetv} \cyclevarCe - \sum_{\ell \in \idposset} \sum_{\parc \in \arcSetL} \chainvarCe \le 0 && \oneCe \in \bigsCe, \vertexb \in \pairSet \label{SecondSSol}\\
    \sum_{\vertexa: \parcrev}\chainvarRCe \le 1 - \smallsCeVx_{\vertexb} && \oneCe \in \bigsCe, \vertexb \in \nddSet \label{FailedNDD}\\
    \sum_{\match \in \cyclesetv} \cyclevarCe - \sum_{\ell \in \idposset} \sum_{\parc \in \arcSetL} \chainvarCe \le 1 - \smallsCeVx_{\vertexb} && \oneCe \in \bigsCe, \vertexb \in \pairSet \label{FailedPair}\\
    \sum_{\match \in \cyclesetuv} \cyclevarCe - \sum_{\ell \in \posset} \chainvarCe \le 1 - \smallsCeArc_{\narc} && \oneCe \in \bigsCe, \parc \in \arcSet \label{FailedArc}\\
    \sum_{\vertexa:\parc \in \arcSetL} \chainvarCe -  \sum_{\vertexa:\parcrev \in \arcSetL} \chainvarNextCe \le 0 && \oneCe \in \bigsCe, \vertexb \in \pairSet, \ell \in \idposset \setminus \{\chainCap - 1\} \label{ChainPos}\\
    \sum_{\vertexa: \parcrev}\chainvarR \le 1 && \vertexb \in \nddSet \label{FailedNDDf}\\
    \sum_{\match \in \cyclesetv} \cyclevar - \sum_{\ell \in \idposset} \sum_{\parc \in \arcSetL} \chainvar \le 1 && \vertexb \in \pairSet\\
    \sum_{\vertexa:\parc \in \arcSetL} \chainvar -  \sum_{\vertexa:\parcrev \in \arcSetL} \chainvarNext \le 0 && \oneCe \in \bigsCe, \vertexb \in \pairSet, \ell \in \idposset \setminus \{\chainCap - 1\} \label{ChainPosf}\\
    \vertxCevar \ge 0 && \oneCe \in \bigsCe, \vertexb \in \pairSet\\
    \varRecoZ_\match, \cyclevarCe \in \{0,1\} && \oneCe \in \bigsCe, \match \in \cycleSet\\
    \chainvar, \chainvarCe \in \{0,1\} && \oneCe \in \bigsCe, \parc \in \arcSet, \ell \in \posset \label{varNature}
\end{align}
\label{for:PICEF_RO}
\end{subequations}

Constraints \eqref{selidxF} are used to determine the scenario binding the number of patients that receive a transplant in both stages. Such scenario is the worst-case scenario. The objective \eqref{objmaxRO} is then equivalent to the maximum number of patients from the first stage that can be recovered in the second stage under the worst-case scenario.  Constraints \eqref{FirstSSol} and \eqref{SecondSSol} assure that a pair $\vertexb$ is counted as recovered in the objective if it is selected in the first-stage solution (Constraints \eqref{FirstSSol}) and it is also selected in the second stage under scenario $\oneCe \in \bigsCe$ (Constraints \eqref{SecondSSol}). Constraints \eqref{FailedNDD} to Constraints \eqref{FailedArc} guarantee that the solution obtained for every scenario $\oneCe \in \bigsCe$ is a matching. Specifically, Constraints \eqref{FailedNDD} assure that if a non-directed donor fails under some scenario $\oneCe \in \bigsCe$, i.e., $\smallsCeVx_{\vertexb} = 1$ then its corresponding arcs in position one cannot be used to trigger a chain. Similarly, Constraints \eqref{FailedPair} guarantee that if a pair fails, then it cannot be present in either a cycle nor a chain. Constraints \eqref{FailedArc} ensure that when an arc $\parc \in \arcSet$ fails under some scenario $\oneCe \in \bigsCe$, it does not get involved in either a cycle or a chain. Constraints \eqref{ChainPos} assure the continuity of a chain by selecting arcs in consecutive positions. Constraints \eqref{FailedNDDf} to Constraints \eqref{ChainPosf} select a solution corresponding to a matching in the first stage. The remaining constraints correspond to the nature of the decision variables.

\section{Robust formulation for first-stage-only recourse}
\label{sec:fisrtOnlyROFormulation}
In addition to the constraints defining formulation \eqref{for:PICEF_RO}, a new one is introduced to limit the recourse solutions to include only vertices that were selected in the first stage under every scenario. 
\begin{subequations}
\label{FirstStageOnlyFormulation}
\begin{align}
    \max \qquad &\objPosIdx \\
    \eqref{selidxF} &-\eqref{varNature}&\\
    \sum_{\match \in \cyclesetv} \cyclevarCe - \sum_{\ell \in \idposset} \sum_{\parc \in \arcSetL} \chainvarCe &\le \sum_{\match \in \cyclesetv} \varRecoZ_\match - \sum_{\ell \in \idposset} \sum_{\parc \in \arcSetL} \chainvar & \oneCe \in \bigsCe, \vertexb \in \pairSet
\end{align}
\end{subequations}

\section{Failure scenarios generated by \ref{MasterBasic} and \ref{MasterExp}}
\label{Example}

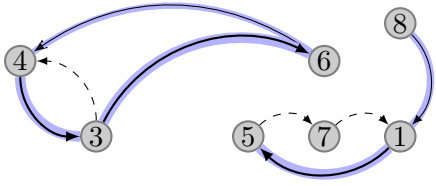
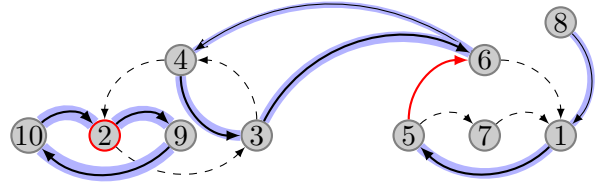
\begin{figure}[tbp]
	\centering	
	\begin{adjustbox}{minipage=\linewidth,scale=1}
		\begin{subfigure}[t]{.48\textwidth}
			\centering
    	    \tikzstyle{place}=[circle,draw=blue!50,fill=blue!20,thick,
			inner sep=0pt,minimum size=4mm]
			\tikzstyle{transition}=[circle,draw=black!50,fill=black!20,thick,
			inner sep=0pt,minimum size=4mm]
			\tikzstyle{dummy}=[rectangle,draw=white,fill=white,thick,
			inner sep=0pt,minimum size=4mm]
			\begin{tikzpicture}
				
				\node[transition] (cuatro) at (0,2) {$4$}; 
				\node[transition] (tres) at (1,1) {$3$};
				
				\node[transition] (cinco) at (3,1) {$5$};
				\node[transition] (uno) at (5,1)  {$1$}; 
				\node[transition] (ocho) at (5,2.5)  {$8$}; 
				\node[transition] (seis) at (4,2)  {$6$}; 
				\node[transition] (siete) at (4,1)  {$7$}; 
				
				
				\draw [-{latex}, thick] (cuatro) to [bend right=45] (tres);
				\draw [-{latex},dashed] (tres) to [bend right=45] (cuatro);
				
				\draw [-{latex}, thick] (uno) to [bend left=45] (cinco);
				\draw [-{latex},dashed] (ocho) to [bend left=45] (uno);
				\draw [-{latex},dashed] (siete) to [bend left=45] (uno);
				\draw [-{latex},dashed] (cinco) to [bend left=45] (siete);
				
				\draw [-{latex}, thick] (tres) to [bend left=45] (seis);
				\draw [-{latex}, thick] (seis) to [bend right=35] (cuatro);
				
				\draw[-latex,preaction = {draw,blue!30,-,double=blue!30,double distance = 5\pgflinewidth}] (ocho) to [bend left=45] (uno);
				\draw[thick, -latex,preaction = {draw,blue!30,-,double=blue!30,double distance = 3\pgflinewidth}] (uno) to [bend left=45] (cinco);
				\draw[thick, -latex,preaction = {draw,blue!30,-,double=blue!30,double distance = 3\pgflinewidth}] (cuatro) to [bend right=45] (tres);
				\draw[-latex,preaction = {draw,blue!30,-,double=blue!30,double distance = 4\pgflinewidth}] (seis) to [bend right=35] (cuatro);
				\draw[thick, -latex,preaction = {draw,blue!30,-,double=blue!30,double distance = 3\pgflinewidth}] (tres) to [bend left=45] (seis);
			\end{tikzpicture}
			\caption{Second-stage compatibility graph at iteration 2.}
			\label{fig:IteOne-ScndSGraph}
		\end{subfigure}
            \hspace{0.5cm}
		\begin{subfigure}[t]{.48\textwidth}
			\centering
    	    \tikzstyle{place}=[circle,draw=blue!50,fill=blue!20,thick,
			inner sep=0pt,minimum size=4mm]
			\tikzstyle{transition}=[circle,draw=black!50,fill=black!20,thick,
			inner sep=0pt,minimum size=4mm]
			\tikzstyle{dummy}=[rectangle,draw=white,fill=white,thick,
			inner sep=0pt,minimum size=4mm]
			\begin{tikzpicture}
				
				\node[transition] (cuatro) at (0,2) {$4$}; 
				\node[transition] (tres) at (1,1) {$3$};
				\node[transition, draw=red] (dos) at (-1,1) {$2$};
				\node[transition] (nueve) at (0,1) {$9$};
				\node[transition] (diez) at (-2,1) {$10$};
				
				\node[transition] (cinco) at (3,1) {$5$};
				\node[transition] (uno) at (5,1)  {$1$}; 
				\node[transition] (ocho) at (5,2.5)  {$8$}; 
				\node[transition] (seis) at (4,2)  {$6$}; 
				\node[transition] (siete) at (4,1)  {$7$}; 
				
				\draw [-{latex},dashed] (dos) to [bend left=45] (nueve);
				\draw [-{latex},dashed] (nueve) to [bend left=45] (diez);
				\draw [-{latex},dashed] (diez) to [bend left=45] (dos);
				
				\draw [-{latex}, dashed] (cuatro) to [bend right=45] (dos);
				\draw [-{latex}, dashed] (dos) to [bend right=45] (tres);
				\draw [-{latex}, thick] (cuatro) to [bend right=45] (tres);
				\draw [-{latex}, dashed] (tres) to [bend right=45] (cuatro);
				
				\draw [-{latex}, dashed] (seis) to [bend left=45] (uno);
				\draw [-{latex}, thick] (uno) to [bend left=45] (cinco);
				\draw [-{latex}, dashed] (ocho) to [bend left=45] (uno);
				\draw [-{latex}, thick, draw=red] (cinco) to [bend left=45] (seis);
				\draw [-{latex},dashed] (siete) to [bend left=45] (uno);
				\draw [-{latex},dashed] (cinco) to [bend left=45] (siete);
				
				\draw [-{latex}, thick] (tres) to [bend left=45] (seis);
				\draw [-{latex}, thick] (seis) to [bend right=35] (cuatro);
				
				\draw[-latex,preaction = {draw,blue!30,-,double=blue!30,double distance = 5\pgflinewidth}] (ocho) to [bend left=45] (uno);
				\draw[thick, -latex,preaction = {draw,blue!30,-,double=blue!30,double distance = 3\pgflinewidth}] (uno) to [bend left=45] (cinco);
				\draw[thick, -latex,preaction = {draw,blue!30,-,double=blue!30,double distance = 3\pgflinewidth}] (cuatro) to [bend right=45] (tres);
				\draw[-latex,preaction = {draw,blue!30,-,double=blue!30,double distance = 4\pgflinewidth}] (seis) to [bend right=35] (cuatro);
				\draw[thick, -latex,preaction = {draw,blue!30,-,double=blue!30,double distance = 3\pgflinewidth}] (tres) to [bend left=45] (seis);

                \draw[thick, -latex,preaction = {draw,blue!30,-,double=blue!30,double distance = 5\pgflinewidth}] (diez) to [bend left=45] (dos);
                \draw[thick, -latex,preaction = {draw,blue!30,-,double=blue!30,double distance = 5\pgflinewidth}] (dos) to [bend left=45] (nueve);
                \draw[thick, -latex,preaction = {draw,blue!30,-,double=blue!30,double distance = 5\pgflinewidth}] (nueve) to [bend left=45] (diez);
				
			\end{tikzpicture}
			\caption{Transitory graph at iteration 2. Elements in red fail at iteration 1.}
			\label{fig:IteOne-TrGraph}
		\end{subfigure}
	\end{adjustbox}
	\caption{Optimal recourse solutions at iteration 2 when vertex 2 and arc (5,6) in the first-stage compatibility graph fail at iteration 1. Shaded arcs are the optimal recourse solution. Dashed arcs exist but are not selected. }
	\label{fig:MS-MT-example}
\end{figure}

The following example demonstrates that \ref{MasterExp} can yield failure scenarios that dominate those in \ref{MasterBasic}. Consider Figure \ref{fig:RecourseExamples} again under full recourse, $\vBudget = 1$, $\aBudget = 1$ and the first-stage solution $\varFirstc \in \setFirst$ involving eight pairs in that figure. Observe, the first-stage compatibility graph also corresponds to the transitory graph $\digraphSndTrc$. At iteration 1, we observe failure scenario $\tildesmallsCeVx_{2} = 1$ and $\tildesmallsCeArc_{56} = 1$, thus, vertices 2, 9 and 10 do not belong to the realization of the second-stage compatibility graph, $\digraphSndSc$ (Figure \ref{fig:IteOne-ScndSGraph}). The optimal objective value to the recourse problem in $\digraphSndSc$, i.e., $\RecoPcdte$, is $\RecoObjVarcdte = 5$ with an optimal solution $\varRecoOpSndSc$ involving cycle (3,4,6) and chain (8,1,5). When we solve the recourse problem in the transitory graph $\digraphSndTrc$ (Figure \ref{fig:IteOne-TrGraph}), the optimal objective value is also 5, but its optimal recourse solution $\VarRecoOpExp = \varRecoOpSndSc \cup \{(2,9,10)\}$ involves cycle (3,4,6), chain (8,1,5), as well as the failed cycle (2,9,10). At iteration 2, \ref{MasterBasic} and \ref{MasterExp} attempt to find a new failure scenario, but the failure scenario that is feasible to \ref{MasterBasic} is infeasible to \ref{MasterExp}. Below, the first two iterations of \ref{MasterBasic} and \ref{MasterExp} are shown. The vertices/arcs in the failure scenario of every iteration are indicated within a box.

\begin{tabular}{p{0.45\textwidth} @{\quad} | @{\quad} p{0.45\textwidth}}
$\text{MasterSecond}(\varFirstc)$ (\ref{MasterBasic}): & $\text{MasterTransitory}(\varFirstc)$ (\ref{MasterExp}):\\
\hline
\uline{Iteration 1}: \ $\UBCovering = 8$\newline
$\boxed{\smallsCeVx_{2}}    + \smallsCeVx_{9}     + \smallsCeVx_{10} +
    \smallsCeArc_{2,9} + \smallsCeArc_{9,10} + \smallsCeArc_{10,2} + \smallsCeVx_{3}    +  \smallsCeVx_{4}    + 
    \smallsCeArc_{3,4} \newline
    \text{\quad} + \smallsCeArc_{4,3}     +
    \smallsCeVx_{1}    + \smallsCeVx_{5}     + \smallsCeVx_{6} +
    \smallsCeVx_{1,5}    + \boxed{\smallsCeVx_{5,6}} +
    \smallsCeVx_{6,1}
    \ge 1$ 
    \newline\newline
\uline{Iteration 2}: \ $\UBCovering = 5$ \newline
    $\smallsCeVx_{2}    + \smallsCeVx_{9}     + \smallsCeVx_{10} +
    \smallsCeArc_{2,9} + \smallsCeArc_{9,10} + \smallsCeArc_{10,2} + \boxed{\smallsCeVx_{3}}     +  \smallsCeVx_{4}    + 
    \smallsCeArc_{3,4} \newline
    \text{\quad} + \smallsCeArc_{4,3}     +
    \smallsCeVx_{1}    + \smallsCeVx_{5}     + \smallsCeVx_{6} +
    \smallsCeVx_{1,5}    + \boxed{ \smallsCeVx_{5,6}} +
    \smallsCeVx_{6,1}
    \ge \cancel{1} \rightarrow 2 \text{ by } \ref{prop:RHSi} \newline
    \boxed{\smallsCeVx_{3}}    + \smallsCeVx_{6}     + \smallsCeVx_{4} +
    \smallsCeArc_{3,6} + \smallsCeArc_{6,4}  + \smallsCeArc_{4,3} + \smallsCeVx_{1}    +  \smallsCeVx_{5}    + \smallsCeVx_{8} \newline
    \text{\quad} + \smallsCeArc_{1,5}    + \smallsCeArc_{8,1} 
    \ge 1 $
&
\uline{Iteration 1}: \ $\UBCovering = 8$\newline
    $\boxed{\smallsCeVx_{2}}    + \smallsCeVx_{9}     + \smallsCeVx_{10} +
    \smallsCeArc_{2,9} + \smallsCeArc_{9,10} + \smallsCeArc_{10,2} + \smallsCeVx_{3}    +  \smallsCeVx_{4} \newline
    \text{\quad} + \smallsCeArc_{3,4}    + \smallsCeArc_{4,3}     +
    \smallsCeVx_{1}    + \smallsCeVx_{5}     + \smallsCeVx_{6} +
    \smallsCeVx_{1,5}    + \boxed{ \smallsCeVx_{5,6}} +
    \smallsCeVx_{6,1}
    \ge 1$ 
    \newline\newline
\uline{Iteration 2}: \ $\UBCovering = 5$\newline
    $\smallsCeVx_{2}    + \smallsCeVx_{9}     + \smallsCeVx_{10} +
    \smallsCeArc_{2,9} + \boxed{\smallsCeArc_{9,10}} + \smallsCeArc_{10,2} + \boxed{\smallsCeVx_{3}}     +  \smallsCeVx_{4}    + 
    \smallsCeArc_{3,4}    \newline
    \text{\quad} + \smallsCeArc_{4,3}     +
    \smallsCeVx_{1}    + \smallsCeVx_{5}     + \smallsCeVx_{6} +
    \smallsCeVx_{1,5}    +  \smallsCeVx_{5,6} +
    \smallsCeVx_{6,1}
    \ge \cancel{1} \rightarrow 2 \text{ by } \ref{prop:RHSi} \newline
    \boxed{\smallsCeVx_{3}}    + \smallsCeVx_{6}     + \smallsCeVx_{4} +
    \smallsCeArc_{3,6} + \smallsCeArc_{6,4}  + \smallsCeArc_{4,3} + \smallsCeVx_{1}    +  \smallsCeVx_{5}    + \smallsCeVx_{8} \newline
    \text{\quad} + \smallsCeArc_{1,5}    + \smallsCeArc_{8,1} \ge 1 \newline
    \boxed{\smallsCeVx_{3}}    + \smallsCeVx_{6}     + \smallsCeVx_{4} +
    \smallsCeArc_{3,6} + \smallsCeArc_{6,4}  + \smallsCeArc_{4,3} + \smallsCeVx_{1}    +  \smallsCeVx_{5}    + \smallsCeVx_{8} \newline
    \text{\quad} +   \smallsCeArc_{1,5}    + \smallsCeArc_{8,1} +
    \smallsCeVx_{2}    +  \smallsCeVx_{9}    + \smallsCeVx_{10} + \smallsCeArc_{2,9}   + \boxed{\smallsCeArc_{9,10}} \newline
    \text{\quad} + \smallsCeArc_{10,2} \ge \cancel{1} \rightarrow  2 \text{ by } \ref{prop:RHSExtdi}$
\end{tabular} \\

The new failure scenario for \ref{MasterBasic} would lead to an optimal recourse solution with again $\RecoObjVarcdte = 5$, whereas the optimal recourse solution for \ref{MasterExp} would lead to $\RecoObjVarcdte = 3$. That is, in the third iteration, the right-hand side of constraints \ref{MasterExp} will get updated with a new upper bound $\UBCovering = 3$, leading \ref{MasterExp} to infeasibility sooner and thus, to the optimal objective value of the second-stage problem. On the other hand, \ref{MasterBasic} in the third iteration will need another attempt to ``discover'' a failure scenario that will bring down the ceiling of $\UBCovering$.

\section{The recourse problem}
\label{Alg:CyChrecourse}

We present in this section two algorithms to enumerate the feasible cycles and chains in $\cycleSet^{\policy}(\varFirstc)$ and $\chainSet^{\policy}(\varFirstc)$, respectively, that lead to a transitory graph $\digraphSndTrc$ and a realization of the second-stage compatibility graph $\digraphSndSc$ under scenario $\cdteOne \in \bigsCec$. 

\subsection{Cycles and chains for the full-recourse policy}

We start by defining some notation. Let $\SelPairsX = \SelVerXc \cap \pairSet$ be an auxiliary set, corresponding to the pairs selected in the first-stage solution. Moreover, we denote by $\cyclechainSet_{\cycleCap}^{\vertexa}$  the set of feasible \textit{cycles} that include vertex $\vertexa \in \SelPairsX$. Lastly, consider $R$ as an auxiliary vertex set. Algorithm \ref{Alg:FullRecoCycleSet} iteratively builds $\cycleSet^{\text{Full}}(\varFirstc)$. At the start, $R$ and $\cycleSet^{\text{Full}}(\varFirstc)$ are empty. Within the While loop, a pair $\vertex$ from the first-stage solution is selected. Then, in line 4, a deep search procedure, starting from vertex $\vertex$, is used to find $\cyclechainSet_{\cycleCap}^{\vertexa}$ in graph $\tilde{D} = (\vertexSet \setminus R, \arcSet)$. Note that if $R \neq \emptyset$, the new cycles are found in a graph where the previously selected vertices (the ones in set $R$) are removed, since otherwise, cycles already in $\cyclechainSet_{\cycleCap}^{\vertexa}$ could be found again. The new cycles are then added to $\cycleSet^{\text{Full}}(\varFirstc)$ and vertex $\vertex$ is removed from $\pairSet(\varFirstc)$. When no more vertices are left in that set, the algorithm ends.

\begin{algorithm}[tbp]
	\algorithmicrequire { Set with pairs from the first-stage solution, $\SelPairsX$}\\
	\algorithmicensure { Set $\cycleSet^{\text{Full}}(\varFirstc)$}
    \begin{algorithmic}[1]
    \STATEnonum \textbf{Step 0: }
	\STATE  $R = \emptyset$; $\cycleSet^{\text{Full}}(\varFirstc) = \emptyset$  \\ 
	\STATEnonum \textbf{Step 1: }
	\WHILE{$\SelPairsX \neq \emptyset$}
	\STATE \text{Select vertex} $\vertexa \in \SelPairsX $
	\STATE \text{Find} $\cyclechainSet_{\cycleCap}^{\vertexa}$ \text{in graph } $\tilde{D} = (\vertexSet \setminus R, \arcSet)$ \text{from vertex} $\vertex$
	\STATE $\cycleSet^{\text{Full}}(\varFirstc) \gets \cycleSet^{\text{Full}}(\varFirstc) \cup \cyclechainSet_{\cycleCap}^{\vertexa}$
	\STATE $\SelPairsX  \gets \SelPairsX  \setminus \{\vertexa\}$
	\STATE $R \gets R \cup \{\vertexa\}$
	\ENDWHILE
	\STATEnonum \textbf{Step 2: }
	\STATE Return $\cycleSet^{\text{Full}}(\varFirstc)$
	\caption{Obtaining $\cycleSet^{\policy}(\varFirstc)$ with $\policy = \text{Full}$}
	\label{Alg:FullRecoCycleSet}
	\end{algorithmic}
\end{algorithm}
The correctness of Algorithm \ref{Alg:FullRecoCycleSet} follows from the fact that only cycles/chains with at least one vertex from the first stage contribute to the weight of a recourse solution. Thus, it suffices to find such cycles per every vertex in $\SelPairsX$. Similar reasoning is used when finding chains. To this end, we include additional notation. Let $\tilde{\nddSet} = \nddSet$ be an auxiliary set that corresponds to the set of non-directed donors. Moreover, let $\cyclechainSet_{\ell}^{\vertexa}$ be the set of chains that include at least one pair in $\SelPairsX$ and are triggered by vertex $\vertexa \in \nddSet$ with exactly $1 \le \ell \le \chainCap$ arcs. Algorithm \ref{Alg:FullRecoChainSet} iteratively builds the set of chains that include at least one pair that is selected in the first stage.

\begin{algorithm}[tbp]
	\algorithmicrequire { A first-stage solution $\varFirstc \in \setFirst$}\\
	\algorithmicensure { Set $\chainSet^{\text{Full}}(\varFirstc)$}
    \begin{algorithmic}[1]
    \STATEnonum \textbf{Step 0: }
	\STATE  $R = \emptyset$; $\chainSet^{\text{Full}}(\varFirstc) = \emptyset$  \\ 
	\STATEnonum \textbf{Step 1: }
	\WHILE{$\tilde{\nddSet} \neq \emptyset$}
    	\STATE \text{Select vertex} $\vertexa \in \tilde{\nddSet}$
    	\FORALL{$ 1 \le \ell \le \chainCap$}
        	\STATE \text{Find} $\cyclechainSet_{\ell}^{\vertexa}$ \text{in graph } $\graphKEP$ \text{from vertex} $\vertex$
        	\STATE $\chainSet^{\text{Full}}(\varFirstc) \gets \chainSet^{\text{Full}}(\varFirstc) \cup \cyclechainSet_{\ell}^{\vertexa}$
        \ENDFOR
        \STATE $\tilde{\nddSet} \gets \tilde{\nddSet} \setminus \{\vertexa\}$
	\ENDWHILE
	\STATEnonum \textbf{Step 2: }
	\STATE Return $\chainSet^{\text{Full}}(\varFirstc)$
	\caption{Obtaining $\chainSet^{\policy}(\varFirstc)$ with $\policy = \text{Full}$}
	\label{Alg:FullRecoChainSet}
	\end{algorithmic}
\end{algorithm}

\subsection{Cycles and chains for the first-stage-only recourse policy}

Algorithm \ref{Alg:FullRecoCycleSet} can be modified to accommodate the first-stage-only recourse for cycles. Specifically, we can replace $\cycleSet^{\text{Full}}(\varFirstc)$ by $\cycleSet^{\text{1stSO}}(\varFirstc)$, where the latter is the set of simple cycles with at least one vertex in $\SelPairsX$ satisfying the first-stage-only recourse policy. In line 4, the vertex set $\vertexSet$ in graph $\tilde{D}$ is replaced by $\SelPairsX$ so that only pairs from the first stage can be part of the allowed cycles. The arc set of $\tilde{D}$ can then be defined such that every arc in it has a starting and terminal vertex in $\SelPairsX$. Likewise, Algorithm \ref{Alg:FullRecoChainSet} can also support chains for the first-stage-only recourse. The only change in Algorithm \ref{Alg:FullRecoChainSet} is to replace graph $\graphKEP$ by graph $\tilde{\graph} = (\SelVerXc, \tilde{A})$ where every arc in the arc set $\tilde{A}$ has both extreme vertices in $\SelVerXc$.

The algorithms just described are used to obtain the cycles and chains that can participate in a recourse solution when the recourse problem is solved as a sub-problem in the robust decomposition presented in the next section.

\subsection{Formulations for the recourse problem}
\label{sec:RecourseFormls}
In this section we present the cycle-and-chain MIP formulations presented in \citet{Blom2021} adapted to the problem we study. However, we note that a MIP formulation for the recourse problem does not require the explicit enumeration of cycles and chains. The advantage of enumeration is that different policies can be easily addressed, specifically, the full recourse and the first-stage-only recourse.

 Recall that $\cycleSetSndSc \cup \chainSetSndSc \subseteq \cycleSetSndTrc \cup \chainSetSndTrc$ is the set of non-failed cycles/chains in a second-stage compatibility graph under scenario $\cdteOne \in \bigsCe$ and policy $\policy \in \policySet$, i.e., $\sum_{\vertexa \in V(\match)} \smallsCeVx_{\vertexa} + \sum_{\parc \in A(\match)} \smallsCeArc_{\narc} = 0 \ \forall \match \in \cycleSetSndSc \cup \chainSetSndSc$. We present the recourse problem based on the so-called cycle formulation \citep{Abraham2007} as follows:
\begin{subequations}
	\label{RecoPModel}
    \begin{align}
    	\RecoP:  \max_{y} \quad&  \sum_{\match \in \cycleSetSndSc \cup \chainSetSndSc}\wTwoStageGralc y_{\match} \label{RecoP:Obj} \tag{R}\\
    	& \sum_{\match: \vertexa \in V(\match)} y_{\match} \le  1 &  \vertexa \in \vertexSet \label{eq:uniqueC}\\
    	&y_{\match} \in \{0,1\} & \match \in \cycleSetSndSc \cup \chainSetSndSc
    \end{align} 
\end{subequations}

\noindent Here, $y_{\match}$ is a decision variable for every cycle/chain $\match \in \cycleSetSndSc \cup \chainSetSndSc$, taking on value one if selected, and zero otherwise. Constraints \eqref{eq:uniqueC} ensure that a vertex belongs to at most one cycle/chain. An optimal solution to formulation $\RecoPcdte$ finds a matching with the greatest number of matched pairs from the first stage after observing failure scenario $\cdteOne$. Thus, by solving formulation \eqref{RecoPModel}, a new Constraint \eqref{eq:solsTwo} can be created.  \cite{Blom2021} proposed to expand recourse solutions by including, in addition to the optimal non-failed cycles/chains found by formulation \eqref{RecoPModel}, failed cycles and chains while guaranteeing that the solution is still a matching. Although an expanded solution does not contribute to more recourse value under the failure scenario in consideration, it may imply other violated constraints. Consider two recourse solutions $\VarRecoOpExp \in \setRecoTrcdte$ and $\varRecoOpSndSc \in \setRecocdte$, one in the transitory compatibility graph and another one in the second-stage compatibility graph, respectively, and also assume that $\varRecoOpSndSc \subseteq \VarRecoOpExp$. Then, constraint \eqref{eq:solsTwo} associated to $\varRecoOpSndSc$ is directly implied by that of $\VarRecoOpExp$, i.e., if the constraint corresponding to $\varRecoOpSndSc$ is violated, so is the constraint corresponding to $\VarRecoOpExp$. We find expanded recourse solutions by solving a deterministic KEP in the transitory graph $\digraphSndTrc$ instead of the second-stage graph $\digraphSndSc$ in formulation \eqref{RecoPModel}, and by assigning new weights to all cycles/chains $\match \in \cycleSetSndTrc \cup \chainSetSndTrc$ as in \citet{Blom2021}. The new weights, assigned to each cycle $\match \in \cycleSetSndTrc \cup \chainSetSndTrc$, are given by
\begin{align}
\label{expRecoObFcn}
\wTwoStageExpanded = \left\{
\begin{array}{ll}
\wTwoStageGralc \lvert \vertexSet \rvert + 1      & \mbox{if } V(\match) \subseteq \vertexSet \setminus \vertexSet_{\oneCe}\\
1      & \mbox{otherwise }
\end{array}
\right.
\end{align}
We denote by $\RecoPexpc$ the resulting recourse problem with expanded solutions,
\begin{subequations}
	\label{RecoPModelexp}
    \begin{align}
    	\RecoPexpc:  \max_{y} \quad&  \sum_{\match \in \cycleSetSndTr \cup \chainSetSndTr}\wTwoStageExpanded y_{\match} \label{RecoP:Objexp} \tag{RE}\\
    	& \sum_{\match: \vertexa \in V(\match)} y_{\match} \le  1 &  \vertexa \in \vertexSet \label{eq:uniqueCexp}\\
    	&y_{\match} \in \{0,1\} & \match \in \cycleSetSndTrc \cup \chainSetSndTrc
    \end{align} 
\end{subequations}
The following lemma, based on \citet{Blom2021}, states that the set of cycles/chains that do not fail in a recourse solution optimal to $\RecoPexpc$ is an optimal recourse solution to the original recourse problem $\RecoPcdte$. 
\begin{lemma}
\label{lemma:exp}
For a recourse solution $\VarRecoOpExp \in \setRecoTrcdte$ that is optimal to $\RecoPexpc$, its set of non-failed cycles/chains, $\VarRecoOpExpcdte \subseteq \VarRecoOpExp$ under scenario $\cdteOne \in \bigsCec$ is an optimal recourse solution to $\RecoPcdte$.
\end{lemma}
\proof
Let $\RecoExpObjVar$ be the optimal objective value to solution $\VarRecoOpExp \in \setRecoTrcdte$ and consider $\vert \vertexSet \rvert/2$ as the maximum number of cycles/chains that can originate in a feasible matching, i.e., the maximum number of decision variables for which $\VarRecoOpExp_{\match} = 1 \ \forall \match \in \cycleSetSndTrc \cup \chainSetSndTrc$. Moreover, let $\mathbbm{1}_{\match}$ be an indicator variable that takes on value one if $\VarRecoOpExp_{\match} = 1 \ \forall \match \notin \cycleSetSndSc \cup \chainSetSndSc$ and zero otherwise. Then, from \eqref{expRecoObFcn} we know that
\begin{subequations}
    \begin{align}
        \RecoExpObjVar = & \sum_{\match \in \cycleSetSndSc \cup \chainSetSndSc: \VarRecoOpExpcdte_{\match} = 1} \left( \wTwoStageGralc \lvert \vertexSet \rvert + 1 \right) + \sum_{\match \notin \cycleSetSndSc \cup \chainSetSndSc: \VarRecoOpExp_{\match} = 1} \mathbbm{1}_{\match}
        \\
        =&\lvert \vertexSet \rvert \sum_{\match \in \cycleSetSndSc \cup \chainSetSndSc: \VarRecoOpExpcdte_{\match} = 1} \wTwoStageGralc  + \sum_{i = 1}^{\vert \vertexSet \rvert/2} 1 + \sum_{i = 1}^{\vert \vertexSet \rvert/2} 1\\
        =&z^{\star}\lvert \vertexSet \rvert  + \vert \vertexSet \rvert \label{result}
    \end{align}
\end{subequations}

Now, suppose $\VarRecoOpExpcclt$ is not optimal to $\RecoPcdte$, therefore $\sum_{\match \in \cycleSetSndSc \cup \chainSetSndSc: \VarRecoOpExp_{\match} = 1} \wTwoStageGralc$ can be at most $z^{\star} - 1$. If that is true, then, when replacing $z^{\star}$ in \eqref{result} by $z^{\star} - 1$, we obtain that  $\RecoExpObjVar = z^{\star}\lvert \vertexSet \rvert  + \vert \vertexSet \rvert = z^{\star}\lvert \vertexSet \rvert$, which is a contradiction. Thus, it follows that $\VarRecoOpExpcclt$ must be optimal to $\RecoPcdte$. \hfill$\square$

It is worth noting that we formulated the recourse problem by means of cycle-and-chain decision variables, but since it is reduced to a deterministic KEP, the recourse problem can be solved via multiple formulations/algorithms from the literature, e.g., \citet{Omer2022, Blom2021, Riascos2020}.

\section{Additional results for non-homogeneous failure}
\label{CycleLengthFourNonHomoFail}
The data  in Table \ref{tab:policiesk4} includes all runs, optimal or not, except for column HSP(\%) where only instances that were solved to optimality are included. The interpretation of the displayed columns correspond to those in Table \ref{tab:policiesk3}.

\begin{table}[htpb]
    \caption{Policies comparison for $\cycleCap = 4$ under HSA\_MB and non-homogeneous failure.}
    \label{tab:policiesk4}
    \begin{adjustbox}{width=\textwidth}
        \centering
        \begin{tabular}{rrrrrrrrrrrrrrrrrrr}
    \toprule
    &  &  &  \multicolumn{16}{c}{Full recourse}\\
    \cmidrule(lr){4-19}
     &  &  &  &   & \multicolumn{6}{c}{Time} & \multicolumn{7}{c}{Avg. total \#}\\
    \cmidrule(lr){6-11} \cmidrule(lr){12-19}
    $\chainCap$   & $\vBudget$   &  $\aBudget$  &      opt &      HSP (\%) &   total     & Alg\ref{Alg:Heuristics} (\%) &   \ref{MasterExp} (\%) &   \ref{RecoP:Objexp} (\%) &   CG (\%)         &\ref{objMPOptQ} (\%)  & 1stS      &  2ndS         &   Alg\ref{Alg:Heuristics}-true  &   \ref{MasterExp} &   \ref{RecoP:Objexp} &      CGtrue &   \ref{objMPOptQ} &   DomS           \\
       \cmidrule(lr){1-19}
          $\lvert\vertexSet\rvert = 50$&             &             &         &           &          &              &                 &                  &                 &               &              &              &              &                &                 &                &              &          \\
               3 &           1 &        1.83 &      30 &   35.98 &      34.86 &         1.99 &           26.29 &             1.67 &           35.91 &         15.72 &         2.77 &       339.50  &       111.70  &         205.20  &            3.63 &         335.87 &        22.60  &         1.43 \\
               3 &           2 &        1.83 &      30 &   29.82 &     186.00    &         1.82 &           16.56 &             0.84 &           21.82 &         49.24 &         3.47 &       561.27 &       260.93 &         218.13 &            5.33 &         555.93 &        82.20  &         1.17 \\
               3 &           1 &        3.03 &      30 &   29.17 &     276.87 &         1.76 &           19.79 &             0.76 &           23.55 &         41.61 &         3.53 &       548.43 &       196.37 &         278.90  &            4.00    &         544.43 &        73.17 &         1.03 \\
               3 &           2 &        3.03 &      23 &   21.18 &     939.48 &         1.81 &           17.07 &             0.77 &           18.88 &         52.71 &         2.43 &       430.67 &       169.97 &         168.93 &            4.17 &         426.50  &        91.77 &         3.00    \\
               4 &           1 &        1.83 &      30 &   35.87 &      84.48 &         1.62 &           14.41 &             1.90  &           49.89 &          9.59 &         3.23 &       438.30  &       165.23 &         245.10  &            7.07 &         431.23 &        27.97 &         1.30  \\
               4 &           2 &        1.83 &      30 &   27.18 &     134.60  &         1.42 &           11.11 &             1.55 &           37.20  &         31.59 &         3.03 &       477.07 &       225.07 &         194.87 &            5.93 &         471.13 &        57.13 &         0.67 \\
               4 &           1 &        3.03 &      29 &   30.96 &     239.39 &         1.23 &           11.49 &             1.14 &           32.29 &         37.17 &         3.53 &       572.13 &       212.87 &         279.93 &            4.50  &         567.63 &        79.33 &         2.63 \\
               4 &           2 &        3.03 &      25 &   25.29 &     859.72 &         1.08 &           10.97 &             0.50  &           23.34 &         52.00    &         2.73 &       477.30  &       214.37 &         169.43 &            3.93 &         473.37 &        93.50  &         2.20  \\
       \cmidrule(lr){1-19}
       $\lvert\vertexSet \rvert =  100$&                &           &          &           &          &              &                 &                  &                 &               &              &              &              &                &                 &                &              &          \\
              3 &           1 &        3.30  &      13 &   48.74  &    2529.89 &         0.19 &            0.72 &             0.84 &           22.99 &         69.94 &         2.93 &       630.63 &       268.40  &         169.10  &           18.00    &         612.63 &       193.13 &         1.03 \\
              3 &           2 &        3.30  &       4 &   38.61 &    3194.79 &         0.09 &            0.33 &             0.22 &            9.56 &         85.78 &         1.70  &       344.23 &       226.43 &          27.70  &            7.43 &         336.80  &        90.10  &         0.83 \\
              3 &           1 &        5.93 &       1 &   12.50 &    3598.01 &         0.03 &            0.06 &             0.19 &            4.96 &         91.89 &         1.37 &       297.00    &       154.10  &          49.93 &            5.57 &         291.43 &        92.97 &         1.13 \\
              3 &           2 &        5.93 &       0 &   - &    3600.02 &         0.03 &            0.01 &             0.18 &            4.11 &         87.98 &         1.20  &       226.93 &       168.63 &          10.37 &            4.50  &         222.43 &        47.93 &         1.17 \\
              4 &           1 &        3.30  &       7 &   45.27 &    2934.19 &         0.10  &            0.51 &             5.66 &           52.70  &         31.00    &         2.40  &       448.87 &       190.20  &         143.90  &           11.70  &         437.17 &       114.77 &         0.30  \\
              4 &           2 &        3.30  &       5 &   33.28 &    3285.81 &         0.05 &            0.16 &             2.75 &           34.74 &         54.91 &         1.50  &       289.47 &       194.30  &          23.07 &            7.63 &         281.83 &        72.10  &         0.23 \\
              4 &           1 &        5.93 &       0 &   - &    3600.67 &         0.03 &            0.13 &             3.79 &           26.10  &         63.64 &         1.30  &       264.90  &       143.40  &          46.30  &            7.23 &         257.67 &        75.20  &         0.67 \\
              4 &           2 &        5.93 &       0 &   - &    3600.96 &         0.04 &            0.10  &             3.05 &           23.88 &         66.15 &         1.30  &       235.77 &       180.00    &          10.17 &            6.13 &         229.63 &        45.60  &         0.80  \\
       \cmidrule(lr){1-19}
    &  &  &  \multicolumn{16}{c}{First-stage-only recourse}\\
    \cmidrule(lr){4-19}
     &  &  &  &  & \multicolumn{6}{c}{Time} & \multicolumn{7}{c}{Avg. total \#}\\
    \cmidrule(lr){6-11} \cmidrule(lr){12-19}
    $\chainCap$   & $\vBudget$   &  $\aBudget$  &      opt &      HSP (\%)    &   total     & Alg\ref{Alg:Heuristics} (\%) &   \ref{MasterExp} (\%) &   \ref{RecoP:Objexp} (\%) &   CG (\%)         &\ref{objMPOptQ} (\%)  & 1stS      &  2ndS         &   Alg\ref{Alg:Heuristics}-true  &   \ref{MasterExp} &   \ref{RecoP:Objexp} &      CGtrue &   \ref{objMPOptQ} &   DomS           \\
       \cmidrule(lr){1-19}
       $\lvert\vertexSet\rvert = 50$&             &             &             &           &          &              &                 &                  &                 &               &              &              &              &                &                 &                &              &          \\
          3 &           1 &        1.83 &      21 &   33.20   &    1093.26  &         2.94 &           20.68 &             0.72 &           28.72 &          7.82 &        57.53 &      5432.10  &      1815.10  &        3564.93 &           23.80  &        5408.30  &        52.07 &       179.57 \\
          3 &           2 &        1.83 &      25 &   26.06   &    706.25   &         3.64 &           19.86 &             0.34 &           24.42 &         37.49 &        19.10  &      3014.67 &      1356.07 &        1431.77 &           19.13 &        2995.53 &       226.83 &        69.43 \\
          3 &           1 &        3.03 &      18 &   30.69   &    1493.51  &         2.80  &           17.62 &             0.28 &           18.80  &         29.54 &        49.50  &      6231.33 &      1765.20  &        4273.47 &           82.97 &        6148.37 &       192.67 &       247.97 \\
          3 &           2 &        3.03 &      16 &   25.68   &    1881.49  &         2.11 &           11.40  &             0.30  &           11.99 &         50.50  &        87.37 &      8591.30  &      3024.77 &        5190.50  &           61.80  &        8529.50  &       376.03 &       635.13 \\
          4 &           1 &        1.83 &      22 &   36.05   &    1012.63  &         1.89 &           13.05 &             0.90  &           39.52 &          8.15 &        35.00    &      3360.53 &      1045.23 &        2268.30  &           18.03 &        3342.50  &        47.00    &        94.10  \\
          4 &           2 &        1.83 &      19 &   26.26   &    1508.41  &         2.52 &           12.60  &             0.49 &           33.78 &         25.87 &        47.30  &      6954.93 &      2394.40  &        4298.37 &           38.10  &        6916.83 &       262.17 &       868.63 \\
          4 &           1 &        3.03 &      16 &   30.74   &     1716.90  &         2.00    &           11.23 &             0.49 &           29.79 &         25.58 &        57.00    &      5617.77 &      1948.63 &        3365.53 &          102.60  &        5515.17 &       303.60  &       216.23 \\
          4 &           2 &        3.03 &      14 &   20.01   &    2124.06  &         1.71 &           10.43 &             0.13 &           18.28 &         51.53 &        45.77 &      6168.70  &      2230.47 &        3574.93 &           36.33 &        6132.37 &       363.30  &       224.50  \\
       \cmidrule(lr){1-19}
       $\lvert\vertexSet \rvert =  100$&             &             &         &           &          &              &                 &                  &                 &               &              &              &              &                &                 &                &              &          \\
          3 &           1 &        3.30  &       7 &   50.42 &    3054.67 &         0.22 &            1.11 &             0.28 &           21.42 &         68.69 &         7.40  &      1404.37 &       623.10  &         468.20  &           58.53 &        1345.83 &       313.07 &        31.50  \\
          3 &           2 &        3.30  &       3 &   28.67 &    3298.74 &         0.09 &            0.14 &             0.15 &            6.56 &         89.89 &         2.20  &       432.33 &       278.57 &          50.53 &            7.93 &         424.40  &       103.23 &         7.40  \\
          3 &           1 &        5.93 &       1 &   12.50 &    3510.17 &         0.07 &            0.14 &             0.11 &            6.15 &         89.18 &         2.83 &       540.10  &       325.20  &          97.10  &            7.67 &         532.43 &       117.80  &        20.47 \\
          3 &           2 &        5.93 &       1 &   12.50  &    3543.01 &         0.05 &            0.05 &             0.08 &            3.80  &         90.79 &         2.13 &       399.13 &       265.67 &          53.37 &            7.87 &         391.27 &        80.10  &        21.57 \\
          4 &           1 &        3.30  &       9 &   49.75 &    2974.72 &         0.13 &            0.56 &             1.93 &           51.37 &         37.58 &         3.90  &       782.93 &       334.43 &         240.17 &           25.10  &         757.83 &       208.33 &         5.97 \\
          4 &           2 &        3.30  &       4 &   52.77 &    3358.23 &         0.04 &            0.09 &             0.95 &           25.03 &         69.53 &         1.83 &       415.17 &       225.97 &          46.23 &           10.10  &         405.07 &       142.97 &        14.07 \\
          4 &           1 &        5.93 &       1 &   12.50 &    3552.61 &         0.04 &            0.13 &             1.79 &           21.07 &         71.97 &         1.80  &       357.90  &       177.60  &          89.83 &            9.13 &         348.77 &        90.47 &         5.23 \\
          4 &           2 &        5.93 &       1 &   25.00 &    3571.12 &         0.04 &            0.07 &             1.08 &           20.54 &         71.86 &         1.83 &       357.37 &       239.67 &          32.73 &            7.10  &         350.27 &        84.97 &         8.30  \\
    \bottomrule
 \end{tabular}
    \end{adjustbox}
\end{table}

\end{document}